\documentclass[english]{amsart}
  \title[Embedding of Toeplitz operators into $C_0$-semigroups]{Embedding of Toeplitz operators with smooth symbols into strongly continuous semigroups}

\date{\today}

\author[Fricain]{Emmanuel Fricain}
 \address{Univ. Lille, CNRS, UMR 8524 - Laboratoire Paul Painlevé, F-59000 Lille, France}
 \email{emmanuel.fricain@univ-lille.fr}

\author[Grivaux]{Sophie Grivaux}
\address{Univ. Lille, CNRS, UMR 8524 - Laboratoire Paul Painlevé, F-59000 Lille, France}
\email{sophie.grivaux@univ-lille.fr}

\author[Ostermann]{Ma\"eva Ostermann}
\address{Univ. Lille, CNRS, UMR 8524 - Laboratoire Paul Painlevé, F-59000 Lille, France}
\email{maeva.ostermann@univ-lille.fr}

\author[Yakubovich]{Dmitry Yakubovich}
\address{Departamento de Matem\'aticas, Universidad Aut\'onoma de Madrid, Canto-Blanco, 28049 Madrid, Spain}
\email{dmitry.yakubovich@uam.es}



\usepackage[utf8]{inputenc}
\usepackage{amssymb,amsmath,amsthm,amsrefs,newlfont,enumerate,array,graphicx,tikz,tkz-tab,xcolor,bbold,float,subcaption}
\usepackage{thmtools,thm-restate}
\usepackage[colorlinks, citecolor=blue, linkcolor=red, pdfstartview=FitB]{hyperref}
\usepackage[shortlabels]{enumitem}
\usepackage[noabbrev]{cleveref} 
\makeatletter
\newcommand{\customlabel}[2]{\def\@currentlabel{#2}\label{#1}}
\makeatother

\newtheorem{theorem}{Theorem}[section]
\newtheorem*{theorem*}{Theorem}
\newtheorem{lemma}[theorem]{Lemma}
\newtheorem{corollary}[theorem]{Corollary}
\newtheorem{fact}[theorem]{Fact}
\newtheorem{proposition}[theorem]{Proposition}

\theoremstyle{definition}

\newtheorem{remark}[theorem]{Remark}
\newtheorem{example}[theorem]{Example}

\newtheorem{claim}[theorem]{Claim}

\newcommand{\sic}{\sigma_c}

\newcommand{\rhoc}{\rho_c}

\newcommand{\sm}{\setminus}

\newcommand{\ga}{\gamma}

\newcommand{\pt}{\partial}

\newcommand{\intsp}{\mathrm{int}\left(\sigma(T_F)\right)}
\newcommand{\ps}[1]{\left<#1\right>}
\let\Im\relax\DeclareMathOperator{\Im}{Im}
\let\Re\relax\DeclareMathOperator{\Re}{Re}
\let\Ran\relax\DeclareMathOperator{\Ran}{Ran}
\newcommand{\bb}{\mathbb}

\DeclareMathOperator{\spa}{span}
\DeclareMathOperator{\card}{card}
\DeclareMathOperator{\codim}{codim}
\DeclareMathOperator{\conv}{conv}
\DeclareMathOperator{\dist}{dist}
\DeclareMathOperator{\w}{wind}
\newcommand{\dual}[2]{\ps{#1,#2}_{p,q}}
\def\D{\ensuremath{\mathbb D}}

\def\T{\ensuremath{\mathbb T}}

\def\Z{\ensuremath{\mathbb Z}}

\def\C{\ensuremath{\mathbb C}}

\begin{document}

\begin{abstract}
Using the model theory for Toeplitz operators with smooth symbols developed by the fourth author in the 80's, we study whether such  operators $T_{F}$ can be embedded  into a $C_{0}$-semigroup
of operators on the Hardy space $H^p$ of the open unit disk, $1<p<\infty$. We show that it is the case as soon as $0$ belongs to the unbounded connected component of $\mathbb{C}$ minus the interior of the spectrum of $T_{F}$.
We provide several conditions on the symbol $F$, both geometric and analytic in nature, ensuring that this sufficient condition is also necessary. 
For a certain class of symbols, where the curve $F(\T)$ is a ``figure eight in a loop" such that $\mathbb{C}\setminus\sigma(T_F)$ has a bounded connected component, we obtain a complete characterization of the embeddability of $T_F$ into a $C_0$-semigroup.
In the last part of the paper, we discuss the embeddability of $T_F$ when the symbol $F$ is not necessarily smooth, using connections with the numerical range and the functional calculus for bounded sectorial operators.
\end{abstract}
\maketitle
\section{Introduction}\label{Section 1}
In this paper, we investigate whether Toeplitz operators with smooth symbols acting on one of the Hardy spaces $H^{p}$ of the open unit disk, $1<p<\infty$, can be embedded into a $C_{0}$-semigroup of bounded operators on $H^{p}$. Recall that a family $(T_{t})_{t>0}$ of bounded operators acting on a Banach space $X$ is called a $C_{0}$\emph{-semigroup} (or a \emph{strongly
continuous semigroup}) if the following two properties are satisfied:

\begin{enumerate}
 \item [(i)] $T_{t}T_{s}=T_{t+s}$ for every $t,s>0$ (semigroup property);\par\smallskip
 
 \item [(ii)] $\|T_{t}x-x\|\longrightarrow 0$ as $t\longrightarrow 0^+$ for every $x\in X$. In other words, $T_{t}$ converges to the identity operator  on $X$ as $t$ tends to $0^+$
 for the pointwise topology on $X$, also called the Strong Operator Topology.
\end{enumerate}

An operator $T$ on $X$ is said to be \emph{embeddable into a $C_{0}$-semigroup} or, shortly, \emph{embeddable} if there exists a $C_{0}$-semigroup $(T_{t})_{t>0}$ of operators on $X$ such that
$T_{1}=T$. In other words, embeddable operators are those which appear as elements of a $C_{0}$-semigroup.
\par\smallskip
The question whether a given operator can be embedded into a $C_{0}$-semigroup is a difficult one, which was proposed by T.~Eisner in \cite{Eisner-E1} (see also \cite{Eisner2010}*{Ch V, Sec. 1}). It originates from an analogous question in ergodic theory, asking for conditions under which an ergodic measure-preserving transformation of
a probability space can be embedded into a flow. It was proved by de~la~Rue and Lazaro in \cite{RL} that generically, an automorphism of $[0,1]$ endowed with the Lebesgue measure
can be embedded into a flow of automorphisms of $[0,1]$. Eisner obtained several results pertaining to the question of the embeddability of operators. She proved for instance that
if an operator $T\in\mathcal{B}(X)$ is embeddable, then the dimension of its kernel $\ker T$, as well as the codimension of its range $\textrm{Ran}(T)$ is either 0 or infinite
\cite{Eisner2010}*{Ch IV, Th. 1.7}. When $T$ is an isometry (or a co-isometry), these conditions turn out to be sufficient \cite{Eisner2010}*{Ch V, Th. 1.19}. The embeddability of
various classes of operators is considered in \cite{Eisner2010}. It is proved in \cite{Eisner-WOT} that a typical contraction (in the Baire Category sense) on a complex separable Hilbert space $H$ in the Weak Operator
Topology is unitary, hence embeddable. As a consequence of their study of typical properties of contractions on $H$ in the Strong Operator Topology, Eisner and Matrai 
obtained in \cite{Eisner-Matrai} that such a typical contraction for SOT is embeddable too.

\par\smallskip
The present work is a contribution to the study of the embeddability problem for the particularly important class of Toeplitz operators on the Hardy spaces  $H^{p}$, $1<p<\infty$: given a function 
$F\in L^{\infty}(\mathbb{T})$, where $\mathbb{T}$ denotes the unit circle, the Toeplitz operator $T_{F}$ with symbol $F$ is defined on $H^{p}$ by
\[
T_{F}f:=P_{+}(F{f}),\qquad f\in H^{p},
\]
where $P_{+}$ is the Riesz projection from $L^{p}(\mathbb{T})$ onto $H^p$. Such an operator is well-defined and bounded on $H^{p}$. Very recently, Chalendar  and Lebreton  proved that if $\varphi$ is a non-constant inner function, the Toeplitz operator $T_\varphi$ is embeddable into a $C_{0}$-semigroup on $H^2$ if and only if $\varphi$ is not a finite Blaschke product \cite{Chalendar-Lebreton}. Note that this result is based on the criterion given by Eisner for embeddability of isometries. Moreover, they also showed (although it is not stated formally in their paper) that for a general analytic symbol $\varphi$, $T_\varphi$ is embeddable into a $C_{0}$-semigroup of Toeplitz operators if and only if $\varphi$ does not vanish on $\mathbb D$.
Our paper goes into a totally new direction and allows us to treat Toeplitz operators which are very far from being isometric, which makes the problem of embeddability quite difficult because we cannot use Eisner's simple criterion for isometries. The other difficulty is that the powers of a Toeplitz operator are in general very difficult to compute and in particular are almost never Toeplitz operators. Thus, even if a Toeplitz operator is embeddable into a $C_0$-semigroup, the members of the semigroup are rarely Toeplitz operators. Let us now explain in more details the content of our work.

\par\smallskip
In all sections of the paper except the last one, we consider symbols $F$ which are smooth, i.e. of class $C^{1+\varepsilon}$ on $\mathbb{T}$ for some $\varepsilon>\max(1/p,1/q)$, where $q$ is the conjugate exponent of $p$. 
Under some additional assumptions on $F$ (the main one being that $\w_{F}
(\lambda)\le 0$ for every $\lambda \in \mathbb{C}\setminus F(\mathbb{T})$, where $\w_{F}(\lambda )$ is the winding number of the curve $F(\mathbb{T})$ around the point $\lambda $), the fourth author developed in \cite{Yakubovich1989}, \cite{Yakubovich1991} and \cite{Yakubovich1996} (see also \cite{Yakubovich1993}) a model theory for such Toeplitz operators, showing that the adjoint $T^{*}_{F}$ of $T_{F}$ is isomorphic to the shift operator (the multiplication operator $M_\lambda$ by the independent variable $\lambda$) on a certain space of holomorphic functions $E_{F}^{q}$. This model space $E_{F}^{q}$ is defined as a direct sum of certain Smirnov spaces canonically associated to $F$, with some additional boundary conditions (see Section~\ref{Section 3} for the precise definition of $E_{F}^{q}$). This model was fruitfully exploited to investigate various properties of Toeplitz operators such as the description of invariant subspaces, properties of the commutant \cite{Yakubovich1989} and
cyclicity \cite{Yakubovich1991}, \cite{Yakubovich1996}, \cite{Yakubovich1993}, and then further by Fricain-Grivaux-Ostermann in \cite{FricainGrivauxOstermann_preprint} to explore properties of Toeplitz operators connected to linear dynamics, like hypercyclicity,
chaos, etc.
\par\smallskip
In this paper, we apply this model to the study of the embeddability problem for this class of Toeplitz operators with smooth symbols. We first describe the multiplier algebra of the model space
$E_{F}^{q}$, and use it to provide a necessary and sufficient condition for the model operator $M_{\lambda }$ on $E_{F}^{q}$ to be embeddable into a $C_{0}$-semigroup of
multiplication operators on $E^{q}_{F}$.

\begin{restatable}{Theorem}{TheoEmbModelMult}\label{Th:EmbModelMult}
The multiplication operator by the independent variable $M_\lambda$ acting on $E_F^q$ is embeddable into a $C_0$-semigroup of multiplication operators on $E_F^q$, i.e. operators of the form $M_g$, $g\in H^\infty(\intsp )$, if and only if $0$ belongs to the unbounded component of $\mathbb C\setminus \intsp $.
\end{restatable}
Here $\intsp $ denotes the interior of the spectrum $\sigma (T_{F})$ of $T_{F}$, and the algebra of analytic and bounded functions on $\intsp $ is denoted by $H^\infty(\intsp )$. Since $F$ is continuous on $\mathbb{T}$, $\sigma (T_{F})$ admits the following geometric
description:
\[
\sigma (T_{F})=\{\lambda \in \mathbb{C}\setminus F(\mathbb{T})\;;\;\w_{F}(\lambda )\neq 0\}\,\cup\,F(\mathbb{T})
\]
and the interior of $\sigma (T_{F})$ is thus easy to visualize.
\par\medskip
In order to state properly a first sufficient condition for $T_{F}$ to be embeddable, we need to give precisely the assumptions under which the model theory of \cite{Yakubovich1991} holds;
see also \cite{FricainGrivauxOstermann_preprint}*{Appendix B} for the $H^{p}$ version. Let $p>1$, and let $q$ be its conjugate exponent (i.e. $1/p+1/q=1$). The dual space of $H^p$ will be canonically identified to the space $H^q$ via the following duality bracket:
\begin{equation}
\ps{x,y}_{p,q}~=~\frac{1}{2\pi}\int_0^{2\pi}x(e^{i\theta})y(e^{-i\theta})\,\mathrm{d}\theta \quad x\in H^p, y\in H^q.
\end{equation}
This duality bracket is linear on both sides. We keep this somewhat unusual convention throughout the whole paper even in the case where $p=2$, except in \Cref{numerical-range} (which deals with the numerical range), where we get back to the usual definition of the duality bracket in the Hilbertian setting.
\par\smallskip
Consider the following three conditions on the symbol $F$:
\begin{enumerate}[(H1)]
    \item the function $F$ belongs to the class $C^{1+\varepsilon}(\mathbb T)$ for some $\varepsilon>\max(1/p,1/q)$, and the derivative $F'$ of $F$ does not vanish on $\mathbb T$ (here  $C^{1+\varepsilon}(\mathbb T)$ denotes the set
    of functions of class $C^{1}$ on $\mathbb{T}$ whose derivative is $\varepsilon$-H\"olderian);
    \par\smallskip
   \item the curve $F(\mathbb T)$ self-intersects a finite number of times, i.e. 
   the unit circle $\mathbb T$ can be partitioned into a finite number of closed arcs $\alpha_1,\dots,\alpha_m$ such that
     \begin{enumerate}[(a)]
       \item $F$ is injective on the interior of each arc $\alpha_j,~1\le j \le m$;
        \item for every $i\neq j,~1\le i,j\le m$, the sets $F(\alpha_j)$ and $F(\alpha_j)$ have disjoint interiors;
    \end{enumerate}
    \par\smallskip
    \item for every $\lambda\in\mathbb C\setminus F(\mathbb T)$, $\w_F(\lambda)\le0$, where $\w_F(\lambda)$ denotes the winding number of the curve $F(\mathbb T)$ around $\lambda$.
 \end{enumerate}
     \par\smallskip

It is an easy observation that given any bounded operator $T$ acting on a complex Banach space $X$, $T$ is embeddable into a $C_0$-semigroup as soon as $0$ belongs to the unbounded component of $\mathbb{C}\setminus \sigma(T) $ (see \Cref{Fact:RemOnPos0}).
Under the three conditions \ref{H1}, \ref{H2}, and \ref{H3}, a much finer property holds: $T^{*}_{F}\in\mathcal{B}(H^{q})$ is isomorphic (via an isomorphism $U:H^q\to E_F^q$) to $M_{\lambda }$ acting on the model space  $E_{F}^{q}$; one can then deduce from
Theorem \ref{Th:EmbModelMult} that $T^{*}_{F}$ is embeddable into a $C_{0}$-semigroup as soon as $0$ belongs to the unbounded component of 
$\mathbb{C}\setminus\intsp $. Now, it is a classical fact (see
for instance \cite{EngelNagel2006}*{Sec. I 1.13}) that whenever $(T_{t})_{t>0}$ is a $C_{0}$-semigroup of operators acting on a reflexive Banach space $X$, the adjoint semigroup $(T^{*}_{t})_{t>0}$ is a
$C_{0}$-semigroup of operators acting on $X^{*}$. In our setting, it follows that when $F$ satisfies \ref{H1}, \ref{H2}, and \ref{H3}, and $0$ belongs to the unbounded component of $\mathbb{C}
\setminus \intsp $, then $T_{F}\in \mathcal{B}(H^{p})$ is embeddable.
\par\medskip
An important observation is that the assumption \ref{H3} can in fact be replaced by the following hypothesis \ref{H3bis}, which requires that the winding number of $F$ has a constant sign on $\mathbb{C}\setminus F(\mathbb{T})$:
\par\smallskip
\begin{enumerate}[(H3bis)]
 \item \label{H3bis} either $\w_{F}(\lambda )\ge 0$ for every $\lambda \in\mathbb{C}\setminus F(\mathbb{T})$, or $\w_{F}(\lambda )\le 0$ for every $\lambda \in\mathbb{C}\setminus F(\mathbb{T})$.
\end{enumerate} 
\par\smallskip
If $\w_{F}(\lambda )\ge 0$ for every $\lambda \in\mathbb{C}\setminus F(\mathbb{T})$, then setting $f(z)=F(1/z)$, $z\in\mathbb{T}$, we have $T_{f}=T_{F}^*\in\mathcal{B}(H^q) $, and $\w_f(\lambda)=-\w_{F}(\lambda)$ for every $\lambda\in \mathbb{C}\setminus F(\mathbb{T})$. Thus $f$ satisfies \ref{H1}, \ref{H2} and \ref{H3}. Hence $T_{f}^{*}=T_{F}\in\mathcal{B}(H^{p})$ is isomorphic to the multiplication operator by $\lambda $ on the model  space
$E_{f}^{p}$, and since $\sigma (T_{f})=\sigma (T_{F})$, it follows that if $0$ belongs to the unbounded component of $\mathbb{C}\setminus \intsp$, then $T_{F}$ is
embeddable. As a consequence of \Cref{Th:EmbModelMult}, we obtain:

\begin{restatable}{Theorem}{CSPlongement}\label{Th:CSforEmb}
    Let $p>1$. Suppose that $F$ satisfies the assumptions \emph{\ref{H1}}, \emph{\ref{H2}} and \emph{\ref{H3bis}}. If 0 belongs to the unbounded component of $\mathbb{C}\setminus\intsp $, then $T_{F}\in\mathcal{B}(H^{p})$ is embeddable into a $C_{0}$-semigroup.
\end{restatable}

In the sequel of the paper, we investigate the converse of Theorem \ref{Th:CSforEmb}, i.e. the question of whether the embeddability of $T_{F}$ into a $C_{0}$-semigroup implies that $0$
belongs to the unbounded component of $\mathbb{C}\setminus \intsp $. We obtain some conditions under which this converse is true. For instance, we show:
\begin{restatable}{Theorem}{CondPourEquiv}\label{th: premiere-CNS}
    Let $p>1$, and let $F$ satisfy \emph{\ref{H1}}, \emph{\ref{H2}} and \emph{\ref{H3bis}}. Suppose that $\mathbb{C}\setminus \intsp $ is connected, and that $0$ is not an
intersection point of the curve $F(\mathbb{T})$. Then the following assertions are equivalent: 
\begin{enumerate}
\item [(1)] $T_{F}$ is embeddable into a $C_0$-semigroup of bounded operators on $H^p$;\par\smallskip
\item [(2)] $0$ belongs to $\mathbb{C}\setminus\intsp $.
\end{enumerate}
\end{restatable}

We obtain several results in this vein. Observe that if an operator $T\in \mathcal B(X)$ is embedded in a $C_0$-semigroup $(T_t)_{t>0}$, then necessarily the operators $T_t$ belong to the commutant of $T$. Hence, this question of embeddability of an operator $T$ is linked in a natural way to the question of describing its commutant. Thus we also study the commutant of the operator $M_\lambda$ acting on the model space $E^q_F$, and observe that, very surprisingly, it may or not consist entirely of multipliers. More precisely, we exhibit a curve admitting two different parametrizations $F_1$ and $F_2$ such that the commutant of $M_\lambda$ acting on  $E^q_{F_1}$ is made of multipliers, while the commutant of $M_\lambda$ acting on  $E^q_{F_2}$ is not (\Cref{exemples-commutant}). Consequently, we provide conditions of an analytic nature implying that the commutant of $T_{F}$ (when seen on $E_{F}^{q}$) consists of multipliers. Under such conditions, $T_{F}$ is embeddable if and only if $0$ belongs to the unbounded connected component of $\mathbb{C}\setminus\intsp $ (see Sections \ref{Section 5} and \ref{Section 6} for details).
 \par\medskip
\Cref{Section 7}, which is the most technical part of the paper, is an attempt to understand how the embeddability of $T_F$ could be characterized for general symbols $F$, under minimal assumptions. In informal terms, here is the general form of the results we obtain: suppose that $F$ satisfies \ref{H1}, \ref{H2} and \ref{H3}.
Given a connected component $\Omega$ of $\C\setminus F(\T)$ with $|\w_F(\Omega)|=2$, we look at the inverse function $\zeta=1/F^{-1}$ of $F$ (which is well-defined on the curve $F(\T)$ 
minus its points of self-intersection) on boundary arcs of $\partial\Omega$. If, whatever the choice of such an arc $\gamma$, the restriction of the function $\zeta$ to $\gamma$ does not coincide a.e. with the boundary limit of a meromorphic function in the Nevanlinna class of $\Omega$, then the embeddability of $T_F$ forces $0$ to belong to the unbounded component of $\C\setminus\intsp$. These are roughly the contents of \Cref{un-th-supplementaire}. Then we study what happens when the function $\zeta$ is a.e. a boundary limit of a meromorphic function in the Nevanlinna class of $\Omega$ on suitable arcs $\gamma\subseteq\partial\Omega$, first on an example (\Cref{l'exemple!!}), and then in a more general situation where the curve $F(\T)$ looks like a "figure-eight inside a loop" -- see \Cref{Fig5}. In this case, where $\C\setminus F(\T)$ has four connected components (one of winding number $-1$, one of winding number $-2$ and two of winding number $0$ - a bounded one and an unbounded one), we completely describe the cases where $T_F$ is embeddable. The full answer to the embeddability problem in this case is given by the following result:

\begin{restatable}{Theorem}{LaCaracterisation}\label{la-caracterisation!!}
Let $1<p<\infty$ and let $F$ satisfy \ref{H1}. Suppose that $F(\mathbb T)$ is given by \Cref{Fig5} and that
$0\notin\mathcal O$. Then $T_F$ is embeddable into a $C_0$-semigroup if and only if 
one of the following two conditions hold:
\begin{itemize}
    \item[(1)] $0$ belongs to the unbounded component of $\mathbb C\setminus \intsp$;
    \item[(2)] $0$ belongs to the bounded component of $\mathbb C\setminus\intsp$ and the following three conditions hold:
    \begin{enumerate}[(i)]
\item $\zeta_{|\gamma_1}$ (resp. $\zeta_{|\gamma_2}$) coincides a.e. on $\gamma_1$ (resp. on $\gamma_2$) with the non-tangential limit of  a meromorphic functions $\zeta_1$ (resp. $\zeta_2$) on $\Omega_2$;
\item the measure $\mu$ on $\intsp$ defined by
\[
\mathrm d\mu(\lambda)=\frac{\mathbf{1}_{\partial\Omega_2}(\lambda)}{|\zeta_1(\lambda)-\zeta_2(\lambda)|^q}|\mathrm d\lambda|
\]
is a Carleson measure for $E^q(\intsp)$;
\item the maps $$Z_1:w\mapsto \frac1{\zeta_1-\zeta_2}(C_\zeta w-\zeta_1 w)\quad Z_2:w\mapsto \frac1{\zeta_1-\zeta_2}(C_\zeta w-\zeta_2 w)$$ $$Z_3:w\mapsto \frac{\zeta_2}{\zeta_1-\zeta_2}(C_\zeta w-\zeta_1 w)\quad Z_4:w\mapsto \frac{\zeta_1}{\zeta_1-\zeta_2}(C_\zeta w-\zeta_2 w)$$ define bounded operators from $E^q(\Omega_2)$ into itself.
 \end{enumerate}
\end{itemize}
\end{restatable}

\par\medskip
In the last section of the paper, we present some results concerning the embeddability of $T_F$ when the symbol $F$ is not necessarily smooth. The methods here are different from those employed in the rest of the paper. Here the functional calculus from \cite{Yakubovich1991} does not apply anymore, and has to be replaced by a suitable functional calculus for bounded sectorial operators. See \Cref{Subsec:defsectorial} for definitions and some consequences of this functional calculus. Whenever we consider operators acting on a  Hilbert space, we will also study the link between embeddability of $T$ and properties of the numerical range of $T$, which is defined by $W(T)=\{\ps{Tx,x}\,;\,x\in H~\text{and}~\|x\|=1\}$.
\par\smallskip
For instance, we obtain the following sufficient condition for embeddability of Toeplitz operators on $H^2$:

\begin{restatable}{Theorem}{ExempleImageNumerique}\label{Theo4}
Let $F\in L^\infty(\mathbb T)$. Suppose that $0$ does not belong to the interior $\emph{int}({W(T_F)})$ of the numerical range $W(T_F)$ of $T_F$. Then $T_F$ is embeddable into a $C_0$-semigroup of operators on $H^2$.
\end{restatable}

As a consequence, it follows that $T_F$ is embeddable into a $C_0$-semigroup of operators on $H^2$ as soon as $\Re F\ge 0$ a.e. on $\T$.
\par\smallskip
We also explore some consequences of results of Peller (valid only for $p=2$), who gave in \cite{Peller1986} conditions implying that the Kreiss constant of the spectrum of a Toeplitz operator with a symbol belonging to certain algebras of functions on $\T$ is finite. In particular, we show the following:

\begin{restatable}{Theorem}{ConsequencePeller}\label{Theo14}
Let $\mathcal X$ be the Wiener algebra (i.e. the algebra of functions on the circle $\mathbb{T}$ with absolutely convergent Fourier series), or the algebra of Dini-continuous functions on $\T$.
Let $F \in \mathcal X$. Suppose that  
there exists an open disk $D$, contained in the 
unbounded component of $\mathbb C\setminus \sigma(T_F)$, such that $0\in \partial D$. 
Then $T_F$ is embeddable into a $C_0$-semigroup of operators on $H^2$.
\end{restatable}

\begin{figure}[ht]
\begin{subfigure}{.48\linewidth}\centering
    \begin{tikzpicture}[scale=1]
        \draw(0,0)node{\includegraphics[scale=1]{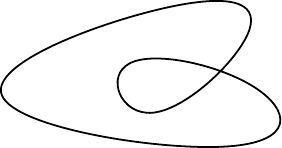}};
        \draw[red](1.35,0.05)node{$+$}(1.4,0.1)node[right]{$0$};
    \end{tikzpicture}
    \caption{}
    \label{SubFig:1a}
\end{subfigure}
~
\begin{subfigure}{.48\linewidth}\centering
    \begin{tikzpicture}[scale=1]
        \draw(0,0)node{\includegraphics[scale=1]{FigIntro.pdf}};
        \draw[red](-1.45,0.5)node{$+$}(-1.5,0.55)node[left]{$0$};
    \end{tikzpicture}
    \caption{}
    \label{SubFig:1b}
\end{subfigure}
    \caption{}
\end{figure}
Note that if $F(\mathbb T)$ is a general curve of the kind which is described by Figure \textsc{\ref{SubFig:1a}},  neither \Cref{Theo4} nor \Cref{Theo14} can be applied to show that $T_F$ is embeddable into a $C_0$-semigroup on $H^2$. However, if $F$ satisfies furthermore  assumption \ref{H1}, then it follows from \Cref{Th:CSforEmb} that $T_F$ is embeddable into a $C_0$-semigroup on $H^2$ (or even on $H^p$). On the other hand, when we are in the situation given by Figure \textsc{\ref{SubFig:1b}},  either \Cref{Theo4} or \Cref{Theo14} yields that $T_F$ is embeddable into a $C_0$-semigroup, without any smoothness assumption on the symbol $F$. 

\par\medskip
The paper is organized as follows: we recall in Sections \ref{Section 2} and \ref{Section 3} some results regarding $C_{0}$-semigroups, Toeplitz operators, and Yakubovich's model space which will be necessary for our
study. In Section \ref{Section:Multipliers}, we first prove that the multiplier algebra of the model space $E_{F}^{q}$ coincides with the set of multiplication operators by functions
$g\in H^{\infty}(\intsp )$. We then prove Theorem \ref{Th:EmbModelMult} 
as well as the sufficient condition for embeddability  given by Theorem \ref{Th:CSforEmb}. 
We study in Sections \ref{Section 5} and \ref{Section 6} necessary conditions for embeddability, and obtain characterizations of embeddability under conditions of different kinds (\Cref{th: premiere-CNS}, \Cref{th:deuxiemeCNS}, \Cref{un-th-supplementaire}, \Cref{th:commutant}). Theorem \ref{th: premiere-CNS} in particular is proved in Section \ref{Section 5}, 
and our study of the commutant of the operator $M_\lambda$ on the model space is carried out in \Cref{Section 6.2,Subsection6.3} (see \Cref{exemples-commutant} and \Cref{th:commutant}).
Section \ref{Section 7} is devoted to the study of the general situation where the curve $F(\T)$ looks like a "figure-eight inside a loop" -- see \Cref{Fig5} -- and to the characterization of the embeddability of $T_F$ in this case (\Cref{la-caracterisation!!}).
The final Section \ref{Section: sectorial} contains further results on the embeddability of $T_F$ when the symbol $F$ is not supposed to be smooth.  After some reminders concerning the functional calculus for sectorial operators, the Kreiss constant and the links between the finiteness of this constant and sectoriality, we prove \Cref{Theo4} as well as (a more general version of) \Cref{Theo14}, which is \Cref{Coro12}.
\par\medskip
\textbf{Thanks:}
We warmly thank Yuri Tomilov for stimulating conversations, and for pointing out to us several useful references.

\section{First reminders and notations}\label{Section 2}
\subsection{Reminders on $C_0$-semigroups}\label{Subsection:Rappel-Co-semigroupes}
In this short subsection, we recall very briefly some notation and results on semigroups of operators which will be used in this paper. The letter $X$ will always denote a complex separable infinite-dimensional Banach space, and $\mathcal B(X)$ the space of linear and continuous operators from $X$ into itself. Let $(T_t)_{t>0}\subset \mathcal B(X)$. The definition of what it means for $(T_t)_{t>0}$ to be a strongly continuous semigroup (or $C_0$-semigroup) was recalled in the Introduction.
Recall also that an operator $T\in \mathcal B(X)$ is embeddable into a $C_0$-semigroup if there exists a $C_0$-semigroup $(T_t)_{t>0}$ such that $T_1=T$.
We have the following easy result:

\begin{proposition}[\cite{EngelNagel2006}*{Sec. I.1.13}]\label{Prop:Reflexif}
Let $X$ be a reflexive Banach space. If $(T_t)_{t>0}$ is a $C_0$-semigroup on $X$, then $(T_t^*)_{t>0}$ is a $C_0$-semigroup on $X^*$.
\end{proposition}

A direct consequence is that if $X$ is a reflexive Banach space and if $T\in \mathcal B(X)$  is embeddable into a $C_0$-semigroup on 
$X$ then $T^*$ is also embeddable into a $C_0$-semigroup on $X^*$.
\par\smallskip
We now recall an important necessary condition for an operator to be embeddable into a $C_0$-semigroup:

\begin{theorem}[\cite{Eisner2010}*{Th. V.1.7}]\label{Th:CN-plong-ker}
Let $X$ be a Banach space and $T\in \mathcal B(X)$. If $T$ is embeddable, then
\[\dim\ker(T)\in\{0,\infty\}\, \textrm{ and } \, \dim\ker(T^*)=\codim\emph{Ran}(T)\in\{0,\infty\}.\]
\end{theorem}

In particular, it follows from \Cref{Th:CN-plong-ker} that a non-bijective Fredholm operator is not embeddable.

\subsection{Reminders on Toeplitz operators}
For $1<p<+\infty$, we denote by $H^p=H^p(\mathbb D)$ the Hardy space of analytic functions $u$ on the open unit disk $\mathbb D$ such that 
\[\|u\|_{H^p}~:=~\sup_{0<r<1}M_p(u,r)<+\infty\quad\text{where}~M_p(u,r)=\left(\int_0^{2\pi}|u(re^{i\theta})|^p\frac{\mathrm d\theta}{2\pi}\right)^{1/p}.\]
A function $u$ belonging to  $H^p$ has non tangential boundary values  $u^*$ almost everywhere on $\mathbb T$. We will still denote this boundary value as $u$. It is well-known that $\|u\|_{H^p}=\|u\|_{L^p(\mathbb T)}$. The dual space of $H^p$ is canonically identified to the space $H^q$, where $q$ is the conjugate exponent of $p$; the duality is given by the formula
\begin{equation}\label{eq:duality-defn}
\ps{x,y}_{p,q}~=~\frac{1}{2\pi}\int_0^{2\pi}x(e^{i\theta})y(e^{-i\theta})\,\mathrm{d}\theta,
\end{equation}
where $x\in H^p$ and $y\in H^q$.
The duality bracket in (\ref{eq:duality-defn}) is linear on both sides, and we keep this convention even in the case where $p=2$; in particular, adjoints of operators on $H^2$  must be understood as Banach space adjoints, and not Hilbert space adjoints.
\par\smallskip
Let $P_+$ denote the Riesz projection from $L^p(\mathbb T)$ onto $H^p$ defined by 
\[(P_+u)(z)~=~\frac{1}{2\pi}\int_0^{2\pi}\frac{u(e^{i\theta})}{1-ze^{-i\theta}}\,d\theta
\quad\text{for }
z\in\mathbb D\text{,}\;
u\in L^p(\mathbb T).\]
Given $F\in L^\infty(\mathbb T)$, the Toeplitz operator $T_F$ with symbol $F$ is defined on $H^p$ by the following formula: $T_Fu=P_+(Fu),~u\in H^p$. It is a bounded operator on $H^p$. Recall that, if $F$ is a continuous function on $\mathbb T$ (which we write as $F\in C(\T)$), then $T_F-\lambda$ is a Fredholm operator on $H^p$ if and only if $\lambda\notin F(\mathbb T)$; when $\lambda\notin F(\mathbb T)$, the Fredholm index of $T_F-\lambda$ is equal to $-\w_F(\lambda)$, and we can describe the spectrum of $T_F$ as 
\[\sigma(T_F)~=~\{\lambda\in\mathbb C\setminus F(\mathbb T)\,;\,\w_F(\lambda)\neq0\}\cup F(\mathbb T),\]
where $\w_F(\lambda)$ is the winding number of the curve $F(\mathbb T)$ with respect to the point $\lambda$. Since, by the Coburn Theorem, a Toeplitz operator is either injective or has dense range,  we have for every $\lambda\in\mathbb C\setminus F(\mathbb T)$
\begin{equation}\label{dimension-noyau-Fredholm}
\dim\ker(T_F-\lambda)=\max(0,-\w_F(\lambda)).
\end{equation}
We refer the reader to \cite{BottcherSilbermann1990} for all basic facts on Toeplitz operators.
\par\smallskip
As a direct consequence on the embedding of Toeplitz operators, we have:

\begin{fact}\label{Fact:RemOnPos0}
Let $F\in C(\mathbb T)$ and $1<p<\infty$. Let $T_F\in\mathcal{B}(H^p)$ be the Toeplitz operator with symbol $F$.
\begin{enumerate}[(i)]
    \item If $0\in \sigma(T_F)\setminus F(\mathbb T)$ then $T_F$ is not embeddable.
    \item If $0$ belongs to the unbounded connected component of $\mathbb C\setminus F(\mathbb T)$, then $T_F$ is embeddable.
\end{enumerate}
\end{fact}

\begin{proof}
Assertion (i) follows immediately from the fact that if $0\in \sigma(T_F)\setminus F(\mathbb T)$, then $T_F$ is a non-bijective Fredholm operator, hence is not embeddable. As to assertion (ii), it suffices to observe that if $0$ belongs to the unbounded component of $\mathbb C\setminus F(\mathbb T)$, then $0$ belongs to the unbounded component of $\mathbb C\setminus \sigma(T_F)$ and thus there exists an analytic determination of the logarithm, denoted by $\log$, on a neighborhood of $\sigma(T_F)$. Hence $\log(T_F)$ is a well-defined and bounded operator on $H^p$ by the Dunford functional calculus and thus $T_F$ can be embedded into the $C_0$-semigroup $(T_F^t)_{t>0}$ where $T_F^t=e^{t\log(T_F)}$.
\end{proof}

\begin{remark}
Note that if $0$ belongs to the unbounded component of $\mathbb C\setminus F(\mathbb T)$, then the proof of Fact~\ref{Fact:RemOnPos0} shows that in fact, $T_F$ can be embedded into a semigroup $(T_t)_{t>0}$ which is even uniformly continuous, meaning that $\|T_t-I\|\to 0$ as $t\to 0^+$. In this case, the semigroup has a generator given by $\log(T_F)\in\mathcal B(H^p)$.
\end{remark}

We will see later on in the paper that, under some conditions on the symbol $F$, $T_F$ is embeddable as soon as $0$ belongs to the unbounded component of the complement of the interior of the spectrum $\sigma(T_F)$. Let us observe here that $T_F$ can be embeddable even if $0$ belongs to some bounded component of $\mathbb C\setminus \sigma(T_F)$. Indeed, let $F$ be a bounded analytic function which does not vanish on $\D$. Then, using the canonical decomposition of $F$ as a product of a singular inner function and an outer function, it is easy to see that $T_F\in\mathcal{B}(H^p)$ can be embedded into a $C_0$-semigroup 
of analytic Toeplitz operators on $H^p$ (see \cite{Chalendar-Lebreton}*{Th. 3.9 and Lem. 3.10}). 
\par\smallskip
In particular, we can give the following concrete  example of an embeddable Toeplitz operator such that $0$ belongs to a bounded component of the complement of its spectrum.

\begin{example}\label{example1:embeddable-bounded-compopent}
Let $\varphi(z)=z+2$, $z\in\D$. Using the principal determination of the logarithm on $\C\setminus (-\infty, 0]$, define, for any $s>0$,
$\varphi^s(z)=(z+2)^s$, $z\in\D$. As $\theta$ grows from $0$ to $2\pi$, the argument of $e^{i\theta}+2$ first grows from 0 to the value $a:=\arctan(1/2)$, then decreases from $a$ to $\pi-a$, and finally increases from $\pi-a$ to $2\pi$. It follows that whenever $s>s_0:=\pi/a$, the curve $\varphi^s(\T)$ looks like in \Cref{fig:ExExt} (the picture is not to scale, since $\max_{z\in\T}|\varphi^s(z)|=3^s$ is extremely large compared to  $\min_{z\in\T}|\varphi^s(z)|=1$):
\begin{figure}[ht]
    \includegraphics[scale=.35]{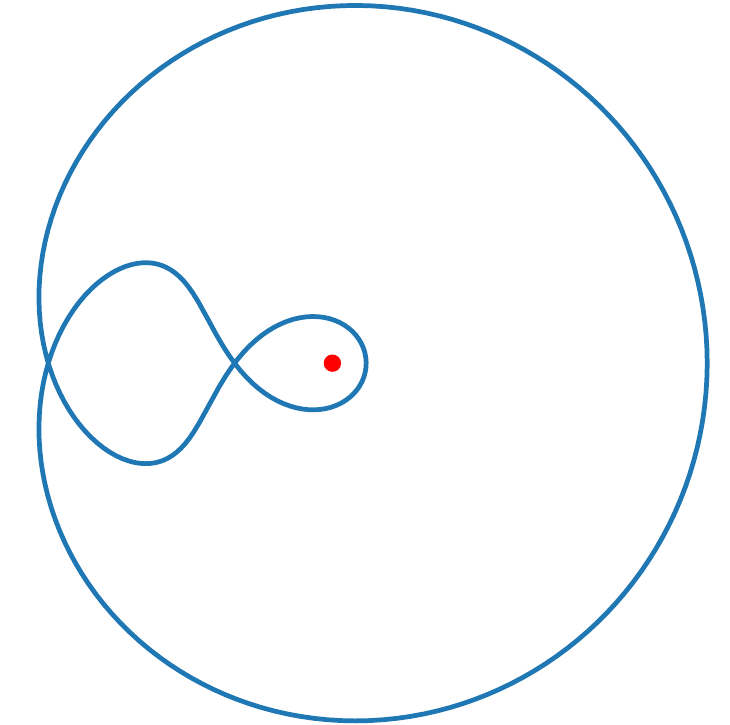}
    \caption{}
    \label{fig:ExExt}
\end{figure}

The point $0$ belongs to the bounded component of $\C\setminus\sigma(T_{\varphi^{s}})$, and nonetheless $T_{\varphi^s}$ is embeddable into a $C_0$-semigroup of analytic Toeplitz operators by the observation above.
\end{example}

\section{A model for Toeplitz operators with smooth symbols}\label{Section 3}

In this section, we recall the model  developed in the hilbertian case ($p=2$) by Yakubovich in \cite{Yakubovich1991}. See Appendix B in \cite{FricainGrivauxOstermann_preprint} for the details of the $H^p$ version of this model. 
Let $p>1$ and let $q$ be its conjugate exponent, i.e. $\frac{1}{p}+\frac{1}{q}=1$. 
In this section, we assume that the symbol $F$ satisfies the following three conditions:

\begin{enumerate}[(H1)]
    \item\label{H1}$F$ belongs to the class $C^{1+\varepsilon}(\mathbb T)$ for some $\varepsilon>\max(1/p,1/q)$, and its derivative $F'$ does not vanish on $\mathbb T$;
    \par\smallskip
   \item\label{H2} the curve $F(\mathbb T)$ self-intersects a finite number of times, i.e. 
   the unit circle $\mathbb T$ can be partitioned into a finite number of closed arcs $\alpha_1,\dots,\alpha_m$ such that
     \begin{enumerate}[(a)]
       \item $F$ is injective on the interior of each arc $\alpha_j,~1\le j \le m$;
        \item for every $i\neq j,~1\le i,j\le m$, the sets $F(\alpha_j)$ and $F(\alpha_j)$ have disjoint interiors;
    \end{enumerate}
    \par\smallskip
    \item\label{H3}for every $\lambda\in\mathbb C\setminus F(\mathbb T)$, $\w_F(\lambda)\le0$, where $\w_F(\lambda)$ denotes the winding number of the curve $F(\mathbb T)$ around $\lambda$.
 \end{enumerate}
Assumption \ref{H2} implies that the curve $F(\mathbb T)$ has only a finite number of points of self-intersection which we denote by $\cal O$. Since $F$ is a bijective map from $\mathbb T\setminus F^{-1}(\cal O)$ onto $F(\mathbb T)\setminus\cal O$, the map
$\zeta=1/F^{-1}$ is well-defined on $F(\mathbb T)\setminus\cal O$. It is also of class ${C}^1$ on each open arc contained in $F(\mathbb T)\setminus\cal O$.

\subsection{Eigenvectors}
It follows from assumption \ref{H3} and from the $H^p$ version of the Coburn Theorem that for every $\lambda\in\mathbb C\setminus F(\mathbb T)$, we have
\[\ker(T_F^*-\lambda)~=~\{0\}\quad\text{and}\quad\dim\,\ker(T_F-\lambda)~=~-\w_F(\lambda).\]
The fourth author of this paper provided in \cite{Yakubovich1991} an explicit expression of a spanning family of elements of the eigenspaces $\ker(T_F-\lambda),~\lambda\notin F(\mathbb T)$. Let $\lambda\in\mathbb C\setminus F(\mathbb T)$, and consider the function $\phi_\lambda$ defined on $\mathbb T$ by 
\[\phi_\lambda(\tau)~=~\tau^{-\w_F(\lambda)}(F(\tau)-\lambda)\quad\text{for }\tau\in \mathbb T.\]
Since $\phi_\lambda$ is of class $C^{1+\varepsilon}$ and does not vanish on $\mathbb T$, and since $\w_{\phi_\lambda}(0)=0$, one can define a logarithm $\log\phi_\lambda$ of $\phi_\lambda$ on $\mathbb T$ that is of class $C^{1+\varepsilon}$ on $\mathbb T$, and set
\begin{equation}\label{eq:vecteur-propre-Yaku33434ZD}
F_\lambda^+~=~\exp(P_+(\log \phi_\lambda)).
\end{equation}
The functions $F_\lambda^+$ and $1/F_\lambda^+$ both belong to the disk algebra $A(\mathbb D)$ (which is the space of holomorphic functions on $\mathbb{D}$ which admit a continuous extension to the closure $\overline{\mathbb{D}}$ of the unit disk, endowed with 
the supremum norm on $\overline{\mathbb{D}}$). For every connected component $\Omega$ of $\mathbb C\setminus F(\mathbb T)$ and  every $z\in \overline{\mathbb D}$, the map $\lambda\mapsto F_\lambda^+(z)$ is analytic on $\Omega$ and continuous on $\overline{\Omega}$. 
For $\lambda\in\sigma(T_F)\setminus F(\mathbb T)$, set 
\begin{equation}\label{eve}
  h_{\lambda,j}(z)~=~z^j\frac{F_\lambda^+(0)}{F_\lambda^+(z)}\quad\text{for every }z\in\mathbb D~\text{and every}~0\le j<|\w_F(\lambda)|.  
\end{equation}
These functions $h_{\lambda,j}$ belong to $A(\mathbb D)$, hence to $H^p$, and it can be checked that $(T_F-\lambda)h_{\lambda,j}=0$ for every $0\le j<|\w_F(\lambda)|$. So, we have
\begin{equation}\label{eq:eigenvector-space}
\ker(T_F-\lambda)~=~\spa\big[h_{\lambda,j}\,;\,0\le j<|\w_F(\lambda)|\big]
\end{equation}
for every $\lambda\in\sigma(T_F)\setminus F(\mathbb T)$.

\subsection{The model space}
In this section, we introduce the definition of the model space for Toeplitz operators with symbols satisfying \ref{H1}, \ref{H2} and \ref{H3}. The construction of this model space is based on the Smirnov spaces, whose main properties we now recall. Even if we do not mention it specifically in every statement of the paper, all domains under consideration will be assumed to have rectifiable boundary.
\subsubsection{Smirnov spaces}\label{Sub-section-Smirnov}
Let $\Omega$ be a bounded simply connected domain of $\mathbb C$ whose boundary $\Gamma $ admits a piecewise $C^1$ parametrization. Given $1<q<\infty$, an analytic function $f$ on $\Omega$ is said to belong to the \emph{Smirnov space} $E^q(\Omega)$ if there exists a sequence  of rectifiable Jordan curves $(C_n)_{n\geq 1}$ included in $\Omega$, tending to the boundary $\Gamma$ (in the sense that $C_n$ eventually surrounds each compact subdomain of $\Omega$), and such that 
\begin{equation}\label{eq-defn1-Ep}
\sup_{n\geq 1}\int_{C_n}|f(z)|^q\,|\mathrm{d}z|~<~\infty.
\end{equation}
Observe that $f$ belongs to  $E^q(\Omega)$ if and only if $(f\circ \varphi)\cdot \varphi'^{1/q}$ belongs to $H^q(\mathbb D)$ for some (equivalently, all) conformal map $\varphi$ from $\mathbb D$ onto $\Omega$. In particular, every function $f\in E^q(\Omega)$ admits a non-tangential limit almost everywhere on $\Gamma$, and this non-tangential limit belongs to $L^q(\Gamma)$; we still denote it by $f$ in order to simplify notation. Note that the non-tangential limit cannot vanish on a set of positive measure unless $f$ is identically $0$ (see for instance \cite{Duren1970}*{Th. 10.3}). 
\par\smallskip
In this paper, we shall also need the extension of the classes $E^q(D)$ to finitely connected bounded domains $D$ whose boundary $C$ consists of finitely many rectifiable Jordan curves. Recall that $D$ is said to be \emph{finitely connected} if its complement in the extended complex plane has finitely many connected components. An analytic function $f$ on $D$ is said to belong to the class $E^q(D)$ if there exists a sequence $(\Delta_n)_{n\geq 1}$ of domains with boundaries $(\Gamma_n)_{n\geq 1}$, each $\Gamma_n$ being a finite union of rectifiable Jordan curves, such that the sequence $(\Delta_n)_{n\geq 1}$ exhausts $D$, the lengths of the curves $\Gamma_n$ are uniformly bounded, and
\[
\sup_{n\geq 1}\int_{\Gamma_n}|f(z)|^q\,|\mathrm dz|<\infty.
\]
Moreover (since $q \geq 1$) the function $f$ can be uniquely recovered from its boundary values by means of the Cauchy integral taken over $C$.
\par\smallskip
If $\Omega$ is the disjoint union of finitely connected bounded domains $O_1,\dots,O_N$ of $\mathbb C$, an analytic function $f$ on $\Omega$ is said to belong to $E^q(\Omega)$ if for every $1\leq j\leq N$, the restriction $f_{|O_j}$ of $f$ to $O_j$ belongs to $E^q(O_j)$.  Then $E^q(\Omega)$ is a Banach space, equipped with the norm
\begin{equation}\label{eq:norm-Smirnov-fc}
\|f\|_{E^q(\Omega)}~=~\left(\sum_{j=1}^N\|f_{|O_j}\|_{E^q(O_j)}^q\right)^{1/q}=~\left(\sum_{j=1}^N\int_{\gamma_j}|f_{|O_j}(z)|^q\,|\mathrm{d}z| \right)^{1/q}
\end{equation}
where $\gamma_j=\partial O_j$, $1\le j\le N$.
\par\smallskip

\par\medskip
We finish this section with some reminders on the Cauchy transform. Let $1<q<\infty$, and let $\Gamma$ be a rectifiable Jordan curve in $\mathbb C$. Denote by $\Omega_0$ the interior of $\Gamma$, 
and by $\Omega_\infty$ its exterior. 
For any function $f\in L^q(\Gamma)$, the Cauchy transform of $f$ is the holomorphic function on
$\mathbb C\setminus\Gamma=\Omega_0\cup\Omega_\infty$ 
defined by 
\[
\mathcal C f(z)=\frac{1}{2i\pi}\int_\Gamma \frac{f(\zeta)}{\zeta-z}\,d\zeta,\quad z\in\Omega_0\cup\Omega_\infty.
\]
We recall a result of G. David \cite{David} concerning the boundedness of the Cauchy transform
from $L^q(\Gamma)$ into 
one the generalized Hardy spaces associated to $\Omega_0$ or $\Omega_\infty$. 
More precisely, denote by 
$\Pi^q(\Omega_0)$ the closure in $L^q(\Gamma)$ of the (analytic) polynomials on $\Omega_0$, and by $\Pi^q(\Omega_\infty)$ the closure in $L^q(\Gamma)$ of functions of the forms $P(\frac{1}{z-a})$, where $a$ is any fixed point in $\Omega_{\infty}$ and $P$ is any polynomial with zero constant term. 
If the curve $\Gamma$ is supposed moreover to be a \emph{Carleson curve} (i.e. if there exists a constant $C>0$ such that, for every $\mu\in\mathbb C$ and every $r>0$, the length of $\Gamma\cap D(\mu,r)$ is bounded by $Cr$), then it is proved in \cite[Thm. 3]{David} that the Cauchy transform $\mathcal C$ defines a bounded operator from $L^q(\Gamma)$ into $\Pi^q(\Omega_0)$,
as well as a bounded operator from $L^q(\Gamma)$ into $\Pi^q(\Omega_\infty)$. 
Moreover, since $\Pi^q(\Omega_0)\subseteq E^q(\Omega_0)$ 
and $\Pi^q(\Omega_\infty)\subseteq E_0^q(\Omega_\infty)$, the embeddings being isometric,
the  Cauchy transform $\mathcal C$ is a bounded operator from $L^q(\Gamma)$ into $E^q(\Omega_0)$
and from $L^q(\Gamma)$ into $E_0^q(\Omega_\infty)$. For the definition of Smirnov spaces on simply connected domains $\Omega$ within the extended complex plane $\widehat{\mathbb C}$, we refer to \cite{MR96803} or \cite{MR83565}. 
\par\smallskip
As a consequence of this result of \cite{David}, we obtain the following lemma, which will be needed in the proof of \Cref{StickingLemma}:

\begin{lemma}\label{lem:technique-bornitude-transformée-Cauchy}
Let $\Omega$ be a finitely connected bounded domain in $\mathbb C$ whose boundary satisfies the following property:
\begin{enumerate}
\item[$(\ast)$] there exist points $\mu_1,\dots,\mu_s\in\partial\Omega$ and positive numbers $r_1,\dots,r_s$ such that $\partial\Omega\subset\bigcup_{i=1}^s D(\mu_i,r_i)$ and for every $i=1,\dots,s$, the open set $\Omega\cap D(\mu_i,r_i)$ has finitely many connected components, each of which is bounded by a Jordan curve which is a Carleson curve. 
\end{enumerate}
Then for every $q>1$, the Cauchy transform defines a bounded operator from $L^q(\partial\Omega)$ into $E^q(\Omega)$. 
\end{lemma}
Note that the hypothesis $(\ast)$ is satisfied as soon as there exists a covering by open disks of the boundary of $\Omega$ such that the intersection of $\partial\Omega$ with each such disk is a union of finitely many simple Carleson curves (not necessarily closed), which intersect only at a finite number of points. In particular, this is the case if $\Omega$ is the interior of a Carleson Jordan curve. The assumption $(\ast)$ is also satisfied as soon as $\partial\Omega$ is contained in $F(\T)$, where the function $F$ satisfies \ref{H1}, \ref{H2} and \ref{H3bis}.

\begin{proof}
Let $f\in L^q(\Omega)$. For every $i=1,\dots,s$, let $f_i$ be the restriction of $f$ to $\partial\Omega\cap D(\mu_i,r_i)$. Denote by 
$\Omega_{1,i},\dots,\Omega_{k_i,i}$ the connected components of $\Omega\cap D(\mu_i,r_i)$ and by $\Gamma_{1,i},\dots,\Gamma_{k_i,i}$ their respective boundaries. Write $f_{k,i}=f_i\mathbf{1}_{\partial\Omega\cap\Omega_{k,i}}$ for $k=1,\dots,k_i$. Since $\Gamma_{k,i}$ is a Jordan curve and a Carleson curve by the hypothesis $(\ast)$, the Cauchy transform is a bounded operator from $L^q(\Gamma_{k,i})$ into $E^q(\Omega_{k,i})$. Hence the restriction of $\mathcal Cf$ to $\Omega_{k,i}$ belongs to $E^q(\Omega_{k,i})$. Since $\bigcup_{i=1}^s\bigcup_{k=1}^{k_i}\Omega_{k,i}$ is a neighborhood of $\partial\Omega$ in $\Omega$, it follows that $\mathcal Cf$ belongs to $E^q(\Omega)$ and, keeping track of the constants in the reasoning above, that $\mathcal C$ is a bounded operator from $L^q(\partial\Omega)$ into $E^q(\Omega)$. 
\end{proof}
\subsubsection{Nevanlinna class}
Let $\Omega$ be a simply connected domain of $\mathbb C$. We say that a meromorphic function $f$ on $\Omega$ belongs to the \emph{Nevanlinna class} of $\Omega$, and we write $f\in \mathcal{N}(\Omega)$, if $f$ can be written as the quotient of two functions in $H^\infty(\Omega)$. Note that if $f=g/h$, where $g,h\in E^1(\Omega)$, then $f\in \mathcal N(\Omega)$. Indeed, if $\varphi:\mathbb D\longmapsto\Omega$ is a conformal map from $\mathbb D$ onto $\Omega$, then 
\[
f\circ\varphi=\frac{g\circ\varphi}{h\circ\varphi}=\frac{(g\circ\varphi)\cdot\varphi'}{(h\circ\varphi)\cdot \varphi'},
\]
and the functions $(g\circ\varphi)\cdot \varphi'$ and $(g\circ\varphi)\cdot \varphi'$ belong to $H^1(\mathbb D)$. But a function in $H^1(\mathbb D)$ can be written as a quotient of two $H^\infty(\mathbb D)$ functions (see for instance \cite[Thm 2.1]{Duren1970}). Thus there exist $h_1,h_2\in H^\infty(\mathbb D)$ such that $f\circ\varphi=h_1/h_2$, which means that 
\[
f=\frac{h_1\circ\varphi^{-1}}{h_2\circ\varphi^{-1}},
\]
with $h_1\circ\varphi^{-1},h_2\circ\varphi^{-1}\in H^\infty(\Omega)$. Thus $f\in\mathcal N(\Omega)$.
\subsubsection{The boundary condition}
Let $\gamma$ be a subarc of $F(\mathbb T)$ containing no point of $\cal O$. Then $\gamma$ is included in the boundary of exactly two connected components $\Omega$ and $\Omega'$ of $\mathbb C\setminus F(\mathbb T)$, and 
\[|\w_F(\Omega)-\w_F(\Omega')|=1.\]
If $|\w_F(\Omega)|>|\w_F(\Omega')|$,  $\Omega$ is called the 
\textit{interior component}  and $\Omega'$ is called the \textit{exterior 
component} (with respect to $\gamma$). 
Let $\lambda_0\in\gamma$, and let $u$ be a continuous function on a neighborhood of $\lambda_0$ in $\C\setminus F(\T)$.  
We define 
(when they exist) the following two non-tangential limits of $u$ at the point $\lambda_0$, which are called respectively the \emph{interior and exterior boundary values} of $u$ at $\lambda_0$:
\[
u^{int}(\lambda_0)~=~\lim_{\substack{\lambda\to\lambda_0\\
\lambda\in \Omega}}u(\lambda)\quad\text{and}\quad u^{ext}(\lambda_0)~=~\lim_{\substack{\lambda\to\lambda_0\\
\lambda\in \Omega'}}u(\lambda),
\]

Functions which belong to a Smirnov space of a domain $\Omega$ of $\mathbb C$ (having a rectifiable boundary $\Gamma=\partial\Omega$) admit non-tangential limits almost everywhere on $\Gamma$. If $\Omega$ and $\Omega'$ are two adjacent domains along an arc $\gamma$, and if $u$ belongs to some Smirnov space $E^q(\Omega\cup\Omega')$, then the interior and exterior boundary values $u^{int}$ and $u^{ext}$ of $u$ exist almost everywhere on the arc $\gamma$.
\par\smallskip
Suppose that $F$ satisfies \ref{H1}, \ref{H2} and \ref{H3} and let $$N=\max\{|\w_F(\lambda)|\,;\,\lambda\notin F(\mathbb T)\}.$$ For each $j=0,\dots,N-1$, consider the open sets $\Omega_j^+$ given by
\begin{equation}\label{eq:def-Omegajplus}
\Omega_j^+~=~\{\lambda\notin F(\mathbb T)\,;\,|\w_F(\lambda)|>j\}.
\end{equation}
Recall that the function $\zeta=1/F^{-1}$ is defined almost everywhere on $F(\mathbb T)$. We endow the direct sum $\bigoplus_{j=0}^{N-1}E^q(\Omega_j^+)$ with the following norm:
\begin{equation}\label{eq:norm-somme-smirnov}
\left\|(u_j)_{0\le j\le N-1}\right\|~=~\left(\sum_{j=0}^{N-1}\|u_j\|_{E^q(\Omega_j^+)}^q\right)^{1/q}
\end{equation}
for every $(u_j)_{0\le j\le N-1}\in \bigoplus_{j=0}^{N-1}E^q(\Omega_j^+)$. This norm turns the space $\bigoplus_{j=0}^{N-1}E^q(\Omega_j^+)$ into a Banach space. The model space $E^q_F$ is defined as the closed subspace of $\bigoplus_{j=0}^{N-1} E^q(\Omega_j^+)$ formed by the $N$-tuples  $(u_j)_{0\le j\le N-1}$ in $\bigoplus_{j=0}^{N-1} E^q(\Omega_j^+)$ satisfying, for all $0\le j< N-1$, the following boundary conditions:
\begin{equation}\label{eq:boundary-conditions-ErF}
u_j^{int}-\zeta u_{j+1}^{int}~=~u_j^{ext}\quad\text{a.e. on}~\partial \Omega_{j+1}^+.
\end{equation}
Remark that this subspace is an invariant subspace for the multiplication operator by the independent variable $M_\lambda:\bigoplus E^q(\Omega_j^+)\to \bigoplus E^q(\Omega_j^+)$ defined for $u=(u_j)_{0\le j\le N-1}\in \bigoplus_{j=0}^{N-1} E^q(\Omega_j^+)$ by
\[
M_\lambda u~=~(v_j)_{0\leq j\leq N-1},\quad \text{with }v_j(\lambda)~=~\lambda u_j(\lambda)\quad \textrm{ for every } \lambda\in\Omega_j^+.
\]
The operator $M_\lambda:E_F^q\to E_F^q$ will be the model operator for $T_F^*\in\mathcal B(H^q)$. See Theorem~\ref{T:Yakbovich_Hp} below. 
\par\smallskip
More generally, let  $h\in H^\infty(\intsp )$ be a bounded analytic function on the interior 
of the spectrum of $T_F$. Then the space $E^q_F$ is invariant by the multiplication operator $M_{h}$ on $\bigoplus
E^q(\Omega_j^+)$ defined by $M_h(u_j)_{0\le j\le N-1}=(hu_j)_{0\le j\le N-1}$, where $u_j\in E^q(\Omega_j^+)$ for every ${0\le j\le N-1}$. See the proof of \Cref{multiplicateurs} for details. 
 This means that $H^\infty(\intsp )$ is contained in the multiplier algebra of $E^q_F$. In \Cref{Section:Multipliers}, we will see that $H^\infty(\intsp )$ is exactly the multiplier algebra of $E^q_F$.
 \par\smallskip
Let us remark here that the interior of the spectrum of $T_F$ can be described in the following way (up to the set $\mathcal O$ of intersection points of the curve):
\[\mathcal O\cup \intsp =\Omega_0^+\cup \partial\Omega_1^+. 
\]
Observe also that, for every $0\leq j\leq N-1$, we have $\partial \Omega_{j+1}^+\subset \partial\Omega_{j}^+$, and so in particular, for every $1\leq j\leq N-1$, we have $\partial\Omega_j^+\setminus \mathcal{O}\subset \intsp .$

\subsection{The model from \cite{Yakubovich1991}}
Let $p>1$, and let $q$ be the conjugate exponent of $p$. Suppose that the symbol $F$ of the Toeplitz operator $T_F\in\cal B(H^p)$ satisfies the three assumptions \ref{H1}, \ref{H2} and \ref{H3}. Let $h_{\lambda,j}$, $0\le j\le N-1$, be given by (\ref{eve}). For every function $g\in H^q$, define $Ug=((Ug)_j)_{0\le j\le N-1}$ by setting
\begin{equation}\label{defn-model-U-section2-5}
(Ug)_j(\lambda)~=~\dual{h_{\lambda,j}}{g}~\quad\textrm{for every }\lambda\in\Omega_j^+.
\end{equation}

Note that since $h_{\lambda,j}$ is an eigenvector of $T_F$ associated to the eigenvalue $\lambda$ as soon as $\lambda\in\Omega_j^+$, we have for every $g\in H^q$ and every $0\le j\le N-1$ that
\[\dual{h_{\lambda,j}}{T_F^*g}~=~\dual{T_Fh_{\lambda,j}}{g}~=~\lambda\dual{h_{\lambda,j}}{g}\quad \text{ for every }\lambda\in\Omega_j^+.\]
In other words, \[U(T_F^*g)~=~M_\lambda(Ug) \quad \text{ for every } g\in H^q.\]
Using the expression of the functions $h_{\lambda,j}$, $\lambda\in\Omega_j^+$, given by \cref{eve}, combined with another deep expression of the function $F_\lambda^+$ (whose proof uses, in particular, tools from quasiconformal mapping theory), one can show that for every $z\in\mathbb{D}$, the function $\lambda\longmapsto h_{\lambda,j}(z)$ is in $E^q(\Omega_j^+)$ and that we have
\begin{equation}\label{eq:boundary-conditions-ErF-eigenvectors}
h_{\lambda,j}^{int}(z)-\zeta(\lambda) h_{\lambda,j+1}^{int}(z)~=~h_{\lambda,j}^{ext}(z)\quad\textrm{for almost every }\lambda\in\partial \Omega_{j+1}^+.
\end{equation}
See Corollary B.12 in \cite{FricainGrivauxOstermann_preprint}. In other words, for every $z\in\mathbb{D}$, the $N$-tuple 
\[
\lambda\longmapsto (h_{\lambda,0}(z), h_{\lambda,1}(z),\ldots, h_{\lambda,N-1}(z))
\]
belongs to $E_F^q$. 
\par\smallskip
It was shown in \cite{Yakubovich1991} that the operator $U$ defined by \eqref{defn-model-U-section2-5} is bounded from $H^q$ into $E^q_F$ (see also \cite{FricainGrivauxOstermann_preprint}*{Th. B.17}). Note that the detailed proof contained in \cite{FricainGrivauxOstermann_preprint} and, more specifically, the combination of Equation (B.31) and Fact B.18 of \cite{FricainGrivauxOstermann_preprint}, shows that there exists a neighborhood $V$ of the curve $F(\mathbb T)$ and two constants $C_1,C_2>0$ such that, for every $\lambda \in V \cap \Omega_0^+$, there exist points $z_1, \dots, z_s \in \mathbb D$ such that for every $g \in H^q$
\begin{equation}\label{Eq:eo3u}
|u_0(\lambda)| \le C_1\|g\|_{H^q} + C_1 \sum_{j=1}^s |g(z_j)|,
\end{equation}
where $Ug = (u_0, \dots, u_{N-1})$. Note that the number $s$ of points $z_j$ involved in this inequality might depend on $\lambda$, but it is uniformly bounded. In particular, this implies that if $g \in H^\infty$, then $u_0$ is bounded on $V\cap\Omega_0^+$. But since $u_0$ is also analytic in $\Omega_0^+$, we finally deduce that $u_0 \in H^\infty(\Omega_0^+)$. This therefore yields the following fact, which will be used several times in the rest of the paper:
\begin{fact}\label{elements-de-l-espace-modele} Let $0\le k<N$ and $u=Uz^k$. Then $u$ has the form
\[u=\begin{cases}
(1,0,\dots,0)&\text{if}~k=0,\\
(u_0,\dots,u_{k-1},1,0,\dots,0)&\text{if}~1\le k<N,
\end{cases}\]
where the function $u_\ell$ belongs to $ H^\infty(\Omega_\ell^+)$ for every $0\le \ell\le k-1$.
In particular, the $N$-tuple $U1=(1, 0,\ldots, 0)$ belongs to $E_F^q$.
\end{fact}

\begin{proof}
Let $0\le k<N$, and $u=Uz^k=(u_0,\dots,u_{N-1})$. Then, according to \eqref{eve}, for every $0\le \ell<N$ and every $\lambda\in\Omega_\ell^+$, we have
\[u_\ell(\lambda)=\ps{h_{\lambda,\ell},z^k}_{p,q}=\ps{z^\ell h_{\lambda,0},z^k}_{p,q}.\]
In particular, $u_k(\lambda)=h_{\lambda,0}(0)=1$ and for every $\ell>k$, $u_\ell=0$. 
Fix now $0\leq \ell <k$. Then, for every $\lambda\in\Omega_\ell^+$, we have $u_\ell(\lambda)=\ps{ h_{\lambda,0},z^{k-\ell}}_{p,q}$. Since $z^{k-\ell}\in H^\infty$, it follows from \eqref{Eq:eo3u} that $u_\ell$ is bounded on $\Omega_\ell^+$ and since $u_\ell\in E^q(\Omega_\ell^+)$, we finally obtain that $u_\ell\in H^\infty(\Omega_\ell^+)$, which concludes the proof.
\end{proof}

Here is now the statement of the main theorem of \cite{Yakubovich1991} for $p=2$. See also  \cite{FricainGrivauxOstermann_preprint}*{Th. B.27} for the $H^p$ case.

\begin{theorem}[Yakubovich \cite{Yakubovich1991}]\label{T:Yakbovich_Hp}
Let $F$ be a symbol satisfying \ref{H1}, \ref{H2} and \ref{H3}, and let $T_F$ be the Toeplitz operator with symbol $F$ acting on $H^p$. Then the operator $U$ defined in \eqref{defn-model-U-section2-5} is a linear isomorphism from $H^q$ onto $E_F^q$. Moreover, we have
\[T_F^*~=~U^{-1}M_\lambda U.\]
\end{theorem}
It follows from this theorem and from \Cref{Prop:Reflexif} is that it is equivalent to study the embeddability of $T_F$ and that of the multiplication operator by $\lambda$ on a suitable model space. More precisely:

\begin{proposition}\label{equivalent}
Let $p>1$, and let $q$ be its conjugate exponent. Suppose that the symbol $F$ satisfies \ref{H1}, \ref{H2} and \ref{H3}, and let $f(z)=F(\frac{1}{z})$, $z\in\mathbb{T}$. The following assertions are equivalent:\par\smallskip
\begin{enumerate}
\item the Toeplitz operator $T_F$ acting on $H^p$ is embeddable into a $C_0$-semigroup of operators on $H^p$;\par\smallskip
\item the Toeplitz operator $T_f=T_F^*$ acting on $H^q$ is embeddable into a $C_0$-semigroup of operators on $H^q$;\par\smallskip
\item the operator $M_\lambda$ acting on $E_F^q$ is embeddable into a $C_0$-semigroup of operators on $E_F^q$.
\end{enumerate}
\end{proposition}

\begin{proof}
Suppose first that $F$ satisfies \ref{H3}. By  \Cref{Prop:Reflexif}, $T_F$ acting on $H^p$ is embeddable if and only if $T_F^*$ acting on $H^q$ is embeddable, and by \Cref{T:Yakbovich_Hp} this is the case if and only if 
$M_\lambda$ acting on  $E^q_F$ is embeddable.
\end{proof}

\begin{remark}
 An important consequence of \Cref{equivalent} is that assumptions \ref{H3} and \ref{H3bis} are morally equivalent when dealing with the question of the embeddability of Toeplitz operators into $C_0$-semigroups. Indeed,  if $F$ satisfies \ref{H3} (all winding numbers are nonpositive), then $T_F$ acting on $H^p$ is embeddable if and only if the model operator $M_\lambda$ acting on the model space $E^q_F$ is embeddable. If $F$ satisfies \ref{H3bis} but not \ref{H3} (all winding numbers are nonnegative), then $T_F$ acting on $H^p$ is embeddable if and only if  $M_\lambda$ acting on $E^p_f$ is embeddable. In the rest of the paper, many results will be stated and proved under the assumption \ref{H3}, but such results hold as well under assumption \ref{H3bis}.   
\end{remark}

\begin{remark}
   We finish this section with a last observation which turns out to be quite useful when constructing concrete examples: let $F_1$ and $F_2$ be two symbols satisfying \ref{H1}, \ref{H2} and \ref{H3} (or, more generally \ref{H3bis}), and such that there exists an orientation-preserving diffeomophism $\phi:\T\to\T$  such that 
$F_2=F_1\circ \phi$. Suppose that for every $z\in\T$ such that $F_1(z)$ does not belong to the boundary of $\sigma(T_{F_1})=\sigma(T_{F_2})$, $\phi(z)=z$ (so that $F_2(z)=F_1(z)$).
Define $\zeta_1=1/F_1^{-1}$ and $\zeta_2=1/F_2^{-1}$ on $F_1(\T)\setminus\mathcal{O}$, where $\mathcal{O}$ is the set of self-intersection points of the curve. Then $\zeta_1(\lambda) = \zeta_2(\lambda)$ for every $\lambda\in F_1(\T)\setminus\mathcal{O}$ not belonging to the boundary of $\sigma(T_{F_1})$, and thus the two model spaces $E^q_{F_1}$ and $E^q_{F_2}$ coincide. It follows then from \Cref{T:Yakbovich_Hp} that $T_{F_1}$ and $T_{F_2}$ are similar, and thus $T_{F_2}$ is embeddable if and only if $T_{F_1}$ is. 
\end{remark}

\subsection{Riemann surfaces models for a special class of rational symbols}
In the work \cite{Yakubovich1989}, which preceded \cite{Yakubovich1991}, a Riemann surface model was constructed for Toeplitz operators whose symbols are positively wound and have the form 
$F=P/Q$, where $P\in A(\D)$, $Q$ is a polynomial 
with roots in $\D$, and $F$ is locally univalent on $\D$ near $\T$. Several other more technical conditions on $F$, which we do not reproduce here, were imposed in \cite{Yakubovich1989}. 
\par\smallskip
This model gave rise to a complete description of the commutant of $T_F$, which permits one to give a comprehensive characterization of embeddable Toeplitz operators for this restricted class of symbols. It was proved in \cite[Ch. 2]{Yakubovich1989} that for this class of symbols $F$, there exists a branched cover $(\sic, \rhoc)$ of $\intsp$ such that the commutant $\{T_F\}'$ of $T_F$ is isomorphic, as a Banach algebra, to 
$H^\infty(\sic)$. The bordered Riemann surface in general is assumed to have a finite number of connected components; it is not asserted that it is connected. It follows, in particular, 
that $\{T_F\}'$ coincides with the set 
$\{g(T_F)\,;\, g\in H^\infty(\intsp)\}$ 
if and only if $\rhoc$ is an isomorphism of $\sic$ onto $\intsp$. 
By using this result, it can be shown that 
$T_F$ is embeddable in a $C_0$-semigroup if and only if $0\notin \intsp$ and for any closed curve 
$\ga$ on $\sic$, the winding number of the closed curve $\rhoc\circ\ga$ is $0$. 
\par\smallskip
Our aim being here to deal with a larger class of smooth symbols, we concentrate on the model from \cite{Yakubovich1991} and its consequences on embeddability, and we do not expose the results above in detail in this paper.

\section{Multipliers}\label{Section:Multipliers}
\subsection{Characterization of the multiplier algebra of the model space}
Let $1<p<+\infty$. Given a function $F$ satisfying the assumptions 
\ref{H1}, \ref{H2} and \ref{H3}, observe that the spectrum of $T_F$ does not depend on the Hardy space $H^p$ on which it acts; so we will simply speak about the spectrum of $T_F$, without mentioning the space. 
Recall also that $\Omega_0^+=\sigma(T_F)\setminus F(\mathbb{T})$.
\par\smallskip
Let $1<p<+\infty$ and $q$ be its conjugate exponent, i.e. $\frac1p+\frac1{q}=1$. Consider a function $g:\Omega_0^+\longrightarrow\mathbb C$. We say that $g$ is a \emph{multiplier} of $E_F^q$ if for any element $u=(u_0,\dots,u_{N-1})$ of $E_F^q$, the $N$-tuple $gu=(gu_0,\dots,gu_{N-1})$ also belongs to $E_F^q$. We denote by $\mathcal Mult(E_F^q)$ the \emph{multiplier algebra} of $E_F^q$, i.e. the set of all multipliers of $E_F^q$.
Whenever $g$ is a multiplier on $E_F^q$, a straightforward application of the closed graph theorem ensures that the multiplication operator by $g$ is bounded on $E_F^q$; we denote this operator by $M_g$. Of course, other natural definitions of multipliers could be considered in this context, but we have chosen this one as it appears to be the most suitable for our study. 
\par\smallskip
The proof of our characterization of the multipliers of $E_F^q$ relies on the following sticking lemma, a version of which is mentioned in \cite{Yakubovich1993} (see also \cite{SolomyakVolberg1989}*{Th. 7.4}) as a consequence of the Cauchy formula.

\begin{lemma}[Sticking lemma]\label{StickingLemma}
Let $q>1$, and let $\Omega_1,\Omega_2$ be disjoint finitely connected domains of $\mathbb{C}$ with piecewise smooth boundaries such that $\alpha=\partial\Omega_1\cap\partial\Omega_2$ has a positive measure. Let $\Omega=\mathrm{int}(\Omega_1\cup\Omega_2\cup \alpha)$ and assume that $\Omega$ is a finitely connected domain which satisfies the hypothesis $(\ast)$ of \Cref{lem:technique-bornitude-transformée-Cauchy}. Let $f_j\in E^q(\Omega_j)$, $j=1,2$, be such that $f_1=f_2$ a.e. on $\alpha$. Then there exists a function $f\in E^q(\Omega)$ such that $f_{|\Omega_j}=f_j$ for $j=1,2$.
\end{lemma}
Note that \Cref{StickingLemma} will always be applied later on in the paper to domains $\Omega_j$ which satisfy $\partial\Omega_j\subseteq F(\mathbb T)$, with symbols $F$ satisfying \ref{H1}, \ref{H2} and \ref{H3bis}. In  this case, the set $\Omega$ coincides with $\mathrm{int}(\overline{\Omega}_1\cup\overline{\Omega}_2)$, and the hypothesis that 
$\Omega$ is finitely connected will always be satisfied. As mentioned after \Cref{lem:technique-bornitude-transformée-Cauchy}, the second assumption $(\ast)$ will also always be satisfied in this setting. Indeed, since $F'$ does not vanish on $\T$,  one can cover $\mathbb T$ with a finite number of small arcs $\alpha_j$ such that $F(\alpha_j)$ are Carleson simple curves. Also, note that \Cref{StickingLemma} easily implies (using induction) a similar statement involving $k$ domains $\Omega_1,\ldots,\Omega_k$ with $k\ge 2$.
\par\smallskip
Now, for completeness' sake, we include a detailed proof of \Cref{StickingLemma}.
\begin{proof}
Let $j=1$ or $j=2$. Since $f_j\in E^q(\Omega_j)$, by Th. 10.4 in \cite{Duren1970},
we have 
\begin{equation}\label{Eq:Sticking1}
    f_j(z)=\frac1{2i\pi}\int_{\partial\Omega_j}\frac{f_j(\zeta)}{\zeta-z}\,\mathrm d\zeta\quad\mbox{for every }z\in\Omega_j.
\end{equation}
See \cite[Section 10.5]{Duren1970}. Moreover, since $\Omega_1\cap\Omega_2=\varnothing$, we also have 
\begin{equation}\label{Eq:Sticking2}
    \forall\,z\in\Omega_2,\,\int_{\partial\Omega_1}\frac{f_1(\zeta)}{\zeta-z}\,\mathrm d\zeta=0 ~\text{and}~\forall\,z\in\Omega_1,\,\int_{\partial\Omega_2}\frac{f_2(\zeta)}{\zeta-z}\,\mathrm d\zeta=0.
\end{equation}
Note also that, since $f_1=f_2$ a.e. on $\alpha=\partial\Omega_1\cap\partial\Omega_2$, we deduce that 
\begin{equation}\label{eq:2ZEDSFSF}
\int_\alpha\frac{f_1(\zeta)}{\zeta-z}\,\mathrm d\zeta=\int_\alpha\frac{f_2(\zeta)}{\zeta-z}\,\mathrm d\zeta\quad\text{for all $z\in\Omega_1\cup\Omega_2$}.
\end{equation}
Let $u=f_j$ a.e. on $\partial\Omega\cap\partial\Omega_j$. Then $u\in L^q(\partial\Omega)$ and we can define an analytic function $f$ on $\Omega$ by setting
\begin{equation}\label{eq:sqdqsd93ZEE}
f(z)=\frac1{2i\pi}\int_{\partial\Omega}\frac{u(\zeta)}{\zeta-z}\,\mathrm d\zeta\quad \text{for every}~z\in \Omega.
\end{equation}
Note that the orientations of $\alpha\subseteq\partial\Omega_1$ and $\alpha\subseteq\partial\Omega_2$ differ. According to \eqref{eq:2ZEDSFSF}, and noting that our definition of $\Omega$ implies that $\partial\Omega_1\cup\partial\Omega_2=\partial\Omega\cup\alpha$,  we get that 
\begin{equation}\label{Eq:Sticking3}
f(z)=\frac1{2i\pi}\int_{\partial\Omega_1}\frac{f_1(\zeta)}{\zeta-z}\,\mathrm d\zeta+\frac1{2i\pi}\int_{\partial\Omega_2}\frac{f_2(\zeta)}{\zeta-z}\,\mathrm d\zeta\quad\text{for all $z\in\Omega_1\cup\Omega_2$.}
\end{equation}
Let $z\in\Omega_j$. By
\eqref{Eq:Sticking2}, the term in \eqref{Eq:Sticking3} which is not equal to 
$\frac1{2i\pi}\int_{\partial\Omega_j}\frac{f_j(\zeta)}{\zeta-z}d\zeta$ is zero, and thus by \eqref{Eq:Sticking1}, we have that, for all $z\in\Omega_j$,
\begin{align*}
f(z)\,
=\,\frac1{2i\pi}\int_{\partial\Omega_j}\frac{f_j(\zeta)}{\zeta-z}\,\mathrm d\zeta
=\,f_j(z).
\end{align*}
We have thus shown that $f_{|\Omega_j}=f_j$. To finish the proof, it remains to prove that $f\in E^q(\Omega)$. But, according to \eqref{eq:sqdqsd93ZEE}, we have $f=\mathcal C u$ on $\Omega$, and, since $u\in L^q(\partial\Omega)$, \Cref{lem:technique-bornitude-transformée-Cauchy} implies that $f\in E^q(\Omega)$.
\end{proof}
We are now ready to prove:
\begin{theorem}\label{multiplicateurs}
The multiplier algebra of $E_F^q$ coincides with the set of restrictions to $\Omega_0^+$ of bounded analytic functions on the interior of the spectrum  of $T_F$, that is 
\[
\mathcal Mult(E_F^q)=\{g_{|\Omega_0^+}:g\in H^\infty(\intsp )\}.
\]
\end{theorem}

\begin{proof}
Let $g\in \mathcal Mult(E_F^q)$. Let us first show that $g$ belongs to $H^\infty(\Omega_0^+)$.
Since $(1,0,\dots,0)$ belongs to $E_F^q$ by \Cref{elements-de-l-espace-modele}, and since $g$ is supposed to be a multiplier of $E_F^q$,  the $N$-tuple $(g,0,\dots,0)=M_g(1,0,\dots,0)$ belongs to $E_F^q$. Hence $g\in E^q(\Omega_0^+)$.
Fix now $\lambda\in \Omega_0^+$. The evaluation map at the point $\lambda$ is a continuous linear form on $E^q(\Omega_0^+)$, whence it follows that the linear form $\psi_\lambda:u=(u_0,\dots,u_{N-1})\longmapsto u_0(\lambda)$
is continuous on $E_F^q$. 
 It is easy to see that 
$
M_g^* \psi_\lambda=g(\lambda)\psi_\lambda
$,  
which yields that
\[
|g(\lambda)|\leq
\frac {\|M_g^*\psi_\lambda\|}{\|\psi_\lambda\|}\le 
\|M_g^*\|. 
\]
Since this estimate does not depend on $\lambda\in \Omega_0^+$ and since
$g$ is analytic on $\Omega_0^+$, we conclude that 
$g\in H^\infty(\Omega_0^+)$ and that $\|g\|_\infty\le\|M_g\|$. 
\par\smallskip
Recall now that 
\[
\intsp\cup\mathcal O =\Omega_0^+\cup \partial\Omega_1^+,
\]
and so the only thing to check is that $g$ can be extended into an analytic function on 
$\intsp $. Using one more time that $(g,0,\dots,0)\in E_F^q$, it follows from the boundary relations in $E_F^q$ that $g$ satisfies $g^{int}=g^{ext}$ a.e. on $\partial\Omega_1^+$. 
Applying \Cref{StickingLemma} to each arc contained in $\partial\Omega_1^+\setminus\mathcal{O}$, we obtain that  $g$ can be extended analytically to $\intsp $, and thus 
\[
\mathcal Mult(E_F^q)\subseteq \{g_{|\Omega_0^+}:g\in H^\infty(\intsp )\}.
\]
\par\smallskip
Now suppose that $g$ is a function lying in $H^\infty(\intsp )$. Consider an element
$u=(u_0,\dots,u_{N-1})$ of $E_F^q$. Since $g$ is bounded, the function $gu_j$ lies in $E^q(\Omega_j^+)$ for every $j=0,\ldots, N-1$. Moreover, since for every $\lambda\in \partial\Omega_{1}^+\setminus\mathcal O$, the function $g$ is continuous at the point $\lambda$, we have that $g=g^{int}=g^{ext}$ on $\partial\Omega_{j+1}^+\setminus\mathcal O$ for any $0\le j< N-1$ and so 
\[(gu_j)^{int}-\zeta (gu_{j+1})^{int}=g(u_j^{int}-\zeta u_{j+1}^{int})=gu_j^{ext}=(gu_j)^{ext}~\text{a.e. on }\partial\Omega_{j+1}^+.\]
Hence $gu$ belongs to $E^q_F$, and we have thus shown that $g$ belongs to $\mathcal Mult(E_F^q)$. 
So we finally deduce that 
\[
\mathcal Mult(E_F^q)\supseteq 
\{g_{|\Omega_0^+}:g\in H^\infty(\intsp )\},
\]
and \Cref{multiplicateurs} is proved.
\end{proof}

\begin{remark}
As a nice consequence of \Cref{multiplicateurs}, we obtain that when $F$ is sufficiently smooth (for instance when $F$ is of class $C^2$ on $\mathbb T$), the set of multipliers of $E_F^q$ depends only on the curve $F(\mathbb T)$. It does not depend on $1<q<\infty$,
nor on the choice of the sufficiently smooth parametrization $F$ of this curve (although the model space $E_F^q$ itself depends on the parametrization $F$ through the boundary conditions).     
\end{remark}

\subsection{Embedding of the model operator in a $C_0$-semigroup of multiplication operators and consequences}
Recall that there exists an analytic determination of the square root on a domain $\Omega$ of $\mathbb{C}$ if and only if $\Omega$ is contained in a simply connected domain which does not contain $0$; this is equivalent to saying that $0$ belongs to an unbounded connected component of $\mathbb C\setminus \Omega$. This remark lies at the core of the proof of the following theorem:
\TheoEmbModelMult*

\begin{proof}
Suppose first that $M_\lambda$ is embeddable into a $C_0$-semigroup $(A_t)_{t>0}$ of multiplication operators on $E_F^q$. In particular, $A_{1/2}$ is a multiplication operator, and thus there exists a function $\delta\in H^\infty(\intsp )$ such that $A_{1/2}=M_\delta$. Then we have
\[M_\lambda=A_1=(A_{1/2})^2=(M_\delta)^2=M_{\delta^2}.\]
By applying this equality to the element $(1,0,\dots,0)$ of $E_F^q$, we get that for every $\lambda\in \Omega_0^+$, we have $\delta^2(\lambda)=\lambda$. Since  $\delta$ is analytic on $\intsp $, it follows that this equality is also true for every $\lambda\in \intsp $, and thus  $\delta$ is an analytic determination of the square root on $\intsp $. Hence $0$ belongs to the unbounded component of $\mathbb C\setminus \intsp $.
\par\smallskip
Conversely, suppose that $0$ belongs to the unbounded component of the set $\mathbb C\setminus \intsp $. Then there exists  an analytic determination of the logarithm on $\intsp $, which we write as $\log$. Let $v_t(\lambda)=\lambda^t:=e^{t\log\lambda}$, $t>0$. Since the function $v_t$ belongs to $H^\infty(\intsp )$, the operator $A_t:=M_{v_t}$ is well-defined and bounded on $E_F^q$. By construction, $(A_t)_{t>0}$ is a semigroup satisfying $A_1=M_\lambda$, so we just need to prove its strong continuity. 
Let $u=(u_0,\dots,u_{N-1})$ be an element of $E_F^q$. Then
\[\|A_tu-u\|^q_{E_F^q}=\sum_{j=0}^{N-1}\|(\lambda^t-1)u_j\|_{E^q(\Omega_j^+)}^q=\sum_{j=0}^{N-1}\sum_{\underset{\text{component}}{\Omega\subset\Omega_j^+}}\|(\lambda^t-1){u_j}_{|\Omega}\|_{E^q(\Omega)}^q.\]
But recall that the norm on the Smirnov space $E^q(\Omega)$ is given by an integral on $\partial\Omega$, so, by Lebesgue dominated convergence theorem, $\|A_tu-u\|_{E_F^q}\longrightarrow0$ as $t\to0$.
\end{proof}
When $0$ belongs to the interior of the unbounded component of $\mathbb C\setminus \intsp$ we have a stronger statement: in this case $\log\in H^\infty(\intsp )$, and thus $A_t=e^{tB}$ with $B=M_{\log}\in B(E_F^q)$. So the semigroup $(A_t)_{t>0}$ has a bounded generator, and $\|A_t-I\|\longrightarrow0$ when $t\to 0$.
\par\smallskip

As a direct consequence of \Cref{Th:EmbModelMult} and \Cref{equivalent}, we obtain the following sufficient condition for a Toeplitz operator with a symbol satisfying \ref{H1}, \ref{H2} and \ref{H3bis} to be embeddable into a $C_0$-semigroup:
\CSPlongement*

We have provided in \Cref{example1:embeddable-bounded-compopent} (see also the forthcoming \Cref{{l'exemple!!}}) a symbol $F$ showing that the converse of \Cref{Th:CSforEmb} is not true in general: $T_F$ may be embeddable although $0$ belongs to a bounded connected component of $\mathbb C\setminus \intsp $.
But right now, we study specific situations where the converse of \Cref{Th:CSforEmb} does hold, i.e. the embeddability of $T_F$ into a $C_0$-semigroup forces $0$ to belong to the unbounded component of $\mathbb C\setminus \intsp $. We first study the case where $F(\mathbb{T})$ is a Jordan curve.

\subsection{The case of a Jordan curve}
In the case where $F(\mathbb{T})$ is a Jordan curve, the requirements \ref{H1}, \ref{H2} and \ref{H3bis} on $F$ are equivalent to 
\ref{H1} (the other two assumptions are automatically true).

\begin{theorem}\label{Th:Jordan}
Let $p>1$, and let $F\in L^{\infty}(\mathbb{T})$ be a symbol satisfying \ref{H1}. Assume furthermore that $F(\mathbb T)$ is a Jordan curve. Then the following assertions are equivalent:\par\smallskip
\begin{enumerate}
    \item $T_{F}\in\mathcal{B}(H^p)$ is embeddable into a $C_0$-semigroup.\par\smallskip
    \item $0$ belongs to $\mathbb{C}\setminus\intsp $.
\end{enumerate}

\end{theorem}
\begin{proof}
The implication $(2)\implies (1)$ follows immediately from \Cref{Th:CSforEmb} since $\mathbb{C}\setminus\intsp $ is connected and unbounded.
\par\smallskip
For the reverse implication $(1)\implies (2)$, note that $\intsp =\sigma(T_F)\setminus F(\mathbb T)$ (because $F(\mathbb T)$ is a Jordan curve), and apply \Cref{Fact:RemOnPos0} (i). 
\end{proof}

The regularity assumptions on $F$ are necessary for \Cref{Th:Jordan} to hold. Indeed, if we only suppose that $F$ is continuous on $\mathbb T$, then we have the following counter-example.

\begin{example}\label{ex:root-3}
Let $F(z)=\frac{(1-z)^{1/3}}z$, $z\in\mathbb{T}$, where the $1/3$-root is defined using the principal determination of the logarithm on $\mathbb{C}\setminus(-\infty, 0]$. Then
$F(\bb T)$ is a negatively wound Jordan curve and the point $0$ lies on $F(\bb T)$.
\par\smallskip
Indeed, we have $1-e^{i\theta}=-2i\sin(\theta/2)e^{i\theta/2}=2\sin(\theta/2)e^{i(\theta/2-\pi/2)}$. Then, if $0\le \theta\le 2\pi$, we have $2\sin(\theta/2)\ge0$, and since $\theta/2-\pi/2\in[-\pi/2,\pi/2]\subseteq (-\pi,\pi)$, we have that \((1-e^{i\theta})^{1/3}=\sqrt[3]{2\sin(\theta/2)}e^{\frac i6(\theta-\pi)}\) and so
\[F(e^{i\theta})~=~\sqrt[3]{2\sin(\theta/2)}e^{\frac{-i}6(5\theta+\pi)}~\textrm{ for every }\, 0\le\theta\le2\pi.\]
Note that the function $u$ defined on $[0, 2\pi]$ by $u(\theta)=\frac{-1}6(5\theta+\pi)$ is strictly decreasing, $u(0)=-\pi/6$ and $u(2\pi)=-11\pi/6=-2\pi+\pi/6$. So $F(\bb T)$ is a negatively wound Jordan curve.

\begin{figure}[ht]
    \centering
    \begin{tikzpicture}[scale=1.5]
     \draw[->,domain=0:pi/2, samples=500, smooth, variable=\x]plot({(2*abs(sin(\x/2 r)))^(1/3)*cos((5*\x+pi)/6 r)},
     {(2*abs(sin(\x/2 r)))^(1/3)*sin(-(5*\x+pi)/6 r))});
     \draw[domain=pi/2:3*pi/2, samples=500, smooth, variable=\x]plot({(2*abs(sin(\x/2 r)))^(1/3)*cos((5*\x+pi)/6 r)},
     {(2*abs(sin(\x/2 r)))^(1/3)*sin(-(5*\x+pi)/6 r))});
      \draw[domain=0:pi/2, samples=500, smooth, variable=\x]plot({(2*abs(sin(\x/2 r)))^(1/3)*cos((5*\x+pi)/6 r)},
     {(2*abs(sin(\x/2 r)))^(1/3)*sin((5*\x+pi)/6 r))});
     \draw[->](-1.5,0)--(1,0);
     \draw[->](0,-1.2)--(0,1.2);
    \end{tikzpicture}
    \caption{}
    \label{fig1}
\end{figure}
Now, consider the Toeplitz operator $T_F$ with symbol $F$ acting on $H^p$, $1<p<3$. Then we have $\dim(\ker T_F)=1$.
Indeed, a function $u\in H^p$ belongs to $\ker(T_F)$ if and only if $ (1-z)^{1/3}u$ belongs to $\ker S^*$ (where $S$ is the multiplication operator by $z$ on $H^p$), i.e. if and only if there exists a constant $c\in\bb C$ such that $u(z)=\frac{c}{(1-z)^{1/3}}$. Since $(1-z)^{-1/3}$ belongs to $H^p$, we obtain that $\ker(T_F)=\spa\big[(1-z)^{-1/3}\big]$.
It then follows from \Cref{Th:CN-plong-ker} that $T_F$ is not embeddable into a $C_0$-semigroup on $H^p$. 
\par\smallskip
If $p>3$, let $q$ denote the conjugate exponent of $p$, so that $1<q<3$. Then $T_F$ acting on $H^{q}$ is not embeddable into a $C_0$-semigroup on $H^{q}$ by the argument above, and thus $T_F^*\in\mathcal{B}(H^p)$ is not embeddable into a $C_0$-semigroup on $H^p$ either. Setting $f(z)=F(1/z)$, $z\in\T$, the symbol $f$ has the same regularity as $F$, and $T_f=T_F^*$ is not embeddable into a $C_0$-semigroup on $H^p$.
\end{example}

We will come back in \Cref{Section: sectorial} to the study of the embeddability problem for $T_F$ in the general situation where the symbol $F$ is not assumed to be $C^1$ smooth.

\section{Embedding and eigenvalues on the curve}\label{Section 5}
\subsection{A particular case of a result of Ahern and Clark}
Ahern and Clark studied in \cite{AhernClark1985} the dimension of the kernel of  non-Fredholm Toeplitz operators on $H^2$. When the symbol $F$ is differentiable on $\T$ with $F'\neq0$ on $\T$, their results apply, and can be reformulated to yield the following statement:
\begin{theorem}[Ahern - Clark \cite{AhernClark1985}]\label{Ahern-Clark}
Let $p>1$, let $F$ be a differentiable function on $\T$ with $F'\neq0$ on $\mathbb T$ and let $T_F\in\mathcal B(H^p)$. Then
$\dim\ker(T_F)=\max(0,-\w_+(F))$, where
\[\w_+(F)=\frac1{2\pi}\left(\Delta\arg(F)+m\pi\right),\text{ with }m=\card\{\zeta\in\mathbb T\,;\,F(\zeta)=0\}.\]
\end{theorem}

In the statement of \Cref{Ahern-Clark}, $\Delta\arg(F)$ denotes the variation of the argument of $F(e^{i\theta})$ as $\theta$ varies from $0$ to $2\pi$, ``forgetting jumps". Let us explain more precisely what is meant here:
\smallskip

$\bullet$\enspace If $0$ does not belong to $F(\mathbb T)$, $\Delta\arg(F)$ is simply the variation of the argument of $F(e^{i\theta})$ as $\theta$ varies from $0$ to $2\pi$ (there are no jumps in this case), and we have $\w_+(F)=\w_F(0)$. So the formula given in \Cref{Ahern-Clark} generalizes the classical formula for the dimension of the kernel of a Fredholm Toeplitz operator (see~\eqref{dimension-noyau-Fredholm}). 

\smallskip
$\bullet$\enspace Suppose now that $0$ belongs to $F(\mathbb T)$.
Under the assumptions of \Cref{Ahern-Clark}, $F$ only has a finite number of zeroes on $\T$, as explained in the following fact:

\begin{fact}\label{nb-zeros-fini}
Under the assumptions of \Cref{Ahern-Clark}, the function $F$ has a finite number of zeroes on $\mathbb T$.
\end{fact}

\begin{proof}
Indeed, if $F$ had infinitely many zeros on $\mathbb T$, 
then they would have an accumulation point $\zeta\in\mathbb T$, which would imply that $F'(\zeta)=0$. 
This contradicts the assumption 
$F'\ne0$ on $\mathbb T$. 
\end{proof}
Let $z_1,\dots, z_m$, with $m\ge 1$ be the zeroes of $F$ on $\T$. Write each such zero as
$z_j=e^{i\theta_j}$, where $0\le\theta_1<\theta_2<\ldots\theta_m<2\pi$, and set $\theta_{m+1}=\theta_1$. The function $\arg F(e^{it})$ on $[0,2\pi]$ has jumps only at the
points $\theta_j$. When computing the variation $\Delta\arg(F)$ of the argument of $F$, we forget about the jumps of the argument at the points $\theta_j$, $j=1, \ldots, m$, and define $\Delta\arg(F)$ as the sum over $j=1, \ldots, m$ of the variations  of the argument of $F(e^{i\theta})$ as $\theta$ varies from $\theta_j$ to $\theta_{j+1}$.
\par\medskip
Note that the result of Ahern-Clark is stated in \cite{AhernClark1985} in the context of $H^2$, but it holds true in $H^p$ for every $p>1$ as well.
For the sake of completeness, we include a proof of \Cref{Ahern-Clark} in the case where  $0$ belongs to $F(\mathbb T)$.

\begin{proof} Suppose that $0\in F(\mathbb T)$. Then $m=\card\{\zeta\in\mathbb T\,;\,F(\zeta)=0\}$ is at least $1$; it is finite by \Cref{nb-zeros-fini}.
So let $\zeta_1,\dots,\zeta_m\in\mathbb T$ be the zeroes of $F$ on $\T$. Since $F$ is differentiable on $\T$ and $F'\neq0$, there exists a continuous function $g$ on $\mathbb T$ that does not vanish on $\mathbb T$ and is such that 
\[F(e^{i\theta})=g(e^{i\theta})\prod_{k=1}^m(e^{i\theta}-\zeta_k)= e^{im\theta}g(e^{i\theta})\prod_{k=1}^m(1-\zeta_ke^{-i\theta}).\]
We write
\(F(e^{i\theta})
=g_0(e^{i\theta})\overline{q(e^{i\theta})},\)
with $g_0(e^{i\theta})=e^{im\theta}g(e^{i\theta})$ and $q(e^{i\theta})=\prod_{k=1}^m(1-\overline{\zeta_k}e^{i\theta})$, $e^{i\theta}\in\T$. 
Since  the function $\overline{q}$ is anti-analytic on $\mathbb D$, by \cite{MR3497010}*{Th. 12.4}, we deduce that 
\[T_F=T_{\overline{q}}T_{g_0}.\]
Moreover, $q$ is a polynomial and does not vanish on $\mathbb D$. So $q$ is an outer function and thus $T_q$ has a dense range, which means that $T_{\overline{q}}$ is injective (see \cite{MR3497010}*{Th. 12.19}). Thus we have that $\ker(T_F)=\ker(T_{g_0})$, and hence
\[\dim\ker(T_F)=\dim\ker(T_{g_0})=\max(-\w_{g_0}(0),0)\]
since $g_0$ is continuous and non-vanishing on $\T$.
Finally, note that $\w_{g_0}(0)=\w_+(F)$. Indeed,
\[
    \w_+(F)=\w_+(\overline{q})+\w_+(g_0)
    =\sum_{k=1}^m\underset{=0}{\underbrace{\w_+(1-\zeta_k\overline{\zeta})}}\,+\,\w_{g_0}(0)=\w_{g_0}(0).
\]
So we finally obtain that $\dim\ker(T_F)=\max(0,-\w_+(F))$.
\end{proof}

Here is a direct consequence of 
\Cref{Ahern-Clark} on the possibility to embed a Toeplitz operator into a $C_0$-semigroup.

\begin{corollary}\label{un-corollaire}
Let $p>1$ and let $F$ be a differentiable function  on $\T$ with $F'\neq0$ on $\mathbb T$. If $\w_+(F)<0$, then $T_F$ is not embeddable into a $C_0$-semigroup on $H^p$.
\end{corollary}

In the rest of this section, we give a simple geometrical interpretation of the  number $\w_+$ (which is a reformulation of the interpretation given by Ahern and Clark) so as to be able to see quickly whether the condition of \Cref{un-corollaire} is satisfied or not.

\subsection{A geometrical interpretation of the  number $\w_+$ and consequences}
In \cite{AhernClark1985}, Ahern and Clark gave the following geometric interpretation of the number $\w_+$: \emph{suppose that $0\in F(\mathbb T)$ and let $\Omega$ be a component of $\mathbb C\setminus \sigma(T_F)$ with $0\in\partial\Omega$. We say that an arc $\gamma=F(\{e^{i\theta},\alpha<\theta<\beta\})$ is \emph{negative} if $0\in\gamma$ and, on a neighborhood of $0$, $\Omega$ remains on the right when we travel along $\gamma$. Let $K$ be the number of negative arcs that intersect only at the point $0$. Then $\w_+(F)=\w_F(\Omega)+K$.} 
In this geometrical interpretation, for two subarcs $\gamma,\gamma'$ of $F(\mathbb T)$ which coincide on a neighborhood of $0$, we just count  one of these subarcs. In other words, $K$ is the number of points $\zeta_0\in\mathbb T$ with $F(\zeta_0)=0$ such that for all sufficiently small $\varepsilon>0$, the curve $F(\{\zeta\in\mathbb T\,;\,|\zeta-\zeta_0|<\varepsilon\})$ is negative.

\begin{figure}[ht]
\begin{subfigure}{.48\linewidth}\centering
\begin{tikzpicture}
 \draw(0,0)node{\includegraphics[scale=.8]{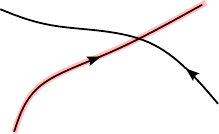}};
 \draw(0,-.4)node{$\Omega$};
\end{tikzpicture}
\caption{$K=1$}
\label{Fig2a}
\end{subfigure}
\begin{subfigure}{.48\linewidth}\centering
\begin{tikzpicture}
 \draw(0,0)node{\includegraphics[scale=.8]{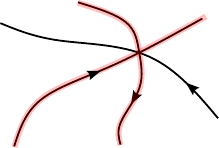}};
 \draw(-.3,-.4)node{$\Omega$};
\end{tikzpicture}
\caption{$K=2$}
\label{Fig2b}
\end{subfigure}
\caption{}
\end{figure}
In our situation, we can interpret $\w_+$ as the winding number in $0$ of a little perturbation of the curve obtained as follows: \emph{consider a very small open arc $\gamma \subset \mathbb T$ such that $0\in F(\gamma)$; then move a little bit the arc $\Gamma:=F(\gamma)$ away from $0$ so as to keep $0$ on the left when traveling on this modified arc $\widetilde{\Gamma}$. Then $\w_+$ is the winding number of the modified curve at the point $0$.}
\par\smallskip
Let $F$ satisfy the assumptions \ref{H1}, \ref{H2} and \ref{H3}, and suppose that $0$ belongs to $F(\mathbb T)\setminus\mathcal O$, i.e. $0$ belongs to the curve $F(\T)$ but is not a point of self-intersection. Let $\Omega_{int}$ and $\Omega_{ext}$ be the interior and exterior components at the point $0$ respectively.  Then the geometrical interpretation above can be illustrated as follows:

\begin{figure}[ht]
\begin{tikzpicture}
 \draw(0,0)node{\includegraphics[scale=1]{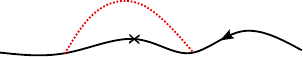}};
 \draw(1,.4)node{$\Omega_{int}$};
 \draw(1,-.7)node{$\Omega_{ext}$};
 \draw(-.3,-.5)node{$0$};
\end{tikzpicture}
    \caption{}
    \label{Fig2}
\end{figure}
This gives directly the following result:

\begin{proposition}\label{prop:calcul-de-w+}
Let $F$ satisfy the assumptions \ref{H1}, \ref{H2} and \ref{H3}, and suppose that $0$ belongs to $F(\mathbb T)\setminus\mathcal O$. Then
\[\w_+(F)=\w_F(\Omega_{ext}),\]
where $\Omega_{ext}$ is the exterior component at the point $0$.
\end{proposition}

Combining \Cref{Fact:RemOnPos0}, \Cref{un-corollaire} and \Cref{prop:calcul-de-w+}, we obtain the following necessary condition for the embeddability of $T_F$:

\begin{proposition}\label{avant le th C}
    Let $p>1$, and let $F$ satisfy the assumptions \ref{H1}, \ref{H2} and \ref{H3bis}. Suppose that $T_F$ is embeddable into a $C_0$-semigroup of bounded operators on $H^p$. Then either $0$ belongs to $\C\setminus \sigma(T_F)$ or, if $0$ belongs to the spectrum of $T_F$, then it belongs to $\partial\sigma(T_F)\cup\mathcal{O}$. In particular, if $0$ belongs to $\intsp$, then necessarily it belongs to $\mathcal{O}$.
\end{proposition}

\begin{proof}
We can suppose without loss of generality that $F$ satisfies \ref{H3}. Suppose that $T_F$ is embeddable. Then, according to \Cref{Fact:RemOnPos0}, either $0$ belongs to $\C\setminus\sigma(T_F)$, or $0$ belongs to $F(\T)$. 
Suppose that $0$ does not belong to $\mathcal{O}$.
By \Cref{un-corollaire}, $\w_+(F)\ge 0$, and then by \Cref{prop:calcul-de-w+}, 
$\w_F(\Omega_{ext})\ge 0$, where $\Omega_{ext}$ is the exterior component at the point $0$. Since $F$ satisfies \ref{H3}, we thus get that 
$\w_F(\Omega_{ext})=0$, and hence $0$ belongs to the boundary of the spectrum of $T_F$.
\end{proof}

As a direct consequence of this proposition, we obtain the converse implication in \Cref{Th:CSforEmb} under the additional assumptions that $\C\setminus\intsp$ is connected and that $0$ does not belong to $\mathcal{O}$:

\CondPourEquiv*

\begin{proof}
Since $\mathbb C\setminus\intsp$ is connected, the implication $(2)\Rightarrow (1)$ is a direct consequence of \Cref{Th:CSforEmb}. As to the implication $(1)\Rightarrow (2)$, it follows immediately from \Cref{avant le th C}.
\end{proof}

\par\smallskip
Consider now the case where $0$ is a simple intersection point of the curve $F(\T)$, and
let us describe what the geometrical interpretation gives in this case. In this situation, there are four possibilities, 
which we will call Type I, Type II, Type III and Type IV intersections; see the pictures below. 

On these pictures, the integer $L$ denotes the maximum of $\w_F(\Omega)$ taken over all the components $\Omega$ of $\sigma(T_F)\setminus F(\T)$ such that $0$ belongs to the boundary of $\Omega$.

\begin{figure}[ht]
\begin{subfigure}{.48\linewidth}\centering
\begin{tikzpicture}
 \draw(0,0)node{\includegraphics[scale=1.2]{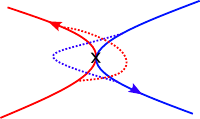}};
 \draw(-1.3,0)node{$L$};
  \draw(1,0)node{$L$};
  \draw(0,1)node{$L-1$};
    \draw(0,-1)node{$L-1$};
\end{tikzpicture}
 \caption{Type I}
 \label{Fig4a}
\end{subfigure}
\begin{subfigure}{.48\linewidth}\centering
\begin{tikzpicture}
 \draw(0,0)node{\includegraphics[scale=1.2]{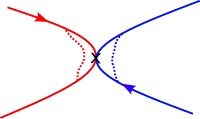}};
  \draw(-1,0)node{$L-1$};
  \draw(1,0)node{$L-1$};
  \draw(0,.6)node{$L$};
    \draw(0,-.6)node{$L$};
\end{tikzpicture}
 \caption{Type II}
 \label{Fig4b}
\end{subfigure}
\begin{subfigure}{.48\linewidth}\centering
\begin{tikzpicture}
 \draw(0,0)node{\includegraphics[scale=1.2]{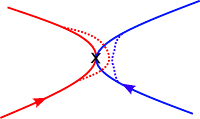}};
  \draw(-.6,0)node{$L$};
  \draw(1,0)node{$L-2$};
  \draw(0,1)node{$L-1$};
    \draw(0,-.8)node{$L-1$};
\end{tikzpicture}
 \caption{Type III}
 \label{Fig4c}
\end{subfigure}
\begin{subfigure}{.48\linewidth}\centering
\begin{tikzpicture}
 \draw(0,0)node{\includegraphics[scale=1.2]{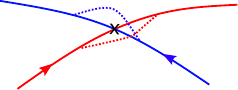}};
  \draw(-1,0.3)node{$L$};
  \draw(1.2,0.3)node{$L-2$};
  \draw(0,1)node{$L-1$};
    \draw(0,-.4)node{$L-1$};
\end{tikzpicture}
 \caption{Type IV}
 \label{Fig4d}
\end{subfigure}
    \caption{}
\end{figure}
When the intersection is of Type II, Type III or Type IV, we have
$\w_+(F)=L$. But when the intersection is of Type I, we have $\w_+(F)=L+1$. 
Suppose now that $F$ satisfies \ref{H1}, \ref{H2} and \ref{H3}, and that all the intersection points of the curve $F(\T)$ are simple.
If $0$ is such an intersection point, it follows from \Cref{Ahern-Clark} that $\ker(T_F)=\{0\}$ if and only if  $0\in\partial\sigma(T_F)$
or if the intersection in $0$ has Type I with $L=-1$. So, if we denote by $\Omega_j$ the set of elements $\lambda\in\mathbb C\setminus F(\mathbb T)$ such that $\w_F(\lambda)=-j$, we can deduce the following result:

\begin{corollary}\label{th:deuxiemeCNS}
Let $p>1$, and let $F$ satisfy \ref{H1}, \ref{H2} and \ref{H3}. Suppose that $\mathbb C\setminus\intsp $ is connected, that all the  intersections of the curve $F(\T)$ are simple, and that this curve has no intersection of Type I on $\partial\Omega_1\cap \partial \Omega_2$.
\par\smallskip
Then $T_{F}\in\mathcal{B}(H^p)$ is embeddable if and only if $0$ belongs to $\mathbb{C}\setminus\intsp $.
\end{corollary}

Note that an intersection of Type I can appear only in the case where the intersection is tangential. In other words, if, for example, $\mathbb C\setminus \intsp $  is connected and all the intersections of the curve $F(\T)$ are simple and transversal, then the assumptions of \Cref{th:deuxiemeCNS} are satisfied, i.e. $T_F$ is embeddable if and only if $0$ belongs to $\mathbb C\setminus \intsp $.
\par\smallskip
In the next section we provide some further, more technical, conditions under which the embeddability of the operator $T_F$ implies that $0$ belongs to the unbounded connected component of $\mathbb C\setminus \intsp $. Our approach goes via a study of the commutant of the model operator acting on the model space. 

\section{A condition on the commutant}\label{Section 6}
Let us recall that 
\(
\Omega_j=\{\lambda\in\mathbb C\setminus F(\mathbb T):|\w_F(\lambda)|=j\}\quad\mbox{for }0\leq j\leq N,
\)
and 
\[
\Omega_j^+=\{\lambda\in\mathbb C\setminus F(\mathbb T):|\w_F(\lambda)|>j\}=\bigcup_{k=j+1}^N\Omega_k\quad\mbox{for }0\leq j\leq N-1.
\]
\subsection{A description of the commutant of the model operator}\label{Section 6.2}
Given a boun\-ded operator $T$ on a separable Banach space $X$, we denote by $\{T\}'$ the commutant of $T$, that is the set
$\{T\}'=\{S\in \mathcal B(X)\,;\,TS=ST\}.$

\begin{lemma}\label{Fact:commutant1}
Let $A$ be a bounded operator on $E_F^q$. Then $A$ commutes with $M_\lambda$ if and only if  for every pair $(i,j)$ of integers with $0\le i,j<N$, there exists a function $a_{i,j}\in E^1(\Omega_{\max(i,j)}^+)$ such that 
\begin{equation}\label{Eq:commutant0}
(Au)_j=\sum_{i=0}^{k-1}a_{i,j}u_i~\text{ on}~\Omega_k
\end{equation}
for every $0\le j<k\le N$ and every $u=(u_0,\dots,u_{N-1})\in E^q_F$.
\end{lemma}

\begin{proof} If $A$ is a bounded operator on $E^q_F$, it is clear that if $A$ is given by \eqref{Eq:commutant0}, then $A$ commutes with $M_\lambda$. So we just need to prove the converse assertion. Let $A\in\{M_\lambda\}'$ 
and let $B$ be the bounded operator on $H^p$ defined by $B=(U^{-1}A U)^*$, where $U: H^q\longrightarrow E_F^q$ is the operator given by \eqref{defn-model-U-section2-5}. It follows from \Cref{T:Yakbovich_Hp} and from the fact that $A$ and $M_\lambda$ commute that
 \begin{align*}
B^*T_F^*~=~U^{-1}AUU^{-1}M_\lambda U~&=~U^{-1}AM_\lambda U\\
&=~U^{-1}M_\lambda AU~=~U^{-1}M_\lambda UU^{-1}AU~=~T_F^*B^*.
 \end{align*}
 Then $BT_F=(T_F^*B^*)^*=(B^*T_F^*)^*=T_FB$, i.e. $B\in \{T_F\}'$. 
\par\smallskip
 Since $B$ commutes with $T_F$, the eigenspaces of $T_F$ are invariant by $B$. But recall that 
 \[\ker(T_F-\lambda)=\spa\big[h_{\lambda,j}\,;\,0\le j<|\w_F(\lambda)|\big]~\text{for all}~\lambda\in\sigma(T_F)\setminus F(\mathbb T),\]
and thus for every $\lambda\in\sigma(T_F)\setminus F(\mathbb T)$ and every pair $(i,j)$ of integers with $0\le i,j<|\w_F(\lambda)|$, there exists a scalar $a_{i,j}(\lambda)\in \mathbb C$ such that 
 \[Bh_{\lambda,j}=\sum_{i=0}^{|\w_F(\lambda)|-1}a_{i,j}(\lambda) h_{\lambda,i}.\]
 Now, let $u=(u_0,\dots,u_{N-1})\in E^q_F$. Then $g:=U^{-1}u$ belongs to $H^q$ and $Au=UB^*g$. Thus, for every $\lambda\in\Omega_k$, $0\le j<k\le N$, we have 
\begin{align*}
    (Au)_j(\lambda)~=~(UB^*g)_j(\lambda)~&=~\ps{B^*g,h_{\lambda,j}}\\
    &=~\ps{g,Bh_{\lambda,j}}~=~\sum_{i=0}^{k-1}a_{i,j}(\lambda)\ps{g,h_{\lambda,i}}=\sum_{i=0}^{k-1}a_{i,j}(\lambda)u_j(\lambda),
\end{align*}
which proves \eqref{Eq:commutant0}. To finish the proof of \Cref{Fact:commutant1}, it remains to prove that each function $a_{i,j}$ belongs to $ E^1(\Omega_{\max(i,j)}^+)$.
\par\smallskip
Let us treat first the case $i=0$. Let $u=U1=(1,0,\dots,0)$, which belongs to $E^q_F$ by \Cref{elements-de-l-espace-modele}. Then, according to \eqref{Eq:commutant2},  we have $Au=(a_{0,0},\dots,a_{0,N-1})$, and since $Au$ is an element of $E^q_F$ each function $a_{0,j}$ belongs to $E^q(\Omega_j^+)$, hence to $E^1(\Omega_j^+).$
\par\smallskip
Now suppose that $0\le k<N$ is such that for every $0\le i\le k$ and every $0\le j<N$, the function $a_{i,j}$ is an element of $E^1(\Omega_{\max(i,j)}^+)$. Let $u:=Uz^{k+1}$, which we write as $u=(u_0,\dots,u_{k},1,0,\dots)$ with $u_l\in H^\infty(\Omega_l^+)$ for $0\le l\le k$ by \Cref{elements-de-l-espace-modele}, and let $v:=Au$, with $v=(v_0,v_1,\dots,v_{N-1})$. Then, for every $0\le j<N$ and every $l$ with $\max(j, k+1)<l<N$, we have
\[v_j=\sum_{i=0}^k a_{i,j}u_i~+~a_{k+1,j}~\text{ on}~\Omega_l.\]
But we already observed that $u_i\in H^\infty(\Omega_l)$, and we know that $a_{i,j}\in E^1(\Omega_l)$ by the induction hypothesis. So it follows that 
\[a_{k+1,j}~=~v_j-\sum_{i=0}^k a_{i,j}u_i~\text{ belongs to } E^1(\Omega_l)~\text{ for every $l>\max(k+1,j)$}.\]
Hence the induction assumption holds for $k+1$ too, and this proves \Cref{Fact:commutant1}.
\end{proof}
\begin{remark}
The proof of \Cref{Fact:commutant1} yields in fact that each function $a_{i,j}$ belongs to $E^q(\Omega_{\max{(i,j)}}^+)$, which is a stronger conclusion than what we stated in \Cref{Fact:commutant1}. However, for simplicity's sake we prefer to state a conclusion to \Cref{Fact:commutant1} which is independent of $q$. 
\end{remark}

We state separately the following consequence of \Cref{Fact:commutant1}, which will be used several times in \Cref{Section 7}:

\begin{lemma}\label{la-forme-de-A_t-sur-Omega_1}
Suppose that $M_\lambda$ is embeddable into a semigroup $(A_t)_{t>0}$ of operators on $E^q_F$. Then there exists an analytic branch $\log$ of the logarithm on $\Omega_1$ such that for every $u=(u_0,\dots,u_{N-1})\in E^q_F$ and every $t>0$, if $v_t=(v_{t,0},\dots,v_{t,N-1})=A_t u$, then $v_{t,0}=\alpha_t u_0$ on $\Omega_1$, with $\alpha_t(\lambda)=e^{t\log(\lambda)}$ for every $\lambda\in\Omega_1$. 
\end{lemma}

\begin{proof}
Since  the operator $A_t$ commutes with $M_\lambda=A_1$ for every $t>0$, \Cref{Fact:commutant1} implies that there exists $\alpha_t\in E^1(\Omega_0^+)$ such that for every  $u=(u_0,\dots,u_{N-1})\in E^q_F$, if $v_t=(v_{t,0},\dots,v_{t,N-1})=A_t u$, then $v_{t,0}=\alpha_t u_0$ on $\Omega_1$. Now, by \Cref{elements-de-l-espace-modele}, we know that $u=(1,0,\dots,0)\in E_F^q$. Then if $v=(v_0,\dots,v_{N-1})=A_{1/2}^2(u)$, we have $v_0(\lambda)=\alpha_{1/2}^2(\lambda)$ for every $\lambda\in\Omega_1$. On the other hand, using that $A_{1/2}^2=A_1=M_\lambda$, we also have $v_0(\lambda)=\lambda$, which gives that $\alpha_{1/2}^2(\lambda)=\lambda$ for every $\lambda\in\Omega_1$. In particular, $\alpha_{1/2}$ is an analytic determination of the square root on $\Omega_1$, which implies that $0\not\in\Omega_1$. 
\par\smallskip
Fix now  $\lambda\in\Omega_1$. Since the family $(\alpha_t(\lambda))_{t>0}$ is a scalar-valued semigroup and $\alpha_1(\lambda)=\lambda\neq 0$, we deduce that $\alpha_t(\lambda)\neq 0$ for every $t>0$, and there exists a complex number $c_\lambda$  such that $\alpha_t(\lambda)=e^{tc_\lambda}$ for every $t>0$. Observe also that, since $A_1=M_\lambda$, we have $e^{c_\lambda}=\lambda$. Since $\Omega_1$ is a finite union of simply connected domains and does not contain $0$, there exists an analytic branch $\varphi$ of the logarithm on $\Omega_1$. For each $\lambda\in\Omega_1$, there exists $k_\lambda\in\Z$ such that 
$c_\lambda=\varphi(\lambda)+2i\pi k_\lambda$. For every $t>0$, $\alpha_t(\lambda)=e^{t(\varphi(\lambda)+2i\pi k_\lambda)}$ on $\Omega_1$, and since the function $\alpha_t$ is analytic on $\Omega_1$, it follows that the map $\lambda\mapsto e^{2i\pi k_\lambda t}$ is analytic on $\Omega_1$. This being true for every $t>0$, the map $\lambda\mapsto k_\lambda$ is constant on each connected component of $\Omega_1$, whence it follows that the map $\lambda\mapsto c_\lambda$ is an analytic branch of the logarithm on $\Omega_1$ which satisfies the required properties.
\end{proof}

\subsection{Links with embeddability}\label{Subsection6.3}
If an operator $T$ is embeddable in a $C_0$-semi\-group $(T_t)_{t>0}$, then obviously $T_t$ commutes with $T$ for every $t>0$. Our motivation for the study of the commutant of Toeplitz operators is the following direct consequence of \Cref{Th:EmbModelMult} and \Cref{equivalent}. In its statement, $M_\lambda$ denotes as usual the multiplication operator by the independent variable $\lambda$ on the model space $E_F^q$. Recall also that by \Cref{multiplicateurs}, multipliers of $E_F^q$ are characterized in the following way: $\mathcal Mult(E_F^q)=\{g_{|\Omega_0^+}:g\in H^\infty(\intsp )\}$.

\begin{proposition}\label{prop:commutant}
Let $p>1$. Let $F$ be a symbol satisfying \ref{H1}, \ref{H2} and \ref{H3}. 
Suppose that the commutant of $M_\lambda$ consists of multiplication operators only, i.e. that $$\{M_\lambda\}'=\left\{M_g\,;\,g\in \mathcal Mult(E^q_F)\right\}.$$
Then $T_F$ is embeddable into a $C_0$-semigroup on $H^p$ if and only if $0$ belongs to the unbounded component of $\mathbb C\setminus \intsp $.
\end{proposition}

\begin{remark}
    For symbols $F$ satisfying \ref{H1}, \ref{H2} and \ref{H3}, knowing that the  only operators in the commutant of $M_\lambda$ are multiplication operators yields a neat description of the commutant of $T_F$ itself. Indeed, under assumptions \ref{H1}, \ref{H2} and \ref{H3} the operator $T_F$ admits an $H^\infty(\intsp )$ functional calculus \cite{Yakubovich1991}. This is a direct consequence of \Cref{T:Yakbovich_Hp}: for any function 
    $g\in H^\infty(\intsp )$, set $g(T_F):=(U^{-1}g(M_\lambda)U)^*=(U^{-1}M_gU)^*$: this defines a bounded functional calculus for $T_F$ on $ H^\infty(\intsp )$.  If the commutant of $M_\lambda$ consists of multiplication operators only, any operator $A\in\{T_F\}'$ can thus be written as $A=g(T_F)$ for some function $g\in H^\infty(\intsp )$.
\end{remark}

\par\smallskip
Hence it is a natural problem to try to determine whether the commutant of $M_\lambda$ coincides with the set of multiplication operators -- but it turns out to be a difficult one. Indeed, we provide below an example of a smooth curve admitting two different parametrizations $F_1$ and $F_2$ such that the commutant of $M_\lambda$ acting on $E^q_{F_1}$ consists of multiplication operators only, while the commutant of $M_\lambda$ acting on $E^q_{F_2}$ does not.

\begin{example}\label{exemples-commutant}
 Let $p>1$, and $\varepsilon>\max(1/p, 1/q)$. Let $F\in C^{1+\varepsilon}(\mathbb T)$ be any function such that $F'$ does not vanish on $\mathbb T$ and the curve 
 $F(\mathbb T)$ is given by \Cref{fig:cercles}. 
 
 \begin{figure}[ht]
\begin{tikzpicture}[scale=.9]
\draw[<-] (0,1)arc(90:450:1);
\draw(0,0)node{\small $\times$}(0,0)node[below]{$0$};
\draw(1,0)node{\small $\times$}(1,0)node[below right]{$1$};
\draw(-1,0)node{\small $\times$}(-1,0)node[below left]{$-1$};
\draw(.3,.4)node{$\Omega_2$};
\draw(-1.7,.5)node{$\Omega_1$};
\draw[<-] (-1,2)arc(90:450:2)node[above]{~};
\draw[-](0,-2.1)node{~};
\end{tikzpicture}
\caption{}
    \label{fig:cercles}
\end{figure}
Because of \Cref{Fact:commutant1}, we know that $A$ commutes with $M_\lambda$ if and only if $A$ belongs to $B(E^q_F)$ and for every $u=(u_0,u_1)\in E^q_F$, $(v_0,v_1)=Au$ satisfies
\begin{equation}\label{Eq:FormAcommut2}
v_0=\begin{cases}
        a_{0,0}u_0&\text{on }~\Omega_1\\
        a_{0,0}u_0+a_{1,0}u_1&\text{on }~\Omega_2
    \end{cases}\quad\text{and}\quad v_1=a_{0,1}u_0+a_{1,1}u_1~
\text{on }~\Omega_2,
\end{equation}
for some functions $a_{i,j}\in E^1(\Omega_{\max(i,j)}^+)$, $0\le i,j\le 1$. So suppose that $A\in\{M_\lambda\}'$, so that $A$ is given by \eqref{Eq:FormAcommut2} for some functions $a_{i,j}\in E^1(\Omega_{\max(i,j)}^+)$. Fix  $u=(u_0,u_1)\in E^q_F$ and write $Au$ as $Au=(v_0,v_1)$. Recall that $(v_0,v_1)\in E^q_F$ if and only if $v_0\in E^q(\Omega_0^+)$, $v_1\in E^q(\Omega_1^+)$ and we have the following boundary relation:
\[
v_0^{int}-\zeta v_1^{int}=v_0^{ext}~\text{a.e. on}~\partial\Omega_1^+=\partial\Omega_2=\mathbb T,
\]
where $\zeta=1/F^{-1}$. When written in terms of the functions $a_{i,j}$, this boundary condition becomes: for every $(u_0,u_1)\in E^q_F$,
\[
    a_{0,0}^{int}u_0^{int}+a_{1,0}^{int}u_1^{int}-\zeta (a_{0,1}^{int}u_0^{int}+a_{1,1}^{int}u_1^{int})=a_{0,0}^{ext}u_0^{ext}~\text{a.e. on }\mathbb T.
\]
Since $u_0^{int}-\zeta u_1^{int}=u_0^{ext} $ a.e. on $\mathbb T=\partial\Omega_2=\partial\Omega_1^+$, this relation is equivalent to
\begin{equation}\label{Eq:Ex2CerclesB1}
(a_{0,0}^{int}-\zeta a_{0,1}^{int}-a_{0,0}^{ext})u_0^{int}+(a_{1,0}^{int}-\zeta(a_{1,1}^{int}-a_{0,0}^{ext}))u_1^{int}=0.
\end{equation}
Using that $U1=(1,0)$ and $Uz=(u_0,1)$ belongs to $E^q_F$, we deduce that  $Au$ satisfies \eqref{Eq:Ex2CerclesB1} a.e. on $\partial \Omega_2$ for every $u\in E^q_F$ if and only if 
\[\begin{cases}
 a_{0,0}^{int}-\zeta a_{0,1}^{int}~=~a_{0,0}^{ext} \\
 a_{1,0}^{int}-\zeta(a_{1,1}^{int}-a_{0,0}^{ext})~=~0
\end{cases}~\text{a.e. on }\mathbb T\]
and this is equivalent to
\begin{equation}\label{Eq:CNS-commut}
\begin{cases}
 a_{0,0}^{int}-\zeta a_{0,1}^{int}~=~a_{0,0}^{ext} \\
 a_{1,0}^{int}-\zeta(a_{1,1}^{int}-a_{0,0}^{int})-\zeta^2a_{0,1}^{int}~=~0
\end{cases}~\text{a.e. on }\mathbb T.
\end{equation}
\par\smallskip
We now consider two different examples of such a function $F$ for which the commutant of $M_\lambda$ acting on $E^q_F$ has different descriptions.
\par\smallskip
\noindent\textbf{Parametrization 1:} Let us first consider the function $F_1$ of the example of Section 4.1 in \cite{FricainGrivauxOstermann_preprint}, which is defined as
\begin{equation}\label{Eq:Ex2CerclesF}
    F_1(e^{i\theta})~=~\begin{cases}
-1+2e^{-i3\theta/2}&\text{if}~0\le \theta<4\pi/3\\
e^{-3i\theta}&\text{if}~4\pi/3\le \theta<2\pi.
\end{cases}
\end{equation}
For this choice of $F_1$,  an operator $A\in\mathcal{B}(E^q_{F_1})$ commuting with $M_\lambda$ \emph{has to be a multiplication operator}, i.e. $a_{1,0}=a_{0,1}=0$ and $a_{0,0}=a_{1,1}$ on $\pt\Omega_2$. According to \eqref{Eq:CNS-commut}, in order to prove this, it is sufficient to check that for all functions $u,v,w\in E^1(\Omega_2)=H^1(\mathbb D)$
\[
u+2\zeta_1 v+\zeta_1^2w=0~\text{a.e. on}~\mathbb T=\partial\Omega_2\implies u=v=w=0~\text{on}~\mathbb D=\Omega_2,
\]
where $\zeta_1:=1/F_1^{-1}$ is defined on $F_1(\mathbb T)\setminus\{1\}$.  So let $u,v,w\in H^1(\mathbb D)$ be such that $u+2\zeta_1 v+\zeta_1^2w=0~\text{a.e. on}~\mathbb T$. We first prove that $w\equiv 0$ on $\mathbb D$. Argue by contradiction, and assume that $w$ is not identically $0$ on $\mathbb D$.
Since $\zeta_1=1/F_1^{-1}$ satisfies
\[\zeta_1(\lambda)~=~\exp\left[\frac {i}{3}\arg_{(0,2\pi)}(\lambda)\right]\quad\text{for every }\lambda\in\partial\Omega_2\setminus\{1\},\]
it admits a bounded analytic extension to $\Omega_2\setminus[0,1]$, given by
\[
\zeta_1(\lambda)~=~\exp\left[\frac 13\log_{(0,2\pi)}(\lambda)\right]\quad\text{for every }\lambda\in\Omega_2\setminus[0,1].
\]
Thus $u+2\zeta_1 v+\zeta_1^2w$ belongs to $E^1(\mathbb D\setminus[0,1])$, and hence the condition $u+2\zeta_1 v+\zeta_1^2w=0~\text{a.e. on}~\mathbb T$ implies that $u+2\zeta_1 v+\zeta_1^2w=0$ on $\mathbb D\setminus[0,1]$. Since $w\not\equiv0$, we have that
\begin{align*}
u+2\zeta_1 v+\zeta_1^2w=0
&\textrm{ if and only if } (w\zeta_1+v)^2-v^2+uw=w(u+2\zeta_1 v+\zeta_1^2w)=0\\
&\textrm{ if and only if }(w\zeta_1+v)^2=v^2-uw.
\end{align*}
In particular it follows that $(w\zeta_1+v)^2$ is analytic and thus continuous on $\mathbb D$. For every $x\in(0,1)$, we have that 
\[\zeta_1(x^+):=\lim_{y\to 0^+}\zeta_1(x+iy)=\sqrt[3]{x}\quad\text{and}\quad\zeta_1(x^-):=\lim_{y\to 0^-}\zeta_1(x+iy)=\sqrt[3]{x}e^{2i\pi/3}.\]
But the continuity of $(w\zeta_1+v)^2$ at the point $x$ gives that $(w(x)\zeta_1(x^+)+v(x))^2=(w(x)\zeta_1(x^-)+v(x))^2$, which means that 
\begin{equation}\label{eq:232ZDSSD}w(x)\zeta_1(x^+)+v(x)=\pm(w(x)\zeta_1(x^-)+v(x))\quad\mbox{for every }x\in (0,1).
\end{equation}
Now, by the uniqueness principle, we can find $0<a<b<1$ such that $w(x)\neq 0$ for every $x\in (a,b)$. Since $\zeta_1(x^+)\neq\zeta_1(x^-)$, we easily see that \eqref{eq:232ZDSSD} necessarily implies that 
for every $x\in(a,b)$, we have 
\[w(x)\zeta_1(x^+)+v(x)=-(w(x)\zeta_1(x^-)+v(x)),\]
and thus
\begin{equation}\label{Eq:Ex2cercles3}
(1+e^{2i\pi/3})w(x)\sqrt[3]{x}=w(x)(\zeta_1(x^+)+\zeta_1(x^-))=-2v(x).
\end{equation}
Let $\widetilde v=-2v$ and $\widetilde w=(1+e^{2i\pi/3})w$. Then \eqref{Eq:Ex2cercles3} yields that $\widetilde v(z)^3=z\widetilde w(z)^3$ for every $z\in(a,b)$, and thus this equality holds also for every $z\in\mathbb D$ by the uniqueness principle. Since $w\not\equiv0$, $0$ is a zero with finite order of $\widetilde{w}$, and thus of $\widetilde{v}$ too. Denote by $n_1$ (respectively $n_2$) the order of multiplicity of this zero for $\widetilde{v}$ (respectively for $\widetilde{w}$). Then the equation $\widetilde v^3=z\widetilde w^3$ gives in particular that $3n_1=1+3n_2$, and this is the desired contradiction. Hence $w=0$ on $\mathbb D$ and $u+2\zeta v=0$ a.e. on $\mathbb T$. It then follows from Proposition 4.9 in \cite{FricainGrivauxOstermann_preprint} that $u=v=0$.
\par\smallskip
Therefore we deduce that if the parametrization of the curve in \Cref{fig:cercles} is given by \eqref{Eq:Ex2CerclesF}, then 
$$\{M_\lambda\}'=\left\{M_g\,;\,g\in \mathcal Mult(E^q_{F_1})\right\}.$$
\par\medskip
\noindent\textbf{Parametrization 2:} 
We now give an another parametrization $F_2$ of the curve in \Cref{fig:cercles} for which there exists an operator $A\in\{M_\lambda\}'$, $A\in\mathcal{B}(E^q_{F_2})$, which is not a multiplication operator.  
\par\smallskip
Let $F_2\in C^{1+\varepsilon}(\mathbb T)$ be such that 
 $F_2(\mathbb T)$ is given by \Cref{fig:cercles} and satisfies
 \[F_2(e^{i\theta})=e^{-2i\theta}~\text{ for every}~-\pi<\theta<0.\]
Then $\zeta_2=1/F_2^{-1}$ satisfies 
    \[\zeta_2(\lambda)=\exp\left[\frac i2\arg_{(0,2\pi)}(\lambda)\right],~\text{for every}~\lambda\in\partial\Omega_2\setminus\{1\}\]
    and thus $\zeta_2$ has an analytic extension to $\Omega_2\setminus[0,1]$ which is given by
    \[\zeta_2(\lambda)=\exp\left[\frac 12\log_{(0,2\pi)}(\lambda)\right],~\text{for every}~\lambda\in\Omega_2\setminus[0,1].\]
    Note that this analytic extension satisfies $\zeta_2(\lambda)^2=\lambda$ for every $\lambda\in\Omega_2\setminus[0,1]$. So if we consider the operator $A$ given by \eqref{Eq:FormAcommut2} with
    \[a_{0,0}(\lambda)=\exp\left[\frac 12\log_{(0,2\pi)}(\lambda)\right]~\text{on}~\Omega_1\]
    and
    \[a_{0,0}=a_{1,1}=0,~a_{0,1}=-1~\text{and}~a_{1,0}(\lambda)=-\lambda~\text{on}~\Omega_2,\]
   by \eqref{Eq:CNS-commut}, we obtain an operator $A\in\mathcal B(E^q_{F_2})$ which commutes with $M_\lambda$. Since $(AU1)_1=-1\neq0$, we deduce that \[A\notin\{M_a\,;\,a\in \mathcal Mult(E^q_{F_2})\},\]
    and thus the commutant of $M_\lambda$ on $E^q_{F_2}$ contains operators which are not multiplication operators.
\end{example}
These two examples show that it might be difficult to find conditions of a geometric nature on the curve $F(\T)$ implying that the commutant of $T_F$ consists of multiplication operators. We thus finish this section by presenting an \emph{analytic} condition, inspired by \Cref{exemples-commutant}, implying that the assumptions of \Cref{prop:commutant} are satisfied. Under this condition, $T_F$ is embeddable into a $C_0$-semigroup on $H^p$ if and only if $0$ belongs to the unbounded component of $\mathbb C\setminus \intsp $.

\begin{theorem}\label{th:commutant}
Let $p>1$, and let $F$ be a symbol satisfying the assumptions \ref{H1}, \ref{H2} and \ref{H3}. Suppose that for every connected component $\Omega$ of $\sigma(T_F)\setminus F(\mathbb T)$, the following condition holds:
\begin{enumerate}
    \item\label{Hyp:commut1} if $|\w_F(\Omega)|=2$, then for any $u,v,w\in E^1(\Omega)$ we have
    \[u+\zeta v+\zeta^2w=0~\text{a.e. on}~\partial\Omega\cap\partial\Omega_1\implies u=v=w=0~\text{on}~\Omega.\]
    \item\label{Hyp:commut2} if $k:=|\w_F(\Omega)|>2$, for any $u,v\in E^1(\Omega)$ we have
    \[u+\zeta v=0~\text{a.e. on}~\partial\Omega\cap\partial\Omega_{k-1}\implies u=v=0~\text{on}~\Omega.\]
\end{enumerate}
Then $\{M_\lambda\}'=\{M_a;\,a\in\mathcal Mult(E^q_F)\}$, i.e. the commutant of $M_\lambda$ consists of multiplication operators only. 
\par\smallskip
As a consequence, $T_F$ is embeddable into a $C_0$-semigroup on $H^p$ if and only if $0$ belongs to the unbounded component of $\mathbb C\setminus \intsp $.
\end{theorem}

Condition (\ref{Hyp:commut2}) of \Cref{th:commutant} has appeared already in \cite{FricainGrivauxOstermann_preprint}, where a geometric condition (called Property (P)) 
implying that (\ref{Hyp:commut2}) holds 
was given. As already mentioned above, we do not know any geometric condition implying that the analytic condition  (\ref{Hyp:commut1}) of \Cref{th:commutant} holds.

\begin{proof}
The fact that all multiplication operators belong to the commutant of $M_\lambda$ is clear. Consider now $A\in\{M_\lambda\}'$, and let us prove that there exists $a\in \mathcal Mult(E^q_F)$ such that $A=M_a$.  By \Cref{Fact:commutant1}, for every  pair $(i,j)$ of integers with $0\le i,j<N$, there exists a function $a_{i,j}\in E^1(\Omega_{\max(i,j)}^+)$ such that for every $\lambda\in\sigma(T_F)\setminus F(\mathbb T)$ and every $u\in E^q_F$,
\begin{equation}\label{Eq:commutant2}
(Au)_j(\lambda)=\sum_{i=0}^{|\w_F(\lambda)|-1}a_{i,j}(\lambda)u_i(\lambda).
\end{equation}
We now apply the hypothesis of \Cref{th:commutant} to certain linear combinations of the functions $a_{i,j}$ to prove by induction on $1\le k<N$ that $a_{0,1}=\dots=a_{0,k}=0$ on $\Omega_{k+1}$.
\par\smallskip
To this aim, let us first consider $u=U1$. Then $v=Au=(a_{0,0},\dots,a_{0,N-1})\in E^q_F$ and thus the functions $a_{0,j}$, $0\le j<N-1$ satisfy
\begin{equation}\label{Eq:commutant3}
    a_{0,j}^{int}-\zeta a_{0,j+1}^{int}=a_{0,j}^{ext}~\text{a.e. on}~\partial\Omega_{j+1}^+.
\end{equation}
We consider next $u=Uz$, written as $u=(u_0,1,0,\dots,0)$ by \Cref{elements-de-l-espace-modele}. Then $u_0^{int}-\zeta=u_0^{ext}$ a.e. on $\partial\Omega_1^+$ and if we set $v=Au$ and write $v=(v_0,v_1,\dots,v_{N-1})$, then 
\[\begin{cases}
    v_0~=~a_{0,0}u_0&\text{on}~\Omega_1\\
    v_j~=~a_{0,j}u_0+a_{1,j}&\text{on}~\Omega_{\max(1,j)}^+ \text{ for every } 0\le j<N.\\
\end{cases}\]
Since $v_0^{int}-\zeta v_1^{int}=v_0^{ext}$ a.e. on $\partial\Omega_1^+$, we have, in particular, a.e. on $\partial\Omega_1\cap\partial\Omega_2$, the following relation:
\[a_{0,0}^{int}u_0^{int}+a_{1,0}^{int}-\zeta (a_{0,1}^{int}u_0^{int}+a_{1,1}^{int})=a_{0,0}^{ext}u_0^{ext}.\]
Since $u_0^{int}-\zeta=u_0^{ext}$ a.e. on $\partial\Omega_1^+$, in particular a.e. on $\partial\Omega_1\cap\partial\Omega_2$,  we deduce that we have, again a.e. on $\partial\Omega_1\cap\partial\Omega_2$,
\[(a_{0,0}^{int}-\zeta a_{0,1}^{int}-a_{0,0}^{ext})u_0^{int}+a_{1,0}^{int}-\zeta(a_{1,1}^{int}-a_{0,0}^{ext})=0.\]
By \eqref{Eq:commutant3}, we know that $a_{0,0}^{ext}=a_{0,0}^{int}-\zeta a_{0,1}^{int}$ a.e. on $\partial\Omega_1^+$, and thus we deduce that 
\[a_{1,0}^{int}-\zeta(a_{1,1}^{int}-a_{0,0}^{int})-\zeta^2a_{0,1}^{int}=0~\text{a.e. on }\partial\Omega_1\cap\partial\Omega_2.\]
So, by the assumption (\ref{Hyp:commut1}) of \Cref{th:commutant} applied to any connected component of $\Omega_2$, we obtain that $a_{0,1}=0$ on $\Omega_2$.
\par\smallskip
Suppose now that $1\le k<N-1$ is such that $a_{0,1}=\dots=a_{0,k}=0$ on $\Omega_{k+1}$. By \eqref{Eq:commutant3}, we have for every $1\le j\le k$ the equality
\[a_{0,j}^{int}-\zeta a_{0,j+1}^{int}=a_{0,j}^{ext}=0~\text{a.e. on}~\partial\Omega_{k+1}\cap\partial\Omega_{k+2}.\]
So, this time by the assumption (\ref{Hyp:commut2}) of \Cref{th:commutant} applied to any connected component of $\Omega_{k+2}$, we deduce that $a_{0,j}=a_{0,j+1}=0$ on $\Omega_{k+2}$, i.e. that $a_{0,1}=\dots=a_{0,k+1}=0$ on $\Omega_{k+2}$. We have thus shown by induction that 
\begin{equation}\label{Eq:commutant4}
    a_{0,j}=0~\text{on}~\Omega_j^+~\text{ for every}~1\le j<N.
\end{equation}

We now need the following lemma:

\begin{lemma}\label{Fact:commutant3}
Let $a:=a_{0,0}$, which is defined on $\Omega_0^+$. Then $a$  belongs to $\mathcal Mult(E_F^q).$
\end{lemma}

\begin{proof}
Note that, by \Cref{StickingLemma}, the $N$-tuple $(u_0,0,\dots,0)$ belongs to $E^q_F$ if and only if $u_0$ lies in $ E^q\left(\intsp \right)$. 
Using standard arguments of the theory of multipliers, it is not difficult to prove that 
\[
\mathcal Mult\left(E^q\left(\intsp \right)\right) = H^\infty\left(\intsp \right).
\]
Hence, by \Cref{multiplicateurs}, we have that
\[\mathcal Mult(E_F^q)=\mathcal Mult\left(E^q\left(\intsp \right)\right).\]
So let $u_0\in E^q\left(\intsp \right)$, so that $u=(u_0,0,\dots,0)\in E^q_F$. According to \eqref{Eq:commutant2}, for every $0\leq j<N$ we have 
\[
(Au)_j=a_{0,j}u_0\quad\text{ on }\Omega_j^+,
\]
and using \eqref{Eq:commutant4}, we deduce that 
$Au=(au_0,0,\dots,0)$. Since $Au\in E^q_F$, we get that $au_0\in E^q\left(\intsp \right)$. So we deduce that 
\[a\in \mathcal Mult\left(E^q\left(\intsp \right)\right)=\mathcal Mult(E_F^q).\qedhere\]
\end{proof}

Now, since $a\in\mathcal Mult(E_F^q)$, if we consider $u\in E^q_F$ and then $v=Au-M_au$ written as $v=(v_0,\dots,v_{N-1})$, we have that $v\in E^q_F$ and by \eqref{Eq:commutant2} and \eqref{Eq:commutant4}, $v_0=0$ on $\Omega_1$.  To finish the proof of \Cref{th:commutant}, we need to prove that $v=0$, and this is a direct consequence of the following lemma.

\begin{lemma}\label{Fact:commutant4}
Let $F$ be a symbol satisfying \ref{H1}, \ref{H2} and \ref{H3}. Suppose that for every connected component $\Omega$ of $\sigma(T_F)\setminus F(\mathbb T)$ with $k:=|\w_F(\Omega)|\ge2$, we have that for all $u,v\in E^1(\Omega)$,
\begin{equation}\label{Eq:ImplicLemma6.9}
u+\zeta v=0~\text{a.e. on}~\partial\Omega\cap\partial\Omega_{k-1}\implies u=v=0~\text{on}~\Omega.
\end{equation}
Then the following property holds for every $u=(u_0,\dots,u_{N-1})\in E^q_F$: if $u_0=0$ on $\Omega_1$, then $u=0$.
\end{lemma}

\begin{proof}
Let us prove by induction on $1\le l\le N$ that $u_0=\dots=u_{l-1}=0$ on $\Omega_l$. 
\par\smallskip
For $l=1$, this is true since $u_0$ is supposed to vanish on $\Omega_1$. Suppose now that  $1\le l< N$ is such that that $u_0=\dots=u_{l-1}=0$ on $\Omega_l$. Then for every $0\le j<l$, we have, a.e. on $\partial\Omega_l\cap\partial\Omega_{l+1}$,
\[u_j^{int}-\zeta u_{j+1}^{int}=u_j^{ext}=0.\]
Then, by the assumption of \Cref{Fact:commutant4} applied to every connected component of $\Omega_{l+1}$, we deduce that $u_j=u_{j+1}=0$ on $\Omega_{l+1}$, i.e. $u_0=\dots=u_l=0$ on $\Omega_{l+1}$. This concludes the proof of \Cref{Fact:commutant4}.
\end{proof}
To finish the proof of \Cref{th:commutant}, remark first that the assumptions of this theorem imply the hypothesis of \Cref{Fact:commutant4}. Indeed, the implication \eqref{Eq:ImplicLemma6.9}  follows from Assumption (1) with $w=0$ if $k=2$, and from Assumption (2) if $k>2$.

Let now $u\in E^q_F$ and $v=Au-M_au$, with $v=(v_0,\dots,v_{N-1})$. As mentioned above, $v\in E^q_F$ and $v_0=0$ on $\Omega_1$. 
By \Cref{Fact:commutant4}, we deduce that $v=0$, and thus $Au=M_au$. So $A$ is a multiplication operator and \Cref{th:commutant} is proved.
\end{proof}

\section{Looking for a characterization of embeddability}\label{Section 7}

We begin this section by presenting in \Cref{un-th-supplementaire} below a 
condition on the function $\zeta$ which ensures that the embeddability of $T_F$ is equivalent to the fact that $0$ belongs to the unbounded component of $\mathbb C\setminus \intsp $. Working on a concrete example where $\max_{\lambda\in\C\setminus F(\T)}|\w_F(\lambda)|=2$ (\Cref{l'exemple!!}), we will see in \Cref{l'exemple}
that when this condition is violated, the operator $T_F$ may be embeddable even though $0$ belongs to a bounded component of $\mathbb C\setminus \intsp $. This example leads us to a characterization of the embeddability of $T_F$ for an interesting class of symbols $F$ such that 
$\max_{\lambda\in\C\setminus F(\T)}|\w_F(\lambda)|=2$, which we present in \Cref{la-caracterisation} (\Cref{la-caracterisation!!}). The proof of \Cref{la-caracterisation!!}
relies in part on the methods used to deal with \Cref{l'exemple!!}.

\subsection{A condition on the function $\zeta$}\label{le-th-6.1}
Our aim in this section is to prove the following result:

\begin{theorem}\label{un-th-supplementaire}
Let $p>1$, and suppose that $F\in L^{\infty}(\T)$ satisfies \ref{H1}, \ref{H2} and \ref{H3}.
Let us denote by $X$  the unbounded component of $\mathbb C\setminus \intsp $. Suppose that
\begin{enumerate}
    \item all the intersection points of the curve $F(\T)$ on $\partial X$ are simple;\par\smallskip
    \item for every connected  component $\Omega$ of $\sigma(T_F)\setminus F(\T)$ such that $\Omega\subset\Omega_2$ and $\partial\Omega\cap\partial X\neq\varnothing$, there exists a non trivial connected subarc $\gamma\subset(\partial\Omega\cap\partial\Omega_1)\setminus\mathcal O$ such that $\zeta$ does not coincide a.e. on $\gamma$ with the non-tangential limit of a function in $\mathcal N(\Omega)$.
\end{enumerate}
Then $T_F\in\mathcal{B}(H^p)$ is embeddable if and only if $0$ belongs to $X$.
\end{theorem}

Before we start the proof of \Cref{un-th-supplementaire}, recall that $F$ is a bijective map from $\mathbb T\setminus F^{-1}(\cal O)$ onto $F(\mathbb T)\setminus\cal O$, and that the map $\zeta=1/F^{-1}$ is well-defined on $F(\mathbb T)\setminus\cal O$ and of class $C^1$ on each open arc contained in $F(\mathbb T)\setminus\cal O$.

\begin{proof}
We already know by \Cref{Th:CSforEmb} that $T_F$ is embeddable into a $C_0$-semigroup on $H^p$ as soon as $0$ belongs to $X$. So let us prove the converse assertion, and assume that $T_F$ is embeddable into a $C_0$-semigroup on $H^p$. 
By \Cref{avant le th C}, we know that if $0$ belongs to the spectrum of $T_F$, then necessarily $0\in\mathcal O\cup\partial\sigma(T_F)$.
\par\smallskip
According to \Cref{equivalent}, there exists an operator $A\in \mathcal B(E^q_F)$ such that $A^2=M_\lambda$. Consider the $N$-tuple $(u_0,\dots,u_{N-1})=A(1,0\dots,0)$ of $E^q_F$. Our goal  is to prove that the function $u_0$ is an analytic determination of the square root on the open set $\mathcal V$, where $\mathcal V$ is defined as the interior of the union of the closures of all the connected components $\Omega$ of $\sigma(T_F)\setminus F(\mathbb T)$ which satisfy $\partial\Omega\cap\partial X\neq\varnothing$. Once this is proved, this will imply that $0$ has to belong to the unbounded connected component of $\mathbb C\setminus \mathcal V$, which is $X$.
\par\smallskip

Consider the operator $B=(U^{-1}AU)^*$ which acts on $H^p$, where $U$ is the operator from $H^q$ onto $E_F^q$ given by \eqref{defn-model-U-section2-5}. By the construction of the operator $A$ and \Cref{T:Yakbovich_Hp}, we have $B^2=T_F$ and hence for every $\lambda\in\mathbb C$, we have $B(\ker(T_F-\lambda))\subseteq\ker(T_F-\lambda)$. Now according to \eqref{eq:eigenvector-space}, for every $\lambda\in \Omega_1$, $\ker(T_F-\lambda)=\spa\big[h_{\lambda,0}\big]$. In particular, for every $\lambda\in\Omega_1$, there exists a complex number $\alpha(\lambda)$ such that $Bh_{\lambda,0}=\alpha(\lambda)h_{\lambda,0}$ and $\alpha(\lambda)^2=\lambda$.
By \Cref{elements-de-l-espace-modele}, $U1=(1,0,\dots,0)$, which implies that $(u_0,\dots,u_{N-1})=UB^*1$. In particular, for every $\lambda\in \Omega_1$, we get
$$
u_0(\lambda)= \ps{B^*1,h_{\lambda, 0}}=\ps{1,Bh_{\lambda, 0}}=\alpha(\lambda)\ps{1,h_{\lambda, 0}}=\alpha(\lambda).
$$
It follows that
${u_0}_{|\Omega_1}=\alpha$ is an analytic determination of the square root on $\Omega_1$. 
Our aim being to prove that $u_0$ is an analytic determination of the square root on $\mathcal V$, we are going to show that $u_1\equiv0$ on $\Omega_2\cap \mathcal V$, i.e. on every component $\Omega$ with $\w_F(\Omega)=-2$ and $\partial\Omega\cap\partial X\not=\varnothing$. This will give that for any two components $\Omega,\widetilde\Omega$ of $\mathcal V$ such that $\partial\Omega\cap\partial\widetilde\Omega$ has a positive measure, the boundary condition becomes $u_0^{int}=u_0^{ext}$ a.e. on $\partial\Omega\cap\partial\widetilde\Omega$ (remember that all the intersection points of $F(\T)$ on $\partial X$ are simple, so that a connected component $\Omega$ of $\sigma(T_F)\setminus F(\T)$ with $\partial\Omega\cap\partial X\not=\varnothing$ is necessarily such that $\w_F(\Omega)\ge -2$). Finally, \Cref{StickingLemma} will give that $u_0$ has an analytic extension to $\mathcal V$, and this extension has to be an analytic determination of the square root by the uniqueness theorem.
\par\smallskip
Suppose that there exists a connected component $\Omega$ of $\sigma(T_F)\setminus F(\mathbb T)$ satisfying  $\w_F(\Omega)=-2$ and $\partial\Omega\cap\partial X\not=\varnothing$, and such that ${u_1}$ is not identically zero on $\Omega$. Let $\gamma$ be  as in the hypothesis of \Cref{un-th-supplementaire}. Then, since $0\in\mathcal O\cup\partial\sigma(T_F)$, we have $0\notin\gamma$, and $\gamma\cap\mathcal O=\varnothing$. The fact that $0$ does not belong to $\Omega$ implies that there exists an analytic determination of the square root $\alpha_\gamma$ on $\Omega$ such that $\alpha_\gamma=\alpha$ on $\gamma$. Note that $\alpha_\gamma\in H^\infty(\Omega)\subset E^1(\Omega)$. By the boundary condition involved in the definition of $E^q_F$ we have that
\[u_0^{int}-\zeta u_1^{int}=u_0^{ext}=\alpha_\gamma~\text{a.e. on }\gamma.\]
Then $\zeta$ coincides a.e. on $\gamma$ with the non-tangential limit of the meromorphic function $w$ defined on $\Omega$ by
 \begin{equation}\label{Eq:ExtMero}
     w=\frac{u_0-\alpha_\gamma}{u_1}.
 \end{equation}
Since the three functions $u_0,u_1$ and $\alpha_\gamma$ belong to $E^1(\Omega)$, this yields a contradiction with the hypothesis (2), and thus such a component $\Omega$ does not exist.
\par\smallskip
Hence $u_1(\lambda)=0$ for every $\lambda\in \Omega_2\cap\mathcal V$, and then the boundary conditions yield that 
\[u_0^{int}=u_0^{ext}~\text{a.e. on }\mathcal V\cap F(\mathbb T).\]
Then, by \Cref{StickingLemma}, the function $u_0$ belongs to $E^q(\mathcal V)$. But $\mathcal V$ is connected, and by construction $u_0^2(\lambda)=\lambda$ for every $\lambda\in\mathcal V\cap \Omega_1$. By the uniqueness theorem, $u_0^2(\lambda)=\lambda$ for every $\lambda\in\mathcal V$ and this yields an analytic determination of the square root on $\mathcal V$. Thus we finally obtain that $0$ belongs to the unbounded component of $\mathbb C\setminus \mathcal V$, which is $X$. This concludes the proof of \Cref{un-th-supplementaire}.
\end{proof}
\begin{remark}\label{Une-remarque}
Note that the proof of \Cref{un-th-supplementaire} shows that if there exists a bounded operator $B$ on $H^p$ such that $B^2=T_F$ and if $\Omega$ is a connected component of $\sigma(T_F)\setminus F(\mathbb T)$ with $\w_F(\Omega)=-2$  then only one of these  two situations could occur:
\begin{itemize}
    \item[(a)] $u_1\equiv0$ on $\Omega$ where $(u_0,\dots,u_{N-1})=UB^*1$ and, in particular there exists an analytic determination of the square root on $\mathrm{int}(\overline{\Omega_1\cup\Omega})$;
    \item[(b)] for every connected arc $\gamma\subset(\partial\Omega\cap\partial\Omega_1)\setminus\mathcal O$,  the function $\zeta$ coincides a.e. on $\gamma$ with the non-tangential limit of the function in $\mathcal N(\Omega)$ given by \eqref{Eq:ExtMero}.
\end{itemize}
\end{remark}

\subsection{An example}\label{l'exemple}
In view of \Cref{un-th-supplementaire}, it is natural to try to understand better the
condition (2): we have seen that when it holds, the embeddability of $T_F$ forces $0$ to belong to $X$, the unbounded component of $\C\setminus \intsp$. If $T_F$ is embeddable, although $0$ does not belong to $X$, is it essentially because 
condition (2) is violated? In this section, we work out a concrete example which points towards such a result. This approach will be developed further in \Cref{la-caracterisation} below, using the intuition from the example as well as some results proved in the particular setting of \Cref{l'exemple!!}, but which hold in greater generality.

\begin{example}\label{l'exemple!!}
Suppose that $F$ satisfies \ref{H1}, and that the curve $F(\mathbb T)$ is given by the following figure:

\begin{figure}[ht]
\begin{tikzpicture}[scale=.8]
\draw(0,0)node{\includegraphics[scale=.8]{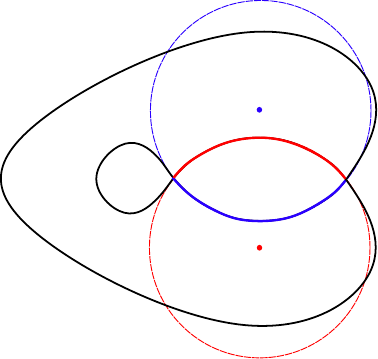}};
\draw[green!50!black](-.6,0)node{\small $\times$}(-.6,0)node[below]{\small$0$};
\draw(-1,.2)node{$\Omega_0$};
\draw(1,0)node{$\Omega_2$};
\draw(-1.7,1)node{$\Omega_1$};
\draw(-2,2)node{$\Omega_\infty$};
\draw[blue](1.3,1.0)node{\small $c_2$};
\draw[red](1.3,-1.0)node{\small $c_1$};
\draw[red](2.4,.5)node{\small$\gamma_1$};
\draw[blue](2.4,-.5)node{\small$\gamma_2$};
\draw[blue,dashed](1.2,1.2)--(1.2,3.2);
\draw[red,dashed](1.2,-1.2)--(1.2,-3.2);
\end{tikzpicture}
    \caption{}
    \label{fig:l'exemple!!}
\end{figure}
The curve $F(\T)$ is negatively oriented. The set $\C\setminus F(\T)$ has four connected components: let $\Omega_1$ and $\Omega_2$ denote the set of $\lambda\in\C$ such that $\w_F(\lambda)=-1$ and $\w_F(\lambda)=-2$, respectively. Let $\Omega_0$ be the bounded component of $\C\setminus \sigma(T_F)$, and $\Omega_{\infty}$ the unbounded component. 
\par\smallskip
For $j=1,2$, the curve $\gamma_j$ is a negatively oriented circle arc of radius $r$ centered at a point $c_j$, where $c_1=\overline{c}_2$ and $-r<\Im(c_1)<0$ (and thus $0<\Im(c_2)<r$). The arc $\gamma_1$ can be written as
$$\gamma_1=\{c_1+re^{i\tau}\,;\,\tau\in (\delta, \pi-\delta)\},$$ where $\delta=\arcsin(\Im(c_2)/r)$, and
$$\gamma_2=\{c_2+re^{i\tau}\,;\,\tau\in (-\pi+\delta, -\delta)\}=\overline{\gamma}_1.$$
The curves $\gamma_1$ and $\gamma_2$ are parametrized by $F$ in the following way: fix $\nu>1$, and consider the two 
 subarcs  $\alpha_1$ and  $\alpha_2$ of $\T$ defined by
\[
\alpha_1=\{e^{i\theta}\in\mathbb T\,;\,\theta\in ((-\pi+\delta)/\nu,-\delta/\nu)\},
\]
and $\alpha_2=\overline{\alpha}_1$. If $\nu$ is sufficiently large, these two arcs are disjoint. We define $F$ on $\T$ in such a way that $F$ is smooth enough and for $j=1,2$, we have
\[
F(e^{i\theta})=c_j+re^{-i\nu\theta}\quad\text{ for every }e^{i\theta}\in\alpha_j.
\]
In particular, $F(\alpha_j)=\gamma_j$.
\par\smallskip
For $j=1,2$, let $\zeta_j$ be the restriction $\zeta_{|\gamma_j}$ of the function $\zeta$ to the curve $\gamma_j$. Then $\zeta_1$ (respectively $\zeta_2$) has an analytic extension to $\Omega_2\setminus (c_1-i\mathbb R_+)$ (respectively on $\Omega_2\setminus (c_2+i\mathbb R_+)$)  given by 
\begin{align*}
\zeta_1(\lambda)&=\exp\left(\frac{1}{\nu}\log_{(-\pi/2,3\pi/2)}\left(\frac{\lambda-c_1}{r}\right)\right)\quad\text{for }\lambda\in\Omega_2\setminus (c_1-i\mathbb R_+)\\
\zeta_2(\lambda)&=\exp\left(\frac{1}{\nu}\log_{(-3\pi/2,\pi/2)}\left(\frac{\lambda-c_2}{r}\right)\right)\quad\text{for }\lambda\in\Omega_2\setminus (c_2+i\mathbb R_+).
\end{align*}
\par\smallskip
The crucial observation in this example is now that, depending on the values of $c_1$, $c_2$ and $r$, the functions $\zeta_1$ and $\zeta_2$ may or not admit an holomorphic extension to $\Omega_2$ - and this changes dramatically the characterization of the embedding property of $T_F$.
\par\smallskip
\textbf{First case: suppose that $c_1$ (and hence $c_2$) belongs to $\Omega_2$}. Then $\zeta_1$ and $\zeta_2$ cannot be extended meromorphically to $\Omega_2$; in fact $F$ satisfies the hypothesis of \Cref{un-th-supplementaire}, and thus $T_F\in\mathcal{B}(H^p)$ is embeddable if and only if $0$ belongs to $\overline{\Omega}_{\infty}$. 
\par\smallskip
\textbf{Second case: suppose that $c_1$ and $c_2$ do not belong to $\overline{\Omega}_2$}. Then both functions $\zeta_1$ and $\zeta_2$ can be extended holomorphically to $\Omega_2$, and their extensions (still denoted by $\zeta_1$ and $\zeta_2$) belong to $A(\Omega_2)$ (the space of continuous functions on $\overline{\Omega}_2$ which are holomorphic  on $\Omega_2$).
Moreover (this will be useful later on) the function $\zeta_1-\zeta_2$ does not vanish on $\overline{\Omega}_2$. Indeed, let $\lambda\in\overline{\Omega}_2$. Then we have 
$$\arg_{(-\pi/2,3\pi/2)}\left(\frac{\lambda-c_1}{r}\right)\in (\delta,\pi-\delta),$$ which implies that $$\frac{1}{\nu}\arg_{(-\pi/2,3\pi/2)}\left(\frac{\lambda-c_1}{r}\right)\in (\delta/\nu,(\pi-\delta)/\nu),$$ and similarly we have $$\arg_{(-3\pi/2,\pi/2)}\left(\frac{\lambda-c_2}{r}\right)\in (-\pi+\delta,-\delta)$$ which implies that $$\frac{1}{\nu}\arg_{(-3\pi/2,\pi/2)}\left(\frac{\lambda-c_2}{r}\right)\in ((-\pi+\delta)/\nu,-\delta/\nu).$$ Thus 
\[
\Im\left(\frac{1}{\nu}\log_{(-\pi/2,3\pi/2)}\left(\frac{\lambda-c_1}{r}\right) 
-\frac{1}{\nu}\log_{(-3\pi/2,\pi/2)}\left(\frac{\lambda-c_2}{r}\right)\right)\in (2\delta/\nu,2(\pi-\delta)/\nu),
\]
and it follows that $\frac{\zeta_1(\lambda)}{\zeta_2(\lambda)}\notin\mathbb R_+$ for every $\lambda\in \overline{\Omega}_2$. In particular, $\zeta_1(\lambda)\neq \zeta_2(\lambda)$ for every $\lambda\in\overline{\Omega}_2$, and moreover the function $\zeta_1-\zeta_2$ is bounded from below on $\overline{\Omega}_2$.
\par\smallskip
Our aim is now to show that in this situation,
 we have the following result:
\begin{claim}\label{une-assertion}
 If   $0$ belongs to $\overline{\Omega}_\infty$ or to $\Omega_0$, then $T_F$ is embeddable into a $C_0$-semigroup.
\end{claim}

\begin{proof}
First note that if $0$ belongs to $\overline{\Omega}_\infty$, then according to \Cref{Th:CSforEmb}, the operator $T_F$ is embeddable into a $C_0$-semigroup of operators on $H^p$. We assume now that $0$ belongs to $\Omega_0$. 
According to \Cref{equivalent}, $T_F$ is embeddable if and only if the multiplication operator $M_\lambda$ by the independent variable $\lambda$, acting on $E^q_F$, embeds into a $C_0$-semigroup of bounded operators $(A_t)_{t>0}$ on $E^q_F$. 
Recall that in this situation, we have
$$E^q_F=\left\{(u_0,u_1)\in E^q(\Omega_1\cup\Omega_2)\oplus E^q(\Omega_2)\,;\,u_0^{int}-\zeta u_1^{int}=u_0^{ext} \;  \textrm{ a.e. on } \gamma_1\cup\gamma_2\right\}.$$
\par\smallskip
\textbf{Step 1:} Suppose that $M_\lambda$ embeds into a semigroup of bounded operators $(A_t)_{t>0}$ on $E^q_F$. Since each operator $A_t$ commutes with $M_\lambda$, \Cref{Fact:commutant1} implies that there exist functions $\alpha_t\in E^1(\Omega_1)$ and $a_t, b_t, c_t, d_t\in E^1(\Omega_2)$ such that if we write, for each $u=(u_0,u_1)\in E^q_F$, $v_t=A_tu$, and $v_t=(v_{t,0},v_{t,1})$, then 
\begin{equation}\label{chose1}
    \begin{cases}
v_{t,0}=\alpha_tu_0&\text{on }\Omega_1\\
v_{t,0}=a_tu_0+b_tu_1&\text{on }\Omega_2\\
v_{t,1}=c_tu_0+d_tu_1&\text{on }\Omega_2\\
\end{cases}.
\end{equation}
Moreover, by \Cref{la-forme-de-A_t-sur-Omega_1}, there exists an analytic determination $\log$ of the logarithm on $\Omega_1$
such that $\alpha_t(\lambda)=e^{t\log(\lambda)}$ for every $\lambda\in\Omega_1$ and every $t>0$.
Since $0\notin \overline{\Omega}_1$, this function $\log$ belongs to  $H^\infty(\Omega_1)$.
It has two different extensions to $\Omega_2$, denoted by $\log_1$ and $\log_2$, respectively, such that $\log_j=\log$ on $\gamma_j$, $j=1,2$. Note that, with the choice of $\gamma_j$ represented on \Cref{Fig5}, we have that $\log_2=\log_1+2i\pi$ on $\Omega_2$.
Then the function $\alpha_t$ has also two analytic extensions to $\Omega_2$ given by 
\begin{equation}\label{eq:3EZZEZDCS8822}
\alpha_{t,j}(\lambda)=e^{t\log_j(\lambda)}, \quad\lambda\in\Omega_2,\quad j=1,2. 
\end{equation}
We have $\alpha_{t,2}=e^{2i\pi t}\alpha_{t,1}$ on $\Omega_2$ and $\alpha_{t,j}^{int}=\alpha_t^{ext}$ on $\gamma_j$. Observe that the function 
\[
\lambda\longmapsto\begin{cases}
\alpha_t(\lambda)&\text{for }\lambda\in\Omega_1\\
\alpha_{t,j}(\lambda)&\text{for }\lambda\in\Omega_2\cup\gamma_j
\end{cases}
\]
is bounded on $\Omega_1\cup\gamma_j\cup\Omega_2$. We will still denote it by $\alpha_{t,j}$.
\par\smallskip
Let 
\begin{equation}\label{eq:definition-Bt-3RDS}
B_t(\lambda)=\begin{pmatrix}
a_t(\lambda)&b_t(\lambda)\\c_t(\lambda)&d_t(\lambda)
\end{pmatrix},\quad\lambda\in\Omega_2,\,t>0.
\end{equation}
Then the conditions in \eqref{chose1} can be rewritten in the following way: $v_{t,0}=\alpha_t u_0 \text{ on }\Omega_1$ and 
$$\begin{pmatrix}
v_{t,0}(\lambda)\\v_{t,1}(\lambda)
\end{pmatrix}=B_t(\lambda)\begin{pmatrix}
u_0(\lambda)\\u_1(\lambda)
\end{pmatrix},\quad\lambda\in\Omega_2. $$
Moreover, we know that for each $u\in E^q_F$, the element $v_t=A_t u$ belongs to $E^q_F$. We are now going to show that the boundary condition defining $E^q_F$ uniquely determines the functions $a_t, b_t, c_t$ and $d_t$.
\par\smallskip
According to \Cref{elements-de-l-espace-modele}, we know that $u=(u_0,u_1)=(1,0)$ belongs to $E_F^q$ and \eqref{chose1} gives $v_{t,0}=\alpha_t u_0$ on $\Omega_1$, and on $\Omega_2$, we have $v_{t,0}=a_t u_0$ and $v_{t,1}=c_t u_0$. Since $a_t$ and $b_t$ belong to $E^1(\Omega_2)$, they have boundary values a.e. on $\partial\Omega_2$, which we denote by $a_t^{int}$ and $c_t^{int}$ respectively. Now the boundary relation of $E_F^q$ gives that $v_{t,0}^{int}-\zeta v_{t,1}^{int}=v_{t,0}^{ext}$ a.e. on $\partial\Omega_2$, that is 
\[
a_t^{int}u_0^{int}-\zeta c_t^{int}u_0^{int}=\alpha_t^{ext}u_0^{ext} \quad\text{a.e. on }\partial\Omega_2.
\]
But $u_0^{int}=u_0^{ext}=1$, which gives 
\begin{equation}\label{chose1bis}
    a_t^{int}-\zeta c_t^{int}=\alpha_t^{ext}~\text{a.e. on}~\partial\Omega_2.
\end{equation}
By using the extensions of the function $\zeta$ through $\gamma_1$ and $\gamma_2$ (which have positive length) and the fact that all the functions $a_t,~\zeta_j c_t$ and $\alpha_{t,j}$ belong to $E^q(\Omega_2)$ (remember that $\zeta_j$ is bounded on $\Omega_2$), we deduce, by the uniqueness property for functions in the Smirnov spaces, that the functions $a_t$ and $c_t$ must necessarily satisfy the following system of equations:
\begin{equation}\label{chose2}
\begin{cases}
a_t-\zeta_1c_t=\alpha_{t,1}\\
a_t-\zeta_2c_t=\alpha_{t,2}=e^{2i\pi t}\alpha_{t,1}
\end{cases}\text{ on }\Omega_2.
\end{equation}
Then we obtain that the functions $a_t$ and $c_t$, if they do exist, must be defined as follows on $\Omega_2$:
\begin{equation}\label{chose3}
a_t=\frac{\zeta_2-e^{2i\pi t}\zeta_1}{\zeta_2-\zeta_1}\,\alpha_{t,1}\quad\text{and}\quad c_t=\frac{1-e^{2i\pi t}}{\zeta_2-\zeta_1}\,\alpha_{t,1}\text{ on }\Omega_2.
\end{equation}
Let us now determine the expressions of the functions $b_t$ and $d_t$ on $\Omega_2$.

According to \Cref{elements-de-l-espace-modele}, we can write $Uz$ as $Uz=(u_0,1)$, where $u_0$ belongs to $H^\infty(\Omega_1\cup\Omega_2)$. Since the pair $(u_0,1)$ belongs to $E^q_F$, we have $u_0^{int}-\zeta=u_0^{ext}$ a.e. on $\partial\Omega_2$. Moreover \eqref{chose1} gives that $v_{t,0}=\alpha_t u_0$ on $\Omega_1$, and on $\Omega_2$ we have $v_{t,0}=a_t u_0+b_t$ and $v_{t,1}=c_t u_0+d_t$. Since $b_t$ and $d_t$ must belong to $E^1(\Omega_2)$, they have boundary  limits a.e. on $\Omega_2$, which we denote by $b_t^{int}$ and $d_t^{int}$ respectively. Now the boundary relation of $E_F^q$ gives that $v_{t,0}^{int}-\zeta v_{t,1}^{int}=v_{t,0}^{ext}$ a.e. on $\partial\Omega_2$, that is 
\[a_t^{int}u_0^{int}+b_t^{int}-\zeta(c_t^{int}u_0^{int}+d_t^{int})=\alpha_t^{ext}u_0^{ext}=\alpha_t^{ext}(u_0^{int}-\zeta)~\text{a.e. on}~\partial\Omega_2.\]
So this means that 
\[(a_t^{int}-\zeta c_t^{int}-\alpha_t^{ext})u_0^{int}+b_t^{int}-\zeta(d_t^{int}-\alpha_t^{ext})=0~\text{a.e. on}~\partial\Omega_2.\]
But we have by (\ref{chose1bis}) that $a_t^{int}-\zeta c_t^{int}=\alpha_t^{ext}$, and hence 
\begin{equation}\label{chose4}
    b_t^{int}-\zeta(d_t^{int}-\alpha_t^{ext})=0~\text{a.e. on}~\partial\Omega_2.
\end{equation}
Then again, by using the extensions $\zeta_1$ and $\zeta_2$ of $\zeta$ to $\Omega_2$ through $\gamma_1$ and $\gamma_2$ respectively, we obtain 
\begin{equation}\label{chose5}
\begin{cases}
b_t-\zeta_1(d_t-\alpha_{t,1})=0\\
b_t-\zeta_2(d_t-e^{2i\pi t}\alpha_{t,1})=0
\end{cases}\text{ on }\Omega_2,
\end{equation}
and thus the functions $b_t$ and $d_t$ are necessarily defined on $\Omega_2$ by the following expressions:
\begin{equation}\label{chose6}
    b_t=\frac{\zeta_1\zeta_2}{\zeta_2-\zeta_1}(e^{2i\pi t}-1)\alpha_{t,1}~\text{ and }~d_t=\frac{e^{2i\pi t}\zeta_2-\zeta_1}{\zeta_2-\zeta_1}\alpha_{t,1}\quad\textrm{ on }\Omega_2.
\end{equation}

Summarizing, we have shown that if $(A_t)_{t>0}$ is a semigroup of bounded operators on $E^q_F$ such that $A_1=M_\lambda$ (no need to suppose here that it is a $C_0$-semigroup), then there exist:

\begin{itemize}
    \item[(A)] an analytic branch $\log$ of the logarithm on $\Omega_1$ with its two different extensions $\log_1$ and $\log_2$ to $\Omega_2$  such that $\log=\log_j$ on $\gamma_j$, $j=1,2$, and if 
    $\alpha_t(\lambda)=e^{t\log(\lambda)}$ for every $\lambda\in\Omega_1$ and every $t>0$, then it has two analytic extensions $\alpha_{t,1}$ and $\alpha_{t,2}$ to $\Omega_2$ given by \eqref{eq:3EZZEZDCS8822};
    \item[(B)] functions $a_t, b_t, c_t$ and $d_t$ on $\Omega_2$ defined by the formulas (\ref{chose3}) and (\ref{chose6})
\end{itemize}
such that for every $t>0$, the action of $A_t$ on a vector $u=(u_0,u_1)\in E^q_F$ is given by $A_tu=v_t=(v_{t,0},v_{t,1})$ defined by using the equations  (\ref{chose1}).

\par\medskip

\textbf{Step 2:} Conversely, since $0\notin\overline{\Omega_1}$, there exists an analytic branch $\log$ of the logarithm which belongs to $H^\infty(\Omega_1)$ and satisfies property (A) above, and for $t>0$, let $a_t,b_t,c_t,d_t$ be defined on $\Omega_2$ by the formulas \ref{chose3}) and (\ref{chose6}). Let $A_t$ be the operator defined on $E_F^q$ by the formula \eqref{chose1}. We claim that

\begin{enumerate}[(a)]
\item $A_t$ is a bounded operator on $E^q_F$;
    
    \item we have $A_1=M_\lambda$;
    
    \item $\|A_t-I\|\longrightarrow 0$ as $t\to0$;
    
    \item $A_{t+s}=A_t A_s$ for every $t,s>0$. 
\end{enumerate}

Since the functions $\alpha_{t,1}$, $\alpha_{t,2}$, $\zeta_1$ and $\zeta_2$ are bounded on ${\Omega}_2$, and since $\zeta_1-\zeta_2$ is bounded from below on ${\Omega}_2$, the functions $a_t, b_t, c_t, d_t$ are bounded on $\Omega_2$. Thus $v_t=A_tu$ belongs to $E^q(\Omega_1\cup\Omega_2)\oplus E^q(\Omega_2)$ for every $u\in E^q_F$. Since $v_t$ satisfies the boundary relation defining $E_F^q$ as well, by construction of the functions $\alpha_t, a_t, b_t, c_t$ and $d_t$, 
it follows that $A_t$ is a bounded linear operator on $E^q_F$. Thus property (a) is satisfied.
By taking $t=1$, we also remark that $A_1=M_\lambda$ (this is (b)), and we finally note that if $B_t$ is defined by \eqref{eq:definition-Bt-3RDS}, then $B_t\longrightarrow I$ uniformly in $\lambda\in \Omega_2$ when $t\to0$; hence (c) is clear as well. So it remains to prove property (d). To this aim, it is sufficient to prove that for every $t,s>0$, $B_{t+s}(\lambda)=B_t(\lambda)B_s(\lambda)$ for every $\lambda\in\Omega_2$. Write
\[a_t=\frac{\alpha_{t,1}}{\zeta_2-\zeta_1}\widetilde{a}_t,~b_t=\frac{\alpha_{t,1}}{\zeta_2-\zeta_1}\widetilde{b}_t,~c_t=\frac{\alpha_{t,1}}{\zeta_2-\zeta_1}\widetilde{c}_t~\textrm{ and }~d_t=\frac{\alpha_{t,1}}{\zeta_2-\zeta_1}\widetilde{d}_t.\]
Since $\alpha_{t+s,1}=\alpha_{t,1}\alpha_{s,1}$ by construction, we just need to check that the following equations hold:
\[\begin{cases}
 \widetilde{a}_t\widetilde{a}_s+\widetilde{b}_t\widetilde{c}_s= \widetilde{a}_{t+s}(\zeta_2-\zeta_1)\\
 \widetilde{a}_t\widetilde{b}_s+\widetilde{b}_t\widetilde{d}_s= \widetilde{b}_{t+s}(\zeta_2-\zeta_1)\\
 \widetilde{c}_t\widetilde{a}_s+\widetilde{d}_t\widetilde{c}_s= \widetilde{c}_{t+s}(\zeta_2-\zeta_1)\\
 \widetilde{c}_t\widetilde{b}_s+\widetilde{d}_t\widetilde{d}_s= \widetilde{d}_{t+s}(\zeta_2-\zeta_1)
\end{cases}\]
\par\smallskip
Let us verify these four equalities. We first have
\begin{align*}
\widetilde{a}_t\widetilde{a}_s+\widetilde{b}_t\widetilde{c}_s~
    &=~(\zeta_2-e^{2i\pi t}\zeta_1)(\zeta_2-e^{2i\pi s}\zeta_1)+\zeta_1\zeta_2(e^{2i\pi t}-1)(1-e^{2i\pi s})\\
    &=~\zeta_2^2-(e^{2i\pi t}+e^{2i\pi s})\zeta_1\zeta_2\\
    &\qquad+e^{2i\pi (s+t)}\zeta_1^2+\zeta_1\zeta_2(e^{2i\pi t}+e^{2i\pi s}-e^{2i\pi (s+t)}-1)\\
    &=~\zeta_2^2-(1+e^{2i\pi (s+t)})\zeta_1\zeta_2+e^{2i\pi(s+t)}\zeta_1^2\\
    &=~(\zeta_2-\zeta_1)(\zeta_2-e^{2i\pi (s+t)}\zeta_1)~=~(\zeta_2-\zeta_1)\widetilde{a}_{t+s}.
\intertext{Then}
\widetilde{a}_t\widetilde{b}_s+\widetilde{b}_t\widetilde{d}_s~
    &=~(\zeta_2-e^{2i\pi t}\zeta_1)\zeta_1\zeta_2(e^{2i\pi s}-1)+\zeta_1\zeta_2(e^{2i\pi t}-1)(e^{2i\pi s}\zeta_2-\zeta_1)\\
    &=~\zeta_1\zeta_2\left[(e^{2i\pi s} -1+e^{2i\pi s}(e^{2i\pi t}-1))\zeta_2\right.\\
    &\qquad\left.-(e^{2i\pi t}(e^{2i\pi s}-1)+e^{2i\pi t}-1)\zeta_1  \right]\\
    &=~\zeta_1\zeta_2(e^{2i\pi (s+t)}-1)(\zeta_2-\zeta_1)~=~(\zeta_2-\zeta_1)\widetilde{b}_{t+s}.
\intertext{Now}
\widetilde{c}_t\widetilde{a}_s+\widetilde{d}_t\widetilde{c}_s~
    &=~(1-e^{2i\pi t})(\zeta_2-e^{2i\pi s}\zeta_1)+(e^{2i\pi t}\zeta_2-\zeta_1)(1-e^{2i\pi s})\\
    &=~(1-e^{2i\pi t}+e^{2i\pi t}(1-e^{2i\pi s}))\zeta_2\\
    &\qquad-(e^{2i\pi s}(1-e^{2i\pi t})+1-e^{2i\pi s})\zeta_1\\
    &=~(1-e^{2i\pi (s+t)})(\zeta_2-\zeta_1)~=~(\zeta_2-\zeta_1)\widetilde{c}_{t+s}.
\intertext{Lastly}
\widetilde{c}_t\widetilde{b}_s+\widetilde{d}_t\widetilde{d}_s~
&=~(1-e^{2i\pi t})\zeta_1\zeta_2(e^{2i\pi s}-1)+(e^{2i\pi t}\zeta_2-\zeta_1)(e^{2i\pi s}\zeta_2-\zeta_1)\\
&=~(e^{2i\pi s}-1-e^{2i\pi (s+t)}+e^{2i\pi t})\zeta_1\zeta_2\\
    &\qquad+e^{2i\pi (s+t)}\zeta_2^2-(e^{2i\pi s}+e^{2i\pi t})\zeta_1\zeta_2+\zeta_1^2\\
&=~e^{2i\pi (s+t)}\zeta_2^2-(1+e^{2i\pi (s+t)})\zeta_1\zeta_2+\zeta_1^2\\
&=~(e^{2i\pi (s+t)}\zeta_2-\zeta_1)(\zeta_2-\zeta_1)~=~(\zeta_2-\zeta_1)\widetilde{d}_{t+s}.
\end{align*}
Hence $B_{s+t}(\lambda)=B_t(\lambda)B_s(\lambda)$ for every $\lambda\in\Omega_2$, and thus the family $(A_t)_{t>0}$ constructed here satisfies (d) as well. 
\par\medskip
\textbf{Conclusion:} We have thus shown that if we define, for each $t>0$, functions $\alpha_t\in E^1(\Omega_1)$ and $a_t, b_t, c_t, d_t\in E^1(\Omega_2)$ which satisfy properties (A) and (B) above, then $M_\lambda$ embeds into the semigroup $(A_t)_{t>0}$ of bounded operators on $E^q_F$ defined by the equations (\ref{chose1}). According to \Cref{equivalent}, this concludes the proof of \Cref{une-assertion}. 
\end{proof}

\begin{remark}
   Observe that whenever $t\in (0,1)$, the function $c_t$ defined in the proof of \Cref{une-assertion} above does not vanish on $\Omega_2$, and thus $A_t$ is not a multiplication operator on $E^q_F$. Hence the construction above provides an example of an embeddable Toeplitz operator satisfying \ref{H1}, \ref{H2} and \ref{H3} such that $M_\lambda$ does not embed into a $C_0$-semigroup of multiplication operators on $E^q_F$.
\end{remark}

\end{example}
\subsection{A characterization}\label{la-caracterisation}
Note that in the previous example, we did not really use the specific form of $F$, nor the parameterizations of the arcs $\gamma_1$ and $\gamma_2$,  but we did use in a crucial way the analytic extensions of the functions 
$\zeta_{|\gamma_1}$ and $\zeta_{|\gamma_2}$
to $\Omega_2$. In the rest of this section, we consider the more general case of a symbol $F$ satisfying \ref{H1} and such that the curve $F(\mathbb T)$ looks topologically like this:

\begin{figure}[ht]\label{figure8}
\begin{tikzpicture}
 \draw(0,0)node{\includegraphics[scale=.78]{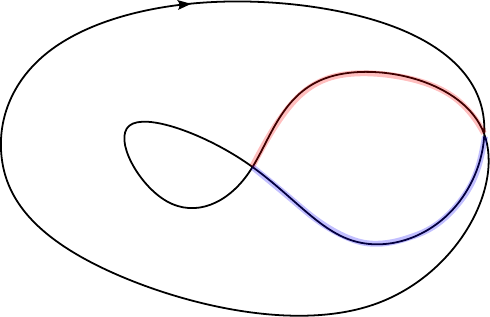}};
 \draw[green!50!black](-1,0)node{$\times$}; \draw[green!50!black](-1,0)node[below]{$0$};
 \draw(2,0)node{$\Omega_2$};
  \draw(0,1.7)node{$\Omega_1$};
  \draw[red](2,1.25)node{$\gamma_1$};
  \draw[blue](2,-1.3)node{$\gamma_2$};
\end{tikzpicture}
    \caption{}
    \label{Fig5}
\end{figure}

\par\smallskip
The following result is essentially a reformulation of what we did in \Cref{l'exemple} in this more general setting:

\begin{proposition}\label{Prop!1}
Let $1<p<\infty$ and let $F$ satisfy \ref{H1}. Suppose that $F(\mathbb T)$ is given by \Cref{Fig5} and that
$0\notin\mathcal O$. Then $T_F$ is embeddable into a $C_0$-semigroup of bounded operators on $H^p$ if and only if 
one of the following two conditions hold:
\begin{itemize}
    \item[(1)] $0$ belongs to the unbounded component of $\mathbb C\setminus \intsp$;
    \item[(2)] $0$ belongs to the bounded component of $\mathbb C\setminus\intsp$ and the following two conditions hold:
    \begin{enumerate}[(i)]
\item $\zeta_{|\gamma_1}$ (resp. $\zeta_{|\gamma_2}$) coincides a.e. on $\gamma_1$ (resp. on $\gamma_2$) with the non-tangential limit of  a meromorphic function $\zeta_1$ (resp. $\zeta_2$) on $\Omega_2$;
\item for every  $(u_0,u_1)\in E^q_F$ and every $t>0$, let $(v_{t,0},v_{t,1}):=A_t(u_0,u_1)$ be defined by the equations \eqref{chose1},
where $\alpha_t(\lambda)=e^{t\log \lambda}$ on $\Omega_1$ for some determination $\log$ of the logarithm on $\Omega_1$  and $a_t,b_t,c_t$ and $d_t$ are given by \eqref{chose3} and \eqref{chose6}. Then the operators $A_t$, $t>0$, defined in this way are bounded on $E^q_F.$
 \end{enumerate}
\end{itemize}
\end{proposition}
Note that the properties from item (ii) are exactly properties (A) and (B) from \Cref{l'exemple!!}.

\begin{proof}First, recall that if $0$ belongs to the unbounded component of $\mathbb C\setminus\intsp$ then $T_F$ is embeddable by \Cref{Th:CSforEmb}. If $0\in\intsp$, then $T_F$ is not embeddable by  \Cref{avant le th C} (recall that $0\notin\mathcal{O}$).
So it is sufficient to show that when $0$ belongs to the bounded component of $\mathbb C\setminus \intsp$, $T_F$ is embeddable if and only if conditions (i) and (ii) of (2) hold.
\par\smallskip
\textbf{Step 1:} Suppose first that $T_F$ is embeddable, i.e. that $M_\lambda$ is embeddable into a $C_0$-semigroup of operators on $E^q_F$.
\par\smallskip
-- Since $0$ belongs to the bounded component of $\mathbb C\setminus \intsp=\mathbb C\setminus\textrm{int}(\overline{\Omega_1\cup\Omega_2)}$, the situation (a) in \Cref{Une-remarque} cannot occur. Hence (i) is satisfied.
\par\smallskip
-- Note that \Cref{Fact:commutant1} and \Cref{la-forme-de-A_t-sur-Omega_1} imply that if $M_\lambda$ is embeddable into a $C_0$-semigroup $(\widetilde A_t)_{t>0}$, there exist $ \alpha_t\in E^1(\Omega_1)$, $ a_t, b_t, c_t$ and $ d_t\in E^1(\Omega_2)$ such that for $(u_0,u_1)\in E^q_F$, and $(v_{t,0},v_{t,1})=\widetilde A_t(u_0,u_1)$ we have
\[v_{t,0}= \alpha_tu_0~\text{on}~\Omega_1,~v_{t,0}= a_tu_0+ b_tu_0~\text{on}~\Omega_2~\text{and}~v_{t,1}= c_tu_0+ d_tu_1~\text{on}~\Omega_2,\]
where $\alpha_t(\lambda)=e^{t\log \lambda}$ for some analytic determination $\log$ of the logarithm on $\Omega_1$.  But using only the fact that $(\widetilde A_t)_{t>0}$ is a $C_0$-semigroup, along with the boundary conditions, we already proved in \Cref{l'exemple!!} that the functions $a_t,b_t,c_t$ and $d_t$ are necessarily given by the formulas \eqref{chose3} and \eqref{chose6},
where $\zeta_1$ and $\zeta_2$ are the meromorphic functions from (i) and the functions $\alpha_{t,1}$ and $\alpha_{t,2}$ are built as in \Cref{l'exemple!!} from the two analytic extensions $\log_1$ and $\log_2$ of the function $\log$ on $\Omega_2$.
In other terms, $T_F$ is embeddable if and only if $M_\lambda$ is embeddable into the $C_0$-semigroup $(A_t)_{t>0}$ defined in (ii). This implies in particular that if $T_F$ is embeddable, these operators $A_t$ are necessarily bounded on $E^q_F$. Hence (ii) is satisfied.

\par\smallskip
\textbf{Step 2:}
Suppose now that (i) and (ii) are satisfied. Then we have seen in \Cref{l'exemple!!} above that 
the family $(A_t)_{t>0}$, defined in (ii), is a semigroup of bounded operators on $E^q_F$ such that $A_1=M_\lambda$. 

To prove that $T_F$ is embeddable,  it is thus sufficient to prove that $A_t\to I$ as $t\rightarrow 0^+$ in the SOT  on $E^q_F$. To this aim, we introduce the operator $C_\zeta$ defined by
\[C_\zeta f(z)=\frac1{2i\pi}\int_{\partial\Omega_2}\frac{f(\lambda)\zeta(\lambda)}{\lambda-z}\,\mathrm d\lambda\quad\textrm{ for every }z\in\Omega_1\cup\Omega_2 \textrm{ and every }f\in E^q(\Omega_2).\]
Since $\partial\Omega_2$ is a Carleson curve, and since $E^q_0(\mathbb C\setminus\overline{\Omega_2})\subseteq E^q(\Omega_1)$, the continuity of the Cauchy transform from $L^q(\Omega_2)$ into $E^q(\Omega_2)$, and into $E^q_0(\mathbb C\setminus\overline{\Omega_2})$ as well, implies the  continuity of the operator $C_\zeta$ from $L^q(\Omega_2)$ into $E^q(\Omega_1\cup\Omega_2)$.

Note also that we have the jump formula (know n as the Sokhotski--Plemelj formula)
\[(\mathcal Cf)^{int}-(\mathcal Cf)^{ext}=f~\text{a.e. on}~\partial\Omega_2~\text{for every}~f\in L^q(\partial\Omega_2),\]
which was obtained by Privalov \cite{MR83565}. In particular, we deduce that
\[(C_\zeta f)^{int}-(C_\zeta f)^{ext}=\zeta f~\text{a.e. on}~\partial\Omega_2\text{ for every }f\in E^q(\Omega_2).\]
In other words, for every $f\in E^q(\Omega_2)$ the pair $\iota f:=(C_\zeta f,f)$ belongs to $ E^q_F$, and the linear map $\iota:E^q(\Omega_2)\to E^q_F$ is bounded. An important fact is that this application $\iota$ allows us to decompose the space $E^q_F$ as a direct sum in the following way:

\begin{fact}\label{fait-decomposition}
Identifying $E^q(\intsp)$ with the closed subspace $E^q(\intsp)\times\{0\}$ of $ E^q_F$, we can decompose $E^q_F$ as the following topological direct sum:
\[E_F^q=E^q(\intsp)\oplus \iota(E^q(\Omega_2)).\]
\end{fact}

\begin{proof}[Proof of \Cref{fait-decomposition}]
Since $C_\zeta0=0$, it is clear that $E^q(\intsp)\cap \iota(E^q(\Omega_2))=\{0\}$. Let $(u_0,u_1)\in E^q_F$ and set $w=u_0- C_\zeta u_1$. Then $w\in E^q(\Omega_1\cup\Omega_2)$ and we have that
\begin{align*}
w^{ext}-w^{int}
&=u_0^{ext}-(C_\zeta u_1)^{ext} -(u_0^{int}-(C_\zeta u_1)^{int}\\
&=(u_0^{ext}-u_0^{int})-((C_\zeta u_1)^{ext}-(C_\zeta u_1)^{int})\\
&=\zeta u_1^{int}-\zeta u_1^{int}=0~\text{a.e. on}~\partial\Omega_2.
\end{align*}
Then by \Cref{StickingLemma}, $w\in E^q(\intsp)$ and $(u_0,u_1)=(w,0)+(C_\zeta u_1,u_1)$. Since the operators $(u_0,u_1)\mapsto (w,0)$ and
$(u_0,u_1)\mapsto (C_\zeta u_1,u_1)$ are bounded from $E_F^q$ into  $E^q(\intsp)$ and $\iota(E^q(\Omega_2))$ respectively, this yields that 
\[E_F^q=E^q(\intsp)\oplus \iota(E^q(\Omega_2)).\qedhere\]
\end{proof}

Let us now finish the proof of \Cref{Prop!1}.
The decomposition given by \Cref{fait-decomposition} implies that $A_t\to I$ as $t\to 0^+$ in the SOT  if and only if the following two properties hold:
\begin{itemize}
    \item[(a)] for every $u_0\in E^q(\intsp)$, $\|A_t(u_0,0)-(u_0,0)\|_{E^q_F}\to0$ as $t\to 0^+$;
    \item[(b)] for every $u_1\in E^q(\Omega_2)$, $\|A_t(C_\zeta u_1,u_1)-(C_\zeta u_1,u_1)\|_{E^q_F}\to0$ as $t\to 0^+$.
\end{itemize}
\par\smallskip
-- Let $u_0\in E^q(\intsp)$ and let $(v_{t,0},v_{t,1})=A_t(u_0,0)$. Let $C>0$ be such that for every $\lambda\in\sigma(T_F)$ and every $0<t<1$, we have $|\lambda|^t\le C$. Then, using \eqref{chose1}, \eqref{chose3} and \eqref{chose6}, we have the following pointwise estimates on $\Omega_1$ and $\Omega_2$ respectively: \[|v_{t,0}|\le C|u_0|~\text{on}~\Omega_1,\]
\[|v_{t,0}|\le C\,\frac{|\zeta_1|+|\zeta_2|}{|\zeta_1-\zeta_2|}\,|u_0|\quad\text{and}\quad |v_{t,1}|\le\frac{2C}{|\zeta_1-\zeta_2|}\,|u_0|~\text{on}~\Omega_2.\]
Since $0\notin\mathcal O$, we have $0\notin\overline{\Omega}_2$, and thus the functions $\alpha_{t,1}$ are bounded and bounded away from $0$ on $\Omega_2$. Dividing $v_{t,0}=a_tu_0$ and $v_{t,1}=c_tu_0$ (on $\Omega_2$) by $\alpha_{t,1}$, we obtain that $\frac{\zeta_1-e^{2i\pi t}\zeta_2}{\zeta_2-\zeta_1}u_0$ and $\frac1{\zeta_2-\zeta_1}u_0$  belong to $E^q(\Omega_2)$ for every $t>0$, and hence $\frac1{\zeta_1-\zeta_2}u_0,~\frac{\zeta_1}{\zeta_1-\zeta_2}u_0$ and $\frac{\zeta_2}{\zeta_1-\zeta_2}u_0$ belong to $E^q(\Omega_2)$. It then follows from \eqref{eq:norm-Smirnov-fc}, \eqref{eq:norm-somme-smirnov}  and from Lebesgue dominated convergence theorem  that
\begin{multline*}
\|A_t(u_0,0)-(u_0,0)\|_{E_F^q}^q=\|(\alpha_t-1)u_0\|_{E^q(\Omega_1)}^q\\+\|v_{t,0}-u_0\|_{E^q(\Omega_2)}^q+\|v_{t,1}\|_{E^q(\Omega_2)}^q\longrightarrow0~\text{as}~t\to0^+,
\end{multline*}
which is the first half of what we wanted to prove.
\par\smallskip
-- Now let $u_1\in E^q(\Omega_2)$ and let $(v_{t,0},v_{t, 1})=A_t(C_\zeta u_1,u_1)$. Then $|v_{t,0}|\le C|C_\zeta u_0|$ on $\Omega_1$. Moreover, observe that on $\Omega_2$, we have 
\begin{align*}
v_{t,0}=a_t C_\zeta u_1+b_t u_1=&\frac{\alpha_{t,1}}{\zeta_2-\zeta_1}\left[(\zeta_2-e^{2i\pi t}\zeta_1)C_\zeta u_1+(e^{2i\pi t}-1)\zeta_1\zeta_2 u_1\right]\\
=&\frac{\alpha_{t,1}}{\zeta_2-\zeta_1}\left[\zeta_2(C_\zeta u_1-\zeta_1 u_1)-e^{2i\pi t}\zeta_1 (C_\zeta u_1-\zeta_2 u_1)\right].
\end{align*}
Again, dividing by $\alpha_{t,1}$, and considering $\tilde v_{t,0}=\frac{v_{t,0}}{\alpha_{1,t}}\in E^q(\Omega_2)$, we deduce that 
\[
\frac{\zeta_2}{\zeta_2-\zeta_1}(C_\zeta u_1-\zeta_1 u_1)=\frac{\tilde v_{1,0}+\tilde v_{1/2,0}}{2}\in E^q(\Omega_2),
\]
and 
\[
\frac{\zeta_1}{\zeta_2-\zeta_1}(C_\zeta u_1-\zeta_2 u_1)=\frac{\tilde v_{1/2,0}-\tilde v_{1,0}}{2}\in E^q(\Omega_2).
\]
Similarly, we have 
\begin{align*}
v_{t,1}=c_t C_\zeta u_1+d_t u_1=&\frac{\alpha_{t,1}}{\zeta_2-\zeta_1}\left[(1-e^{2i\pi t})C_\zeta u_1+(e^{2i\pi t}\zeta_2-\zeta_1)u_1\right]\\
=&\frac{\alpha_{t,1}}{\zeta_2-\zeta_1}\left[C_\zeta u_1-\zeta_1 u_1-e^{2i\pi t}(C_\zeta u_1-\zeta_2 u_1)\right].
\end{align*}
Considering 
$\tilde v_{t,1}=\frac{v_{t,1}}{\alpha_{1,t}}\in E^q(\Omega_2)$, we deduce that 
\[
\frac{1}{\zeta_2-\zeta_1}(C_\zeta u_1-\zeta_1 u_1)=\frac{\tilde v_{1,1}+\tilde v_{1/2,1}}{2}\in E^q(\Omega_2),
\]
and 
\[
\frac{1}{\zeta_2-\zeta_1}(C_\zeta u_1-\zeta_2 u_1)=\frac{\tilde v_{1/2,1}-\tilde v_{1,1}}{2}\in E^q(\Omega_2).
\]
Moreover, for every $0<t<1$, we have the following pointwise estimates on $\Omega_2$:
\begin{align*}
|v_{t,0}|
&=\left|\frac1{\zeta_2-\zeta_1}\Big(\zeta_2(C_\zeta u_1-\zeta_1 u_1)-\zeta_1e^{2i\pi t}(C_\zeta u_1-\zeta_2 u_1)\Big)\alpha_{t,1}\right|\\
&\le \frac C{|\zeta_2-\zeta_1|}\Big(|\zeta_2(C_\zeta u_1-\zeta_1 u_1)|+|\zeta_1(C_\zeta u_1-\zeta_2 u_1)|\Big)
\end{align*}
and 
\begin{align*}
|v_{t,1}|
&=\left|\frac1{\zeta_2-\zeta_1}\Big((C_\zeta u_1-\zeta_1 u_1)-e^{2i\pi t}(C_\zeta u_1-\zeta_2 u_1)\Big)\alpha_{t,1}\right|\\
&\le \frac C{|\zeta_2-\zeta_1|}\Big(|C_\zeta u_1-\zeta_1 u_1|+|C_\zeta u_1-\zeta_2 u_1|\Big)
\end{align*}
By Lebesgue dominated convergence theorem again, we have that
\begin{multline*}
\|A_t(C_\zeta u_1,u_1)\|_{E^q_F}^q=\|(\alpha_t-1)C_\zeta u_1\|^q_{E^q(\Omega_1)}\\+\|v_{t,0}-C_\zeta u_1\|_{E^q(\Omega_2)}^q+\|v_{t,1}-u_1\|_{E^q(\Omega_2)}^q\longrightarrow0~\text{as}~t\to0^+.
\end{multline*}
We conclude that $A_t$ converges to $I$ as $t\to 0^+$ in the Strong Operator Topology, and this terminates the proof of \Cref{Prop!1}.
\end{proof}

\begin{remark}
Let $F$ satisfy \ref{H1} such that $F(\mathbb T)$ is given by \Cref{Fig5} and assume that $T_F$ is embeddable into a $C_0$-semigroup of bounded operators on $H^p$. Assume also that $0$ belongs to the bounded component of $\mathbb C\setminus\intsp$ and $0\notin\mathcal O$. Then,  since $\alpha_{t,1},a_t,b_t,c_t$ and $d_t$ defined in \eqref{eq:3EZZEZDCS8822}, \eqref{chose3} and \eqref{chose6} belong to $E^1(\Omega_2)$, it follows that the four functions
$$\frac1{\zeta_1-\zeta_2}\alpha_{t,1},~\frac{\zeta_1}{\zeta_1-\zeta_2}\alpha_{t,1},~\frac{\zeta_2}{\zeta_1-\zeta_2}\alpha_{t,1} \quad\textrm{ and }\quad\frac{\zeta_1\zeta_2}{\zeta_1-\zeta_2}\alpha_{t,1}
$$ 
also belong to $E^1(\Omega_2)$. Moreover, since $0\notin \mathcal O$, the function $\alpha_{t,1}$ is bounded and bounded away from $0$ on $\Omega_2$. Dividing by $\alpha_{t,1}$, this gives that the functions
\[
\frac1{\zeta_1-\zeta_2},~\frac{\zeta_1}{\zeta_1-\zeta_2},~\frac{\zeta_2}{\zeta_1-\zeta_2}\quad\textrm{ and }\quad\frac{\zeta_1\zeta_2}{\zeta_1-\zeta_2}
\]
belong to $E^1(\Omega_2)$ as well. In particular,  the functions $\zeta_1$ and $\zeta_2$ are quotients of functions in $E^1(\Omega_2)$. Hence they belong to $\mathcal N(\Omega_2)$.
\end{remark}

Our aim is now to reformulate the boundedness condition on the operators $A_t$ in (ii) of \Cref{Prop!1} in a more explicit way, depending only on the functions $\zeta$, $\zeta_1$ and $\zeta_2$. To this purpose, we need to introduce the following Borel measure $\mu$ on $\intsp$ defined as 
\[
\mathrm d\mu(\lambda)=\frac{\mathbf{1}_{\partial\Omega_2}(\lambda)}{|\zeta_1(\lambda)-\zeta_2(\lambda)|^q}|\mathrm d\lambda|,
\]
that is 
\begin{equation}\label{def:measuremucarleson-bornitudeAt}
\mu(A)=\int_{A\cap\partial\Omega_2}\frac{1}{|\zeta_1(\lambda)-\zeta_2(\lambda)|^q}|\mathrm d\lambda|
\end{equation}
for every Borel subset $A$ of $\intsp$.

We recall that $\mu$ is a Carleson measure for $E^q(\intsp)$ if we have the following embedding $E^q(\intsp)\subset L^q(\mu)$, which means that there exists a constant $C>0$ such that, for every $w\in E^q(\intsp)$, we have 
\[
\int_{\partial\Omega_2}\frac{|w(\lambda)|^q}{|\zeta_1(\lambda)-\zeta_2(\lambda)|^q}|\mathrm d\lambda|\le C\|w\|_{E^q(\intsp)}^q.
\]
\begin{proposition}\label{Prop!b}
Fix $0<t<1$, and let $A_t$ be the operator defined in the condition (ii) of \Cref{Prop!1}. Then the following assertions are equivalent.
\begin{enumerate}
    \item The map $A_t:E^q(\intsp)\times\{0\}\longrightarrow E^q_F$ is bounded;
    \item the measure $\mu$ defined by \eqref{def:measuremucarleson-bornitudeAt} is a Carleson measure for $E^q(\intsp)$. 
\end{enumerate}
\end{proposition}
\begin{proof}
Suppose first that $A_t$ is bounded. Consider a function $u_0\in E^q(\intsp)$ and let $(v_{t,0},v_{t,1})=A_t(u_0,0)$. We have $v_{t,1}=c_t u_0$ on $\Omega_2$, and $\|v_{t,1}\|_{E^q(\Omega_2)}\le \|A_t\|\|(u_0,0)\|_{E^q_F}$.
Using the expression of $c_t$ given by (\ref{chose3}), the fact that $\alpha_{t,1}$ is bounded away from $0$ on $\Omega_2$ since
$0\notin\overline{\Omega}_2$, and the fact that $1-e^{2i\pi t}\neq 0$ since $0<t<1$, we obtain that 
\begin{align*}
   \left\| \frac{u_0}{\zeta_2-\zeta_1}\right\|_{E^q(\Omega_2)}&=\left\|\frac{v_{t,1}}{\alpha_t(1-e^{2i\pi t})}\right\|_{E^q(\Omega_2)}\\
    &\le \frac{\|A_t\|}{|1-e^{2i\pi t}|\inf_{\lambda\in \Omega_2}|\lambda|^t}\|(u_0,0)\|_{E^q_F}.
\end{align*}
Moreover, note that $E^q(\intsp)\times \{0\}$ is a closed subspace of $E^q_F$, so that $E^q(\intsp)$ is a Banach space when endowed with the norm $\|u_0\|=\|(u_0,0)\|_{E^q_F}$. Also,  $\|u_0\|_{E^q(\intsp)}\le\|(u_0,0)\|_{E^q_F}=\|u_0\|$ for every $u_0\in E^q(\intsp)$, 
so that by the Banach isomorphism theorem, the two norms $\|\,.\,\|$ and $\|\,.\,\|_{E^q(\intsp)}$ are equivalent
on $E^q(\intsp)$. Hence there exists a positive constant $a$ such that $\|u_0\|\le a\|u_0\|_{E^q(\intsp)}$ for every $u_0\in E^q(\intsp)$. It follows that for every $u_0\in E^q(\intsp)$, we have
\begin{align*}
\left(\int_{\partial\Omega_2}\frac{|u_0(\lambda)|^q}{|\zeta_1(\lambda)-\zeta_2(\lambda)|^q}|\mathrm d\lambda|\right)^{1/q}
&=\left\|\frac{u_0}{\zeta_2-\zeta_1}\right\|_{E^q(\Omega_2)}\\
&\le \frac{\|A_t\|}{|1-e^{2i\pi t}|\inf_{\lambda\in \Omega_2}|\lambda|^t}\|(u_0,0)\|_{E^q_F}\\
&\le \frac{a\,\|A_t\|}{|1-e^{2i\pi t}|\inf_{\lambda\in \Omega_2}|\lambda|^t}\|u_0\|_{E^q(\intsp)}.
\end{align*}
This means exactly that the measure $\mu$ defined by \eqref{def:measuremucarleson-bornitudeAt} is a Carleson measure for $E^q(\intsp)$, and (2) is proved.

\par\smallskip
Conversely, suppose now that the measure $\mu$ defined by \eqref{def:measuremucarleson-bornitudeAt} is a Carleson measure for $E^q(\intsp)$, and let $C>0$ be such that for every $w\in E^q(\intsp)$,
\begin{equation}\label{etoile}
\int_{\partial\Omega_2}\frac{|w(\lambda)|^q}{|\zeta_1(\lambda)-\zeta_2(\lambda)|^q}|\mathrm d\lambda|\le C\|w\|_{E^q(\intsp)}^q.
\end{equation}
Recall that the non-tangential limit of $\zeta_j$ coincides with $\zeta$ a.e. on $\gamma_j$. In particular we have that $|\zeta_j|=1$ a.e. on $\gamma_j$. Now, remark that 
\[\frac{\zeta_2-e^{2i\pi t}\zeta_1}{\zeta_2-\zeta_1}=1+\frac{(1-e^{2i\pi t})\zeta_1}{\zeta_2-\zeta_1}=e^{2i\pi t}+\frac{(1-e^{2i\pi t})\zeta_2}{\zeta_2-\zeta_1}.\]
We will decompose the integral on $\partial\Omega_2$ in (\ref{etoile}) as a sum of two integrals over $\gamma_1$ and $\gamma_2$ respectively,  and use the two forms of $\frac{\zeta_2-e^{2i\pi t}\zeta_1}{\zeta_2-\zeta_1}$ above to estimate these two integrals. Let $u_0\in E^q_F$ and $(v_{t,0},v_{t,1})=A_t(u_0,0)$. Then $v_{t,0}=\alpha_t u_0$ on $\Omega_1$, so that
\[\|v_{t,0}\|_{E^q(\Omega_1)}\le \sup_{\lambda\in \Omega_1}|\lambda|^t\|u_0\|_{E^q(\Omega_1)}\le\sup_{\lambda\in \Omega_1}|\lambda|^t\|(u_0,0)\|_{E^q_F}.\]
Moreover, $v_{t,0}=a_t u_0$ on $\Omega_2$, so that
\begin{align*}
\|v_{t,0}\|_{E^q(\Omega_2)}^q&=\int_{\partial\Omega_2}\left|\frac{\zeta_2-e^{2i\pi t}\zeta_1}{\zeta_2-\zeta_1}\alpha_tu_0\right|^q|\mathrm d\lambda|\\
&=\int_{\gamma_1}\left|1+\frac{(1-e^{2i\pi t})\zeta_1}{\zeta_2-\zeta_1}\right|^q|\alpha_tu_0|^q|\mathrm d\lambda|\\
&+\int_{\gamma_2}\left|e^{2i\pi t}+\frac{(1-e^{2i\pi t})\zeta_2}{\zeta_2-\zeta_1}\right|^q|\alpha_tu_0|^q|\mathrm d\lambda|\\
&\le C_q\, \sup_{\lambda\in\Omega_2}|\lambda|^{qt}\left(\|u_0\|_{E^q(\Omega_2)}^q+2^q\int_{\partial\Omega_2}\left|\frac{u_0}{|\zeta_1-\zeta_2}\right|^q|\mathrm d\lambda|\right)\\
&\le C_q\, \sup_{\lambda\in\Omega_2}|\lambda|^{qt}(\|u_0\|_{E^q(\Omega_2)}^q+2^qC\|u_0\|_{E^q(\intsp)}^q)\\
&\le C_q\,\sup_{\lambda\in\Omega_2}|\lambda|^{qt}(1+2^qC)\|(u_0,0)\|_{E^q_F}^q,
\end{align*}
where $C_q$ is a positive constant such that $(1+x)^q\le C_q(1+ x^q)$ for every $x>0$.
Also, $v_{t,1}=c_t u_0$ on $\Omega_2$. Here the estimate is more direct:
\begin{align*}
\|v_{t,1}\|_{E^q(\Omega_2)}
&\le 2\sup_{\lambda\in\Omega_2}|\lambda|^t\left\|\frac{u_0}{\zeta_1-\zeta_2}\right\|_{E^q(\Omega_2)}\\
&\le 2C^{1/q}\sup_{\lambda\in\Omega_2}|\lambda|^t\|u_0\|_{E^q(\intsp)}\\
&\le 2C^{1/q}\sup_{\lambda\in\Omega_2}|\lambda|^t\|(u_0,0)\|_{E^q_F}.
\end{align*}
Hence we deduce that 
\[\|A_t(u_0,0)\|_{E^q_F}^q\le \sup_{\lambda\in\sigma(T_F)}|\lambda|^{qt}(1+2^qC)(1+C_q)\|(u_0,0)\|_{E^q_F}^q\]
for every $u_0\in E^q(\intsp)$,
which means that the restriction of the  operator $A_t$ to $E^q(\intsp)\times\{0\}$ is bounded.
\end{proof}

\begin{remark}
\Cref{Prop!b} shows the following: if $A_{t_0}$ is a bounded operator from $E^q(\intsp)\times\{0\}$ into $E^q_F$ for \emph{some} $t_0\in (0,1)$, then it is bounded for \emph{all} $t\in(0,1)$, and hence for all $t>0$ by the semigroup property. The same remark holds for \Cref{Prop!a} below.
\end{remark}

\begin{proposition}\label{Prop!a}
Fix $0<t<1$ and let $A_t$ be defined as in (ii) of \Cref{Prop!1}. Then the following assertions are equivalent:
\begin{enumerate}
    \item The map $A_t:\iota(E^q(\Omega_2))\longrightarrow E^q_F$ is bounded;
    \item The four maps $$Z_1:w\mapsto \frac1{\zeta_1-\zeta_2}\,(C_\zeta w-\zeta_1 w)\quad Z_2:w\mapsto \frac1{\zeta_1-\zeta_2}\,(C_\zeta w-\zeta_2 w)$$ $$Z_3:w\mapsto \frac{\zeta_2}{\zeta_1-\zeta_2}\,(C_\zeta w-\zeta_1 w)\quad  Z_4:w\mapsto \frac{\zeta_1}{\zeta_1-\zeta_2}\,(C_\zeta w-\zeta_2 w)$$ define bounded operators from $E^q(\Omega_2)$ into itself.
\end{enumerate}
\end{proposition}

\begin{proof}
Let us first remark that the following equalities hold:
\begin{equation}\label{Eq37}
\frac{\zeta_2-e^{2i\pi t}\zeta_1}{\zeta_2-\zeta_1}=1+\frac{(1-e^{2i\pi t})\zeta_1}{\zeta_2-\zeta_1}=e^{2i\pi t}+\frac{(1-e^{2i\pi t})\zeta_2}{\zeta_2-\zeta_1},
\end{equation}
and
\begin{equation}\label{Eq38}
  \frac{e^{2i\pi t}\zeta_2-\zeta_1}{\zeta_2-\zeta_1}=e^{2i\pi t}-\frac{(1-e^{2i\pi t})\zeta_1}{\zeta_2-\zeta_1}=1-\frac{(1-e^{2i\pi t})\zeta_2}{\zeta_2-\zeta_1}. 
\end{equation}
Let now $u_1\in E^q(\Omega_2)$ and $(v_{t,0},v_{t,1})=A_t(C_\zeta u_1,u_1)$. Then
$v_{t,0}=a_t C_\zeta u_1+b_t u_1$ and $v_{t,1}=c_t C_\zeta u_1+d_t u_1$ on $\Omega_2$ by (\ref{chose1}), and thus
using \eqref{Eq37}, we have that
\begin{equation}\label{eq:sdfsds3E43ZDSD0}
v_{t,0}=\alpha_{t,1}\Big(C_\zeta u_1-(1-e^{2i\pi t})Z_4u_1\Big)=\alpha_{t,1}\Big(e^{2i\pi t}C_\zeta u_1-(1-e^{2i\pi t})Z_3u_1\Big)
\end{equation}
on $\Omega_2$,.
Using \eqref{Eq38}, we have also
\begin{equation}\label{sdssds2E3ZESD323D}
v_{t,1}=\alpha_{t,1}\Big( u_1-(1-e^{2i\pi t})Z_2u_1\Big)=\alpha_{t,1}\Big(e^{2i\pi t} u_1-(1-e^{2i\pi t})Z_1u_1\Big)
\end{equation}
on $\Omega_2$.
Putting together these equalities, we deduce that, on $\Omega_2$, we have
\begin{align}\label{etoile2}
A_t(C_\zeta u_1,u_1)
&=(\alpha_{t,1}C_\zeta u_1,\alpha_{t,1}u_1)+(1-e^{2i\pi t})(\alpha_{t,1}Z_4u_1,\alpha_{t,1}Z_2u_1).
\end{align}
Since $0\notin \overline{\Omega}_2$, we obtain the equivalence of \Cref{Prop!a}. Indeed, suppose that $A_t$ is bounded from $\iota(E^q(\Omega_2))$ to $E^q_F$ and let $u_1\in E^q(\Omega_2)$ and $(v_{t,0},v_{t,1})=A_t(C_\zeta u_1,u_1)$. Then 
\begin{align*}
\|Z_1u_1\|_{E^q(\Omega_2)}
&=\frac{1}{|1-e^{2i\pi t}|}\left\|e^{2i\pi t}u_1-\frac{v_{t,1}}{\alpha_{t,1}}\right\|_{E^q(\Omega_2)}\\
&\le \frac{1}{|1-e^{2i\pi t}|}\left(\|u_1\|_{E^q(\Omega_2)}+\frac{1}{\inf_{\lambda\in\Omega_2}|\lambda|^t}\|A_t(C_\zeta u_1,u_1)\|_{E^q_F}\right)\\
&\le\frac{1}{|1-e^{2i\pi t}|}\left(1+\frac{\|A_t\|\|\iota\|}{\inf_{\lambda\in\Omega_2}|\lambda|^t}\right)\|u_1\|_{E^q(\Omega_2)}.
\end{align*}
So $Z_1$ is indeed a bounded operator from $E^q(\Omega_2)$ into itself. The operators $Z_2,Z_3$ and $Z_4$ are shown to be bounded in exactly the same way, using the equalities in \eqref{eq:sdfsds3E43ZDSD0} and \eqref{sdssds2E3ZESD323D}.
\par\smallskip
Suppose now that the four operators $Z_1,Z_2,Z_3$ and $Z_4$ are bounded on $E^q(\Omega_2)$. Let $u_1\in E^q(\Omega_2)$. Then, using \eqref{etoile2}, we have
\begin{align*}
\|A_t(C_\zeta u_1,u_1)\|_{E^q_F}^q
&=\|\alpha_tC_\zeta u_1\|_{E^q(\Omega_1)}^q\\
&\qquad\qquad+\left\|\alpha_{t,1}\Big(C_\zeta u_1-(1-e^{2i\pi t})Z_4u_1\Big)\right\|^q_{E^q(\Omega_2)}\\
&\qquad\qquad+\left\|\alpha_{t,1}\Big( u_1-(1-e^{2i\pi t})Z_2 u_1\Big)\right\|_{E^q(\Omega_2)}^q\\
&\le \sup_{\lambda\in\sigma(T_F)}|\lambda|^{tq}\Big(1+(1+2\|Z_4\|)^q\\
&\qquad\qquad+(1+2\|Z_2\|)^q\Big)\|(C_\zeta u_1,u_1)\|_{E^q_F}^q.
\end{align*}
So we deduce that $A_t$ is bounded from $\iota(E^q(\Omega_2)$ into $E^q_F$, and this terminates the proof of \Cref{Prop!a}
\end{proof}

\begin{remark}
It follows from \eqref{eq:sdfsds3E43ZDSD0} and \eqref{sdssds2E3ZESD323D} that the operator $Z_1$ is bounded if and only if $Z_2$ is, and that $Z_3$ is bounded if and only if $Z_4$ is.
\end{remark}

Combining \Cref{Prop!1,Prop!b,Prop!a}, this gives the following result, which is our final characterization of the embeddability of $T_F$ when the curve $F(\T)$ looks as in \Cref{Fig5}:

\LaCaracterisation*

\section{Sectorial Toeplitz operators}\label{Section: sectorial}
The most natural way to prove that an operator $T$ is embeddable into a $C_0$-semigroup is to construct a semigroup $(T_t)_{t>0}$ with $T_1=T$ via a functional calculus, and to prove that this semigroup $(T_t)_{t>0}$ converges to the identity operator in the Strong Operator Topology when $t\to0^+$. For example, the analytic functional calculus on a neighborhood of the spectrum of an operator $T$ gives the embeddability as soon as $0$ belongs to the unbounded component of $\mathbb C\setminus \sigma(T)$. In the first part of this paper, using the  functional calculus from \cite{Yakubovich1991}, we extended this criterion for embeddability and showed that a Toeplitz operator $T_F$ with a smooth symbol $F$ is embeddable into a $C_0$-semigroup on $H^p$ as soon as $0$ belongs to the unbounded component of $\mathbb C\setminus\intsp$ (see \Cref{Th:CSforEmb}). Even without these additional smoothness conditions on the symbol, one may attempt to use alternative tools, such as the numerical range, to define a functional calculus.
\par\smallskip
Recall that the numerical range $W(T)$ of a bounded operator $T$ on a Hilbert space $H$ is the convex set defined by 
\[W(T)=\{\ps{Tx,x}\,;\,x\in H~\text{and}~\|x\|=1\},\]
and that its closure contains the spectrum of $T$. See \cite{GustafsonRao1997} for an account of the properties of the numerical range. In 1999, B. and F. Delyon proved that for every bounded and convex domain $\Omega$ of $\mathbb{C}$, for every operator $T$ such that $W(T)\subseteq \Omega$ and for every polynomial $P$, we have 
\[
\|P(T)\|\le C_\Omega\sup_\Omega|P|\quad \text{where}~C_\Omega=\left(\frac{2\pi \cdot diam(\Omega)^2}{Area(\Omega)}\right)^3+3.\]
In other terms, the set $\Omega$ is $C_\Omega$-spectral for $T$ \cite{DelyonDelyon1999}. Some years later, Crouzeix proved in \cite{Crouzeix2007} that $C_\Omega$ can be replaced by a universal constant (which is $11.08$) and conjectured in \cite{Crouzeix2004} that the numerical range is always $2$-spectral. 
In 2017, Crouzeix and Palencia improved the universal constant $11.08$ and obtained that the numerical range is always $(1+\sqrt2)$-spectral \cite{CrouzeixPalencia2017} (see also \cite{RansfordSchwenninger2018} for a simpler proof of the Crouzeix-Palencia result and \cite{ClouatreOstermannRansford2023} for an abstract version of it). Note that, very recently, Malman, Mashreghi, O'Laughlin and Ransford obtained that for an operator $T$, the best constant $K_T$ such that $W(T)$ is $K_T$-spectral for $T$ must satisfy $K_T<1+\sqrt2$ \cite{MalmanMashreghiOLoughlinRansford2025}.
\par\smallskip
A natural interest for the numerical range in the context of embeddability is that, thanks to the results of B. and F. Delyon and those related to the Crouzeix conjecture, the analytic functional calculus on a neighborhood of the spectrum of $T$ can be extended continuously, for the sup-norm on $W(T)$, to functions that are analytic on $\mathrm{int}(W(T))$ and continuous on $\overline{W(T)}$. In particular, if $0$ belongs to $\mathbb C\setminus \mathrm{int}(W(T)),$ this functional calculus allow us to construct a semigroup $(T_t)_{t>0}$ such that $T_1=T$. Unfortunately, when $0$ belongs to $\partial W(T)$, there is no guarantee that  $(T_t)_{t\ge 0}$ is a $C_0$-semigroup, i.e. that $(T_t)_{t>0}$ converge to $I$ in the Strong Operator Topology when $t\to0^+$.
\par\smallskip
To bypass this problem, we will  consider instead, in \Cref{Subsec:defsectorial}, the functional calculus for sectorial operators. In \Cref{numerical-range}, thanks to the link with the numerical range, we will obtain the SOT convergence of $(T_t)_{t>0}$ to $I$ as soon as $0\in\mathbb C\setminus W(T)$. Thanks to Coburn's lemma, we also improve this condition in the case where $T$ is a Toeplitz operator, and prove that $T_F$ is embeddable as soon as $0\in \mathbb C\setminus \mathrm{int}(W(T_F))$ (this is \Cref{Theo4}). We will finish this section by studying the link between sectorial operators and other tools such as the Kreiss constant of sectors.

\subsection{Definition and embedding of sectorial operators}\label{Subsec:defsectorial}
Given $\omega\in[0,\pi]$ denote by $S_\omega$ the subset of the complex plane defined by
\[S_\omega=\begin{cases}
\{z\in\mathbb C\,;\,z\neq0~\text{and}~|\arg(z)|<\omega\} &\text{if}~\omega\in(0,\pi]\\
(0,+\infty)&\text{if}~\omega=0.
\end{cases}\]
When $\omega\in(0,\pi]$, the set $S_\omega$ is an open sector of the complex plane whose opening is $2\omega$ and with vertex at the origin. 
\par\smallskip
Let $X$ be a Banach space, and let $T\in B(X)$ be a bounded linear operator on $X$. The operator $T$ is said to be \emph{sectorial of angle} $\omega$ for some $\omega\in[0,\pi)$ if the following two properties hold:
\begin{enumerate}
    \item $\sigma(T)\subseteq\overline{S_\omega}$, where $\overline{S_\omega}$ is the closure of the sector $S_\omega$;
    \item for every $\phi\in(\omega,\pi)$, we have $$\sup\{\|\lambda(\lambda-T)^{-1}\|\,;\,\lambda\in \mathbb C\setminus \overline{S_{\phi}}\}<\infty.$$
\end{enumerate}
We refer the reader to \cite{Haase2006} for a comprehensive presentation of sectorial operators (which may in general be unbounded) and their numerous applications. One of the main interests of sectorial operators is the fact that they admit useful functional calculi.
We remind here a few facts concerning the so-called natural functional calculus for bounded sectorial operators (see \cite{Haase2006}*{Subsection 2.5.2} for details).
Let $T\in B(X)$ be a bounded  sectorial operator on $X$ of angle $\omega$, and let $\phi\in(\omega,\pi)$. 
Consider the following class of holomorphic functions on $S_\phi$, satisfying a decay condition at $0$:
\[\mathcal E_0(S_\phi)=\{f\in Hol(S_\phi)\,;\,f(z)=O(|z|^\alpha)~\text{when}~z\longrightarrow0~\text{for some}~\alpha>0\}.\]
Note that in the definition above, the power $\alpha$ which appears depends on $f$.

It turns out that $T$ admits an $\mathcal E_0(S_\phi)$-functional calculus, defined as a Cauchy integral on a certain contour bounding a sector $S_{\omega'}$, $\omega'\in(\omega, \phi)$ except for the region near $\infty$, where it avoids $\infty$ and stays away from $\sigma(T)$. More precisely, let $\omega'\in (\omega,\phi)$, $R>\|T\|$, and let $\Gamma$ be the positively oriented contour $\Gamma=\partial(S_{\omega'}\cap R\mathbb D)$. Then the $\mathcal E_0(S_\phi)$-functional calculus is given by
\[f(T):=\int_{\Gamma}f(z)(z-T)^{-1}\,\frac{\mathrm dz}{2i\pi}\quad\text{for every}~ f\in \mathcal E_0(S_\phi).\]
It is not difficult to check that if $E$ is a closed subspace of $X$ which is hyperinvariant with respect to $T$ (i.e. invariant with respect to every operator in the commutant of $T$), then $E$ is also invariant with respect to $f(T)$ for every $f\in\mathcal E_0(S_\phi)$.
\par\smallskip
Since the function $z\mapsto z^t$ belongs to $\mathcal E_0(S_\phi)$ for every $t>0$, this functional calculus allows us to construct a semigroup $(T_t)_{t>0}$ given by $T_t=z^t(T)$. Then Proposition 3.1.15 in \cite{Haase2006} yields the following result:

\begin{proposition}\label{Prop1}
Let $T\in B(X)$ be a sectorial operator of angle $\omega$, and let $(T_t)_{t>0}$ be the semigroup constructed thanks to the $\mathcal E_0(S_\phi)$-functional calculus for $T$, for some $\phi\in(\omega,\pi)$. Let $x\in X.$ Then 
\[x\in\overline{\Ran(T)}\textrm{ if and only if } \|T_tx-x\|\longrightarrow0~\text{when}~t\to0.\]
In particular $(T_t)_{t>0}$ is a $C_0$-semigroup if and only if $T$ has a dense range.
\end{proposition}
Note that if the Banach space is reflexive and if $T\in B(X)$ is a sectorial operator, then the space $X$ can be decomposed as
\begin{equation}\label{Eq:decompEspRefl}
    X=\ker(T)\oplus\overline{\Ran(T)}.
\end{equation}
See \cite{Haase2006}*{Prop. 2.1.1} for details. Thus when $X$ is a separable Hilbert space, \Cref{Prop1} yields a characterization  of sectorial operators which can be embedded in a $C_0$-semigroup. This result is essentially a consequence of \cite{Eisner-E1}*{Th. 1.4} and \cite{Eisner-E1}*{Prop. 1.13}, but we provide a proof for completeness's sake.

\begin{proposition}\label{plongement-sectoriel}
Let $H$ be a separable Hilbert space, and let $T\in B(H)$ be a sectorial operator. Then $T$ is embeddable into a $C_0$-semigroup if and only if $\dim(\ker(T))=0$ or $\dim(\ker(T))=\infty$.
\end{proposition}

\begin{proof}
If $\dim\ker(T)=0$, then it follows from  \Cref{Prop1} that $T$ is embeddable into a $C_0$-semigroup. If $\ker(T)$ is finite-dimensional and non-zero, then $T$ is not embeddable  by \Cref{Th:CN-plong-ker}. So it remains to consider the case where $\ker(T)$ is infinite-dimensional. 
\par\smallskip
Since $\ker(T)$ is an infinite dimensional separable Hilbert space, we can apply \cite{Eisner2010}*{Lemma V.1.12} to deduce that the zero operator  on $\ker(T)$ is embeddable into a $C_0$-semigroup $(A_t)_{t>0}$ of bounded operators on $\ker(T)$. 
Now let $(S_t)_{t>0}$ be the semigroup of operators on $H$ constructed thanks to the functional calculus for $T$ (i.e. $S_t=z^t(T)$ for each $t>0$). Since $\overline{\Ran(T)}$ is hyperinvariant by $T$,  for every $t>0$, the operator $B_t$, defined as the restriction of $S_t$ to $\overline{\Ran(T)}$, belongs to $\mathcal B(\overline{\Ran (T)})$. 
Whence $(B_t)_{t>0}$ is a $C_0$-semigroup on $\overline{\Ran(T)}$ by \Cref{Prop1}. Using the decomposition $H=\ker(T)\oplus\overline{\Ran(T)}$, define for each $t>0$ an operator $T_t=A_t\oplus B_t$ on $H$, i.e. set
\[T_t(x+y)=A_tx+B_ty\quad\text{for every}~x\in \ker(T)~\text{and every}~y\in\overline{\Ran(T)}.\]
Then it is clear that $(T_t)_{t>0}$ is a semigroup. Let $x\in H$, and let $y\in\ker(T), z\in\overline{\Ran(T)}$ be such that $x=y+z$. Then
\[\|T_tx-x\|=\|A_ty-y+B_tz-z\|\le\|A_ty-y\|+\|B_tz-z\|\longrightarrow0~\text{when}~t\to 0. \]
Hence $(T_t)_{t>0}$ is a $C_0$-semigroup, and  $T_1=A_1\oplus B_1=0\oplus {S_1}_{|\overline{\Ran(T)}}=T.$ We have thus proved that $T$ is embeddable.
\end{proof}

\begin{remark}\label{zero-plonge}
The proof of \Cref{plongement-sectoriel} actually shows that if $T$ is a sectorial operator on a Banach space $X$, $T$ is embeddable into a $C_0$-semigroup if and only if the zero operator on $\ker(T)$  is embeddable. Let us point out that there exist Banach spaces $X$ on which the zero operator is not embeddable: it was shown by Lotz in \cite{Lotz} that whenever $X$ is a Grothendieck space with the Dunford-Pettis property, every $C_0$-semigroup of operators on $X$ is uniformly continuous. Obviously, the zero operator cannot be embedded into a uniformly continuous semigroup. Examples of Grothendieck spaces with the Dunford-Pettis property are the spaces $\ell_{\infty}$ and $L^{\infty}(\Omega,\Sigma,\mu)$, as well as $C(K)$-spaces when $K$ is a compact $\sigma$-Stonian space. These spaces are necessarily non-separable. If $X$ is a Banach lattice with a quasi-interior point, $X$ has the property that every $C_0$-semigroup of operators on $X$ is uniformly continuous if and only if $X$ is a Grothendieck space with the Dunford-Pettis property \cite{VanNeerven}.
\end{remark}

In the setting of Toeplitz operators, recall that Coburn's lemma asserts that for every symbol $F\in L^{\infty}(\mathbb T)$, either $T_F$ or $T_F^*$ is injective  and so $T_F$ is either injective or has dense range. In particular, if $T_F$ is sectorial then $T_F$ is injective with dense range on $H^p$ by \eqref{Eq:decompEspRefl}. So we obtain the following direct consequence of \Cref{Prop1}:

\begin{theorem}\label{Theo2}
Let $p\in(1,+\infty)$ and let $F\in L^\infty(\mathbb T)$. If there exists a constant $a\in\mathbb C\setminus\{0\}$ such that $aT_F$ is a sectorial operator on $H^p$, then $T_F$ is embeddable into a $C_0$-semigroup of operators on $H^p$.
\end{theorem}

\subsection{Link with the numerical range}\label{numerical-range}
Let us begin this section with a word of caution: in this section, we depart from the (unusual) choice of the scalar product on a complex separable Hilbert space made in the rest of the paper - which was linear in both variables. Here the scalar product will be as usual linear in the first variable and antilinear in the second variable.
\par\smallskip
Let $H$ be a complex separable Hilbert space, and let $T\in B(H)$.  The numerical range of $T$ is defined as
\[W(T)=\{\ps{Tx,x}\,;\,x\in H~\text{and}~\|x\|=1\}.\]
Then $W(T)$ is a bounded convex subset of $\mathbb C$ which satisfies $\overline{\conv}(\sigma(T))\subset\overline{W(T)}$, with equality for normal operators. See \cite{GuftasonRao1997} for a detailed account on the properties of the numerical range. It follows from von Neumann equality that  a closed half-plane $A$ of $\mathbb C$ is spectral for $T$ (i.e. for every rational function $f$ which is bounded on $A$, we have
$\|f(T)\|\le \sup_{z\in A}|f(z)|$)
if and only if it contains $W(T)$. 
\par\smallskip
Indeed, let $A$ be a closed half-plane. Without loss of generality,  we can assume that $A=\{z\in\mathbb C:~\Re(z)\ge0\}$ and $\sigma(T)\subset A$. Let $S=\phi(T)$ with $\phi(z)=\frac{1-z}{1+z}$, $z\in\mathbb{C}\setminus\{-1\}$. Since $\phi$ is a conformal map from $A$ onto $\mathbb D$, it follows that $A$ is spectral for  $T$ if and only if the open unit disk $\mathbb D$ is spectral for $S$, which is equivalent to the condition $\|S\|\le1$ by von Neumann inequality. But note now that $S=(I-T)(I+T)^{-1}$ is a contraction if and only if $W(T)\subset A$. Indeed, the operator $I+T$ is invertible, and for every $x\in H$ and $y:=(I+T)x$, we have $Sy=(I-T)x$. Hence it follows that $\|S\|\le 1$ if and only if 
\begin{equation}\label{Eq:vN-W(T)1}
    \|(I-T)x\|\le\|(I+T)x\|~\text{ for every}~x\in H.
\end{equation}
Since $\|(I\pm T)x\|^2=\|x\|^2+\|Tx\|^2\pm 2\Re\ps{Tx,x}$, \eqref{Eq:vN-W(T)1} is equivalent to the condition $\Re\ps{Tx,x}\ge0$ for all $k\in H$, i.e. $W(T)\subset A$.
\par\smallskip
This observation combined with \Cref{Prop1} yields the following sufficient condition for embeddability:

\begin{theorem}\label{Theo3}
Let $T\in B(H)$, where $H$ is a complex separable Hilbert space. If $0$ does not belong to $ W(T)$, then $T$ is embeddable into a $C_0$-semigroup of operators on $H$.
\end{theorem}

\begin{proof}
If $0\notin\overline{W(T)}$, then there exists a determination of the logarithm which is analytic on a neighborhood of $\overline{W(T)}$, and hence on a neighborhood of $\sigma(T)$. The embeddability of $T$ is clear in this case.
\par\smallskip
Assume now that $0\in\partial W(T)$. Since $W(T)$ is convex, there exists a  closed half plane $A$ such that $0\in\partial A$ and $W(T)\subseteq A$. Multiplying if necessary $T$ by a unimodular constant, we can assume that $A=\{z:\Re(z)\ge0\}$. The fact that $A=\overline{S_{\pi/2}}$ contains the numerical range implies that $T$ is sectorial of angle $\pi/2$.
Indeed, we have \[\sigma(T)\subset\overline{W(T)}\subset \overline{S_{\pi/2}.}\]
Let now $\lambda\in \{z:\Re(z)<0\}$. We have
\[\sup_{z\in A}|(\lambda-z)^{-1}|=\frac{1}{\dist(\lambda,A)}=\frac{1}{|\Re(\lambda)|}\cdot\]
It follows from the observation above that
\[\|(\lambda-T)^{-1}\|\le\frac1{\dist(\lambda,A)}=\frac1{|\Re(\lambda)|}\cdot\]
So let $\phi\in(\pi/2,\pi]$ and $\lambda\in \mathbb C\setminus\overline{S_\phi}$. Then $|\arg(\lambda)|\in(\phi,\pi]$ and thus
\[\|\lambda(\lambda-T)^{-1}\|\le\frac{|\lambda|}{|\Re(\lambda)|}\le\frac1{|\cos\phi|}<\infty. \]
Hence $T$ is sectorial of angle $\pi/2$.
\par\smallskip
Remark now that $T$ has a dense range. Indeed if we suppose on the contrary that $\overline{\Ran(T)}\neq H$, then there exists a vector $x\in \Ran(T)^\perp$ such that $\|x\|=1$, and this implies that $0=\ps{Tx,x}\in W(T)$. This gives a contradiction with the hypothesis that $0\notin W(T)$.
\par\smallskip
It now suffices to apply \Cref{Prop1} to deduce that $T$ is embeddable in a $C_0$-semigroup, and this concludes the proof of \Cref{Prop1}.
\end{proof}

As mentioned before, the density of the range of the operator 
is automatic for sectorial Toeplitz operators. So the last part of the proof of \Cref{Theo3} is not necessary in this context. In other words, exactly the same proof yields the following result, which was already stated in the \Cref{Section 1} as \Cref{Theo4}:

\ExempleImageNumerique*

The numerical range of Toeplitz operators is well known: its closure is the closed convex hull of the spectrum; see \cites{BrownHalmos1963,Klein1972} for instance.  Thus we have:

\begin{corollary}
Let $F\in L^\infty(\mathbb T)$. Suppose that there exists $a\in\mathbb C\setminus \{0\}$ such that for almost every $\tau \in\mathbb T$, $\Re(aF(\tau))\ge0$. Then $T_F$ is embeddable into a $C_0$-semigroup of operators on $H^2$.
\end{corollary}

\begin{proof}
Without loss of generality, we can suppose  that $a=1$. As mentioned above, $\overline{W(T_F)}=\overline{\conv}(\sigma(T_F))$, and thus $W(T_F)$ is contained in the half plane $\overline{S_{\pi/2}}$. Alternatively, we can observe that for every $f\in H^2$, we have
\begin{align*}
    \Re\ps{T_Ff,f}=\Re\ps{Ff,f}&=\Re\left(\int_0^{2\pi}F(e^{i\theta})|f(e^{i\theta})|^2\,\mathrm d\theta\right)\\&=\int_0^{2\pi}\Re(F(e^{i\theta}))|f(e^{i\theta})|^2\,\mathrm d\theta\ge0.
\end{align*}
So $0$ does not belong to the interior of ${W(T_F)}$, and it follows from \Cref{Theo4} that $T_F$ is embeddable.
\end{proof}

\subsection{Link with the Kreiss constant, and circularly convex domains}
Let $\Omega$ be a subset of $\mathbb C$, with $\overline{\Omega}\neq\mathbb C$, and let $T\in B(X)$. The Kreiss constant of $T$ with respect to the subset $\Omega$ is defined as
\[\mathcal K_T(\Omega)~=~\sup_{z\notin\overline{\Omega}}\dist(z,\Omega)\|(z-T)^{-1}\|,\]
where we make the convention that $\|(z-T)^{-1}\|=\infty$ if $z\in\sigma(T_F)$.
This  Kreiss constant satisfies the following properties (see \cite{TohTrefethen1999} or \cite{Ostermann2021} for details):

\begin{proposition}\label{Prop-Kreiss-Constant}
Let $T$ be a bounded operator on a Banach space $X$.
\begin{enumerate}
    \item  Let $\Omega\subseteq \mathbb C$. If $\mathcal K_T(\Omega)$ is finite then $\sigma(T)\subset\overline{\Omega}$.\par\smallskip
    \item Let $\Omega_1,\Omega_2\subseteq \mathbb C$. If $\Omega_1\subseteq\Omega_2$, then $\mathcal K_T(\Omega_2)\le \mathcal K_T(\Omega_1)$.\par\smallskip
    \item If $X$ is a Hilbert space, then $\mathcal K_T(W(T))=1$.
\end{enumerate}
\end{proposition}

Here is a standard fact which will be used in the sequel.

\begin{proposition}\label{PropAj2}
Let $\Omega$ be a subset of $\mathbb C$ with $\overline{\Omega}\neq\mathbb C$, and let $T$ be a bounded operator on a Banach space $X$. The Kreiss constant $\mathcal K_T(\Omega)$ is finite if and only if the following two properties hold: $\sigma(T)\subseteq\overline\Omega$, and there exists an open neighborhood $U$ of $\overline\Omega$ and a constant $C\ge0$ such that 
\begin{equation}\label{Eqn1}
\|(z-T)^{-1}\|\le\frac{C}{\dist(z,\Omega)}~\text{ for every}~z\in U\setminus\overline\Omega.
\end{equation}
\end{proposition}

\begin{proof}
The direct implication is clear. Suppose conversely that for some constant $C\ge0$  and some open neighborhood $U$ of $\overline\Omega$, the inequality (\ref{Eqn1}) holds.
\par\smallskip
We decompose $\mathbb C\setminus\overline{\Omega}$ as follows. Let $R>\|T\|$. Since $U$ is an open neighborhood of $\overline\Omega$, we have
\begin{equation}\label{Eq:decomp39}
\mathbb C\setminus\overline{\Omega}~=\big(\mathbb C\setminus U\big)\cup\big(U\setminus\overline{\Omega)}~=~\big(R\overline{\mathbb D}\setminus U)\cup\big(\mathbb C\setminus(R\overline{\mathbb D}\cup U)\big)\cup \big(U\setminus\overline\Omega\big).
\end{equation}
Thus in order to prove that $\mathcal K_T(\Omega)<\infty$, we need to prove an analogue of the inequality \eqref{Eqn1}, where $U\setminus \overline\Omega$ is replaced first by $\mathbb C\setminus(R\overline{\mathbb D}\cup U)$, and then by $R\overline{\mathbb D}\setminus U$.
\par\smallskip
Let $z\in\mathbb C\setminus(R\overline{\mathbb D}\cup U)$. Note that $|z|>R>\|T\|$ so $z-T$ is invertible and \[\|(z-T)^{-1}\|\le\frac1{|z|-\|T\|}\cdot\]
Moreover, since $\sigma(T)\subset\overline{\Omega}$ and $R>\|T\|\ge\rho(T)$, we get 
\[\dist(z,\Omega)\le\dist(z,\sigma(T))\le|z|+\rho(T)\le |z|+R.\]
Thus, for every $z\in\mathbb C\setminus(R\overline{\mathbb D}\cup U)$,
we obtain
\begin{equation}\label{eq:dec49}
    \dist(z,\Omega)\|(z-T)^{-1}\|\le\frac{|z|+R}{|z|-\|T\|}=1+\frac{R+\|T\|}{|z|-\|T\|}\le \frac{2R}{R-\|T\|}\cdot
\end{equation}

Now let $K=R\overline{\mathbb D}\setminus U$ and suppose that $K\neq\varnothing$. The set $K$ is a compact subset of $\mathbb C\setminus\overline\Omega\subset\mathbb C\setminus\sigma(T)$ and since the map $z\mapsto (z-T)^{-1}$ is analytic on $\mathbb C\setminus\sigma(T)$, the application $\phi:z\longmapsto\dist(z,\Omega)\|(z-T)^{-1}\|$ is continuous on $K$ and thus $\phi$ is bounded on $K$. Let $C'=\max_K\phi$. Then we have
\begin{equation}\label{eq:dec79}
   \|(z-T)^{-1}\|\le\frac{C'}{ \dist(z,\Omega)}~\text{for every}~z\in K=R\overline{\mathbb D}\setminus U.
\end{equation}
If $K=\varnothing$, take any constant $C'\ge0$ in the rest of the proof.
Now let
\[\mathcal K=\max\left(C,\,\frac{2R}{R-\|T\|},\,C'\right).\]
Using the decomposition \eqref{Eq:decomp39} and the inequalities \eqref{Eqn1}, \eqref{eq:dec49} and \eqref{eq:dec79}, it follows that
\[\|(z-T)^{-1}\|\le\frac{\mathcal K}{ \dist(z,\Omega)}~\text{for every}~z\in \mathbb C\setminus\overline{\Omega},\]
i.e. the Kreiss constant of $T$ with respect to $\Omega$ is finite and satisfies $\mathcal K_T(\Omega)\le \mathcal K$.
\end{proof}

We will say that a subset $D$ of $\mathbb{C}$ is 
{\it a Riemann sphere disk} if it is either an open disk or the exterior of a closed disk. The next lemma will be useful in the proof of \Cref{Coro12}.

\begin{lemma}\label{PropAj1}
Let $\phi(z)=(az+b)(cz+d)^{-1}$, $z\in\mathbb{C}\setminus \{-d/c\}$, be a M\"obius transformation, and let $T$ be a bounded operator on a Banach space $X$. Let also $D$ be a Riemann sphere disk. Suppose that $\phi(T)$ is a well-defined bounded operator on $X $(that is, 
$cT+d$ is invertible), and that $\phi(D)$ is not a half-plane (and therefore is a Riemann sphere disk). 
If $\mathcal K_T(D)$ is finite, then 
$\mathcal K_{\phi(T)}(\phi(D))$ is finite as well. 
\end{lemma}
This lemma is a particular case of Lemma 2.1 in \cite{Bakaev1998}. We include here a simpler proof of this particular case.
\begin{proof}
Since $\phi(D)$ is not a half-plane, $-d/c\notin\partial D$. So let $U$ be a bounded open set satisfying $\partial D\subset U\subset\overline{U}\subset\mathbb C\setminus\{-d/c\}$, and set $V=\phi(U)$. Then $\phi$ is a conformal mapping from a neighborhood of  $\overline U$ onto a neighborhood of  $\overline V$. In particular there exists a constant $C_1>0$ such that for all $z,w\in U$, we have $|\phi(z)-\phi(w)|\le C_1|z-w|$. Observe that for $z\in\mathbb C\setminus \overline{D}$, since $D$ is an open disk or the exterior of a closed disk, we have $\dist(z,D)=\dist(z,\partial D)$. Hence
\begin{equation}
	\label{eqn:Dm1}
\frac{1}{\dist(z,D)} \le C_1\,\frac1{\dist(\phi(z),\phi(D))}\;\textrm{ for every }
z\in U\setminus \overline D. 
\end{equation}
Let $z\in U\setminus \overline D$. Then 
\begin{align*}
\phi(z)-\phi(T)
&=\frac{az+b}{cz+d}-(aT+b)(cT+d)^{-1}\\
&=\Big[(az+b)(cT+d)-(cz+d)(aT+b)\Big](cz+d)^{-1}(cT+d)^{-1}\\
&=(ad-bc)(z-T)(cz+d)^{-1}(cT+d)^{-1}.
\end{align*}
Since $\mathcal K_T(D)$ is finite, we know by assertion (1) of \Cref{Prop-Kreiss-Constant} that $\sigma(T)\subset\overline{D}$. In particular, for every $z\in U\setminus \overline{D}$, the operator $(z-T)$ is invertible and we have   
\begin{equation}
	\label{eqn:Dm2}
\|(\phi(z)-\phi(T))^{-1}\|
\le 
\frac{\sup_{z\in U\setminus\overline{D}}|cz+d|}{|ad-bc|} \|cT+d\| \|(z-T)^{-1}\|.   
\end{equation}
This implies that there exists a constant $C_2>0$ such that for every $z\in V\setminus\overline{\phi(D)}$, we have
\[\|(z-\phi(T))^{-1}\|\le C_2\|(w-T)^{-1}\|\le C_2\frac{\mathcal K_T(D)}{\dist(w,D)}\le C_1C_2\frac{\mathcal K_T(D)}{\dist(z,\phi(D))}\]
(it suffices to apply (\ref{eqn:Dm2}) above to $w=\phi^{-1}(z)$).
Since $\sigma(T)\subseteq\overline{D}$ we have the inclusion $\sigma(\phi(T))\subseteq\phi(\overline{D})$. Since $\partial \phi (D)\subset V$, it follows that $V\cup\overline{\phi(D)}=V\cup \phi(D)$ is an open neighborhood of $\overline{\phi(D)}$ and thus, by \Cref{PropAj2}, we finally conclude that $\mathcal K_{\phi(T)}(\phi(D))$ is finite.
\end{proof}

Here is now an important characterization of sectoriality.

\begin{theorem}\label{Thm-sectorial-equivalent-Kreiss-fini}
Let $T$ be a bounded operator on a Banach space $X$. The following assertions are equivalent:
\begin{enumerate}
    \item there exists $\omega\in[0,\pi)$ such that  $T$ is sectorial of angle $\omega$;\par\smallskip
    \item there exists $\omega\in[0,\pi)$ such that $\mathcal K_T(S_\omega)$ is finite.
\end{enumerate}
\end{theorem}
\begin{proof}
$(1)\implies (2)$: Assume first that $T$ is sectorial of angle $\omega$ and let $\phi\in(\omega,\pi)$. Since $0\in\overline {S_\phi}$, we have  $|z|\ge \dist(z,S_\phi)$ for all $z\in\mathbb C\setminus\overline {S_\phi}$. Thus
    \[\mathcal K_T(S_\phi)=\sup_{z\notin\overline{S_\phi}}\dist(z,S_\phi)\|(z-T)^{-1}\|\le \sup_{z\notin \overline{S_\phi}}\|z(z-T)^{-1}\|<\infty,\]
which gives (2).
\par\smallskip
$(2)\implies (1)$: Suppose now that there exists $\omega\in[0,\pi)$ such that $\mathcal K_T(S_\omega)$ is finite. Then by assertion (1) of  \Cref{Prop-Kreiss-Constant}, we have $\sigma(T)\subseteq \overline{S_\omega}$.

Moreover, note that for all $z=re^{i\theta}$ with $r>0$ and $\omega<|\theta|\le\pi$, we have that
\[\dist(z,S_\omega)=\begin{cases}
r\sin(|\theta|-\omega)&\text{if }|\theta|\le \omega+\pi/2\\
r&\text{else.}
\end{cases}
\]
So, set $C_\phi:=\sin(\phi-\omega)$ if $\omega<\phi\le\omega+\pi/2$ and $C_\phi:=1$ if $\phi\ge \omega+\pi/2$. Then for all $\phi\in(\omega,\pi]$ and all $z\in \mathbb C\setminus \overline{S_\phi}$, we have $\dist(z, S_\omega)\ge C_\phi|z|$. Using that 
$\mathbb C\setminus\overline{S_\phi}\subseteq \mathbb C\setminus\overline{S_\omega}$, this implies that
\[\sup_{z\notin\overline{S_\phi}}\|z(z-T)^{-1}\|\le \frac 1{C_\phi}\sup_{z\notin\overline{S_\phi}}\dist(z,S_\omega)\|(z-T)\|^{-1}\le \frac{\mathcal K_T(S_\omega)}{C_\phi}<\infty.\]
Thus we deduce that $T$ is sectorial of angle $\omega$, which gives (1).
\end{proof}

\begin{remark}\label{rem-sup-inf-sectorial-Kreiss}
Note that the proof of \Cref{Thm-sectorial-equivalent-Kreiss-fini} implies that if $\mathcal K_T(S_\omega)<\infty$ for some $\omega\in [0,\pi)$, then $T$ is sectorial of angle $\omega$.
\end{remark}

Combining Theorem \ref{Thm-sectorial-equivalent-Kreiss-fini} and Theorem \ref{Theo2}, we immediately obtain the following corollary.

\begin{corollary}\label{cor:Kreiss-Toeplitz}
Let $F\in L^\infty(\mathbb T)$ and $p>1$. Suppose that there exists $\omega\in[0,\pi)$ such that $\mathcal K_{T_F}(S_\omega)$ is finite, where $T_F$ is viewed as a bounded operator on $H^p$. Then $T_F$ is embeddable into a $C_0$-semigroup on $H^p$.
\end{corollary}

When $p=2$, we have $\mathcal K_{T_F}(W(T_F))=1$ by assertion (3) of \Cref{Prop-Kreiss-Constant}. Thus \Cref{cor:Kreiss-Toeplitz} is a generalization of \Cref{Theo4}.
\par\medskip
The estimation of the resolvent norm of a Toeplitz operator in terms of the distance to the spectrum is studied in the papers \cites{MR4849682,Peller1986}. In particular, conditions on the symbols are given ensuring that the Kreiss constant of the spectrum is finite. In \cite{MR4849682}, the authors study the case of Laurent polynomials, and in \cite{Peller1986}, Peller studies the case of 
 sufficiently regular symbols $F$ such that $\sigma(T_F)$ is a so-called \emph{circularly convex set} (or, more generally, the case where $\sigma(T_F)$ is contained in a circularly convex set). We now present some consequences of Peller's results concerning the embedding problem, in the case where $p=2$.
\par\smallskip
We say that a compact subset $\Omega$ of $\mathbb{C}$ is \emph{circularly convex} if there exists a radius $r>0$ such that for every $\lambda\in\mathbb C\setminus \Omega$ with $\dist(\lambda, \Omega)<r$, there exist two points $\mu\in\partial \Omega$ and $\nu\in\mathbb C\setminus \Omega$ for which $|\mu-\nu|=r, \lambda\in(\mu,\nu)$ and such that $\{\zeta\in\mathbb C\,;\,|\nu-\zeta|<r\}\cap \Omega=\varnothing.$ In other words, one can roll a disk of radius $r$ along the boundary of $\Omega$ while remaining in the complement of $\Omega$. Note that if $\Omega$ is a convex set, or if its boundary is $C^2$-smooth, then $\Omega$ is circularly convex.
\par\smallskip
The result quoted below is not formally present in \cite{Peller1986}, and it is stated there only in the case where $\Omega=\sigma(T_F)$, but by repeating literally the arguments used in the proof of \cite{Peller1986}*{Th. 4}, we obtain \Cref{Theo9}. 

\begin{theorem}[\cite{Peller1986}*{Th. 4}]\label{Theo9}
Let $\cal X$ be a Banach algebra of functions on $\mathbb T$ satisfying the following properties:
\begin{enumerate}[i)]
    \item $\cal X$ is continuously embedded in $C(\mathbb T)$;\par\smallskip
    \item $P_+(\mathcal X)\subseteq L^\infty(\mathbb T)$;\par\smallskip
    \item Every multiplicative linear function on $\cal X$ coincides with the function evaluation at some point $\zeta\in\mathbb T$, i.e. $f\mapsto f(\zeta)$.
\end{enumerate}
Let $F\in \cal X$. 
Suppose that $\Omega$ is a circularly convex compact set which contains $\sigma(T_F)$. Then $\mathcal K_{T_F}(\Omega)<\infty$, where the Toeplitz operator $T_F$ is viewed as an operator on $H^2$.
\end{theorem}

Peller gave in \cite{Peller1986} some examples of Banach algebras satisfying the hypothesis of \Cref{Theo9}.  These conditions are satisfied when for example $\cal X$ is the Wiener algebra, or when $\cal X$ is the space of Dini-continuous functions, i.e. of functions $f\in C(\mathbb T)$ such that
\[\int_0^1\frac{w_f(t)}t\mathrm dt<\infty\]
where $w_f(t)=\sup\{|f(z)-f(z')|, |z-z'|\le t\}$ is the modulus of continuity of $f$. 
\par\smallskip
A first consequence of \Cref{Theo9} is:

\begin{corollary}\label{cor:contained-in-sector}
Let $\cal X$ be a Banach algebra of functions on $\mathbb{T}$ satisfying the assumptions of \Cref{Theo9}.
Let $F\in \cal X$ be such that $\sigma(T_F)$ is circularly convex, and contained in a closed sector $\overline{S_\omega}$ for some $\omega\in[0,\pi)$. Then $T_F$ is embeddable into a $C_0$-semigroup on $H^2$.
\end{corollary}
\begin{proof}
It follows from \Cref{Theo9} that $\mathcal K_{T_F}(\sigma(T_F))<\infty$. Then, according to assertion (2) of \Cref{Prop-Kreiss-Constant}, $\mathcal K_{T_F}(S_\omega)=K_{T_F}(\overline{S_\omega})<\infty$, and it remains to apply \Cref{cor:Kreiss-Toeplitz} to conclude that $T_F$ is embeddable into a $C_0$-semigroup on $H^2$.
\end{proof}

Another consequence of \Cref{Theo9} is the following result, which is more general than ~\Cref{cor:contained-in-sector} 
and the statement of \Cref{Theo14} in \Cref{Section 1}.

\begin{theorem}\label{Coro12}
Let $\mathcal X$ be a Banach algebra of functions on $\mathbb{T}$ satisfying the assumptions of  \Cref{Theo9}.
Let $F \in \mathcal X$. Suppose that  
there exists an open disk $D$, contained in the 
unbounded component of $\mathbb C\setminus \sigma(T_F)$, such that $0\in \partial D$. 
Then $T_F$ is embeddable into a $C_0$-semigroup of operators on $H^2$.
\end{theorem}

\begin{proof}Let $\Delta=\mathbb C\setminus \overline{D}$, 
and $\Delta'=\overline{\Delta}\cap R\overline{\mathbb D}=R\overline{\mathbb D}\setminus D$ where $R$ is any positive radius with $R>\|T_F\|$. Then $\Delta'$ is circularly convex, and contains the spectrum of $T_F$. By \Cref{Theo9}, the Kreiss constant $\mathcal K_{T_F}(\Delta')$ is finite, hence $\mathcal K_{T_F}(\Delta)$ is also finite by assertion (2) of \Cref{Prop-Kreiss-Constant}. Let $a\in D$
and set $\phi(z)=\frac z{z-a}$, $z\in\mathbb{C}\setminus\{a\}$. Then the operator $T-a$ is invertible. Since $\Delta$ is a Riemann sphere disk and $\mathcal K_{T_F}(\Delta)<+\infty$, by \Cref{PropAj1}, the operator $S:=T_F(T_F-a)^{-1}=\phi(T_F)$ is such that $\mathcal K_S(\phi(\Delta))$ is also finite. Since $\phi(\partial D)$ is a circle, 
and since $\phi(a)=\infty$ and $\phi(0)=0$, we see that $\phi(\Delta)$ is a disk, and $0$ belongs to $\partial\phi(\Delta)$. So there exists $\alpha\in\mathbb R$ such that $e^{i\alpha}\phi(\Delta)\subseteq\{z\in\mathbb C\,;\,\Re(z)\ge0\}=\overline{S_{\pi/2}}$. But $\mathcal K_{e^{i\alpha}S}(e^{i\alpha}\phi(\Delta))=\mathcal K_{S}(\phi(\Delta))<\infty$, and thus $\mathcal K_{e^{i\alpha}S}(S_{\pi/2})$ is finite. By \Cref{Thm-sectorial-equivalent-Kreiss-fini} (or \Cref{rem-sup-inf-sectorial-Kreiss}), it follows that $e^{i\alpha}S$ is sectorial of angle $\pi/2$. 
Since $S=T_F(T_F-a)^{-1}=(T_F-a)^{-1}T_F$ and $(T_F-a)^{-1}$ is invertible, we have that $\ker(S)=\ker(T_F)$ and $\Ran(S)=\Ran(T_F)$. So again, the combination of Coburn's Lemma and the decomposition given by \eqref{Eq:decompEspRefl} implies that $S$ has dense range.
So, by \Cref{Prop1}, $S$ is embeddable into a $C_0$-semigroup $(A_t)_{t>0}$ of operators on $H^2$ which is given by
\[A_t=e^{-i\alpha t}\int_{\Gamma}z^t(z-e^{i\alpha}S)^{-1}\,\frac{\mathrm dz}{2i\pi}=e^{-i\alpha t}\int_{\Gamma}z^t(z-e^{i\alpha}T_F(T_F-\alpha)^{-1})^{-1}\,\frac{\mathrm dz}{2i\pi},\] where $\Gamma=\partial (S_{3\pi/4}\cap R\mathbb D)$ for some $R>\|S\|$.  Note that if an operator $A$ commutes with $T_F$, it will commute with $R(T_F)$ for every rational function $R$ without poles in $\sigma(T_F)$ and thus $A$ will commute also with any operator of the semigroup $(A_t)_{t>0}$.
\par\smallskip
Now, since $a$ belongs to the unbounded component of $\mathbb C\setminus \sigma(T_F)$, it follows that $0$ belongs to the unbounded component of $\mathbb C \setminus \sigma(T_F-a)$ and thus there exists an analytic determination of the logarithm on a neighborhood of $\sigma(T_F-a)$, denoted by $\log$. Let $(B_t)_{t>0}$ be the $C_0$-semigroup of operators defined by $B_t=e^{t\log(T_F-a)}$, $t>0$, and set $S_t:=A_tB_t$. Since $B_t$ commutes with $T_F$, it also commutes with $A_t$ by the observation above, 
thus $(S_t)_{t>0}$ is a semigroup of operators on $H^2$ which satisfies $S_1=A_1B_1=S(T_F-a)=T_F$. Finally, since $\log(T_F-a)\in\mathcal B(H^2)$, it follows that $(B_t)_{0<t<1}$ is uniformly bounded and so $(S_t)_{t>0}$ is a $C_0$-semigroup. 
\end{proof}

Note that the symbol $F$ considered in~\Cref{ex:root-3} is Dini-continuous, and the corresponding Toeplitz operator is not embeddable (for $p=2$). 
It follows that the above theorem will no longer be true if the hypothesis on the existence of 
an open disc $D$ as above is replaced by the condition that the unbounded component of $\C\sm F(\T)$ should contain an open sector with vertex at $0$.

\par\medskip
\noindent {\bf Acknowledgment:} This work was supported in part by the project COMOP of the French National Research Agency (grant ANR-24-CE40-0892-01). The authors acknowledge the support of the CDP C$^2$EMPI, as well as of the French State under the France-2030 program, the University of Lille, the Initiative of Excellence of the University of Lille, and the European Metropolis of Lille for their funding and support of the R-CDP-24-004-C2EMPI project. The third author also acknowledges the support of the CNRS. 
The fourth author was partially supported by Plan Nacional  I+D grant no. PID2022-137294NB-I00, Spain. Part of this work 
was done during his visit to Laboratoire Paul Painlevé, Univ. of Lille in 2025, and he expressed his gratitute for hospitality during this visit. 
\bibliographystyle{plain}
\bibliography{biblio}
\end{document}